\documentstyle[12pt]{article}
\input amssym.def
\begin {document}
\def \Z{{\bf Z}}
\def \C{{\bf C}}

\def \Q{\Bbb Q}
\def \N{{\bf N}}

\def \bZ{{\bf Z}}

\def \wt{{\rm wt}\;}
\def \fg{\frak g}

\def \Res{{\rm Res}}
\def \End{{\rm End}\;}
\def \Ker{{\rm Ker}\;}

\def \Aut{{\rm Aut}}
\def \Hom{{\rm Hom}}
\def \mod{{\rm mod}\;}

\def \<{\langle} 
\def \>{\rangle} 

\def \a{\alpha }
\def \e{\epsilon }

\def \b{\beta }

\def \be{\begin{equation}\label}
\def \ee{\end{equation}}
\def \bl{\begin{lem}\label}
\def \el{\end{lem}}
\def \bt{\begin{thm}\label}
\def \et{\end{thm}}
\def \bp{\begin{prop}\label}
\def \ep{\end{prop}}
\def \br{\begin{rem}\label}
\def \er{\end{rem}}
\def \bc{\begin{coro}\label}
\def \ec{\end{coro}}
\def \bd{\begin{de}\label}
\def \ed{\end{de}}
\def \pf{{\bf Proof. }}

\newtheorem{thm}{Theorem}[section]
\newtheorem{prop}[thm]{Proposition}
\newtheorem{coro}[thm]{Corollary}

\newtheorem{lem}[thm]{Lemma}
\newtheorem{rem}[thm]{Remark}
\newtheorem{de}[thm]{Definition}

\makeatletter
\@addtoreset{equation}{section}

\makeatother
\makeatletter

\begin{center}{\Large \bf On abelian coset 
generalized vertex algebras}

\vspace{0.5cm}
Haisheng Li\footnote{Partially supported by NSF grant
DMS-9970496 and by a grant from the Rutgers University Research Council}\\ 
Department of Mathematical Sciences, Rutgers University,
Camden, NJ 08102
\end{center}

{\bf Abstract} 
This paper studies the algebraic aspect of a general abelian 
coset theory with a work of Dong and Lepowsky as our main motivation.
It is proved that the vacuum space $\Omega_{V}$
(or the space of highest weight vectors) of a Heisenberg algebra  
in a general vertex operator algebra $V$ has a natural 
generalized vertex algebra structure in the sense of Dong and Lepowsky
 and that 
the vacuum space $\Omega_{W}$ of a $V$-module $W$ is a natural
$\Omega_{V}$-module. The automorphism group 
$\Aut_{\Omega_{V}}\Omega_{V}$ of the adjoint $\Omega_{V}$-module
is studied and it is proved to be a central extension 
of a certain torsion free abelian group by $\C^{\times}$.
For certain subgroups $A$ of $\Aut_{\Omega_{V}}\Omega_{V}$,
certain quotient algebras $\Omega_{V}^{A}$ of $\Omega_{V}$ are 
constructed.
Furthermore, certain functors among the category of $V$-modules, 
the category of 
$\Omega_{V}$-modules and the category of $\Omega_{V}^{A}$-modules
are constructed and irreducible $\Omega_{V}$-modules and 
$\Omega_{V}^{A}$-modules are classified in terms of irreducible 
$V$-modules. If the category of $V$-modules is semisimple,
then it is proved that the category of $\Omega_{V}^{A}$-modules 
is semisimple.

\section{Introduction}
Vertex operator algebras are a family of new ``algebras,'' which arose
in mathematics in the 80's.
Affine Lie algebras are attached in a certain way
to the family of vertex operator algebras.
More specifically, for an affine Lie algebra $\hat{\fg}$ 
with a fixed level $\ell$, a certain generalized Verma module $M(\ell,0)$
has a canonical vertex operator algebra structure and
the category of $M(\ell,0)$-modules is canonically isomorphic to 
the category of suitably restricted $\hat{\fg}$-modules of level $\ell$
(see for example [DL2], [FF], [FZ], [Li2], [MP2]).
In the study of the representation theory for affine Lie algebras,
$Z$-algebras, introduced in [LW4-7] and [LP1-3], have played an important role
(see for example [LW4-7], [LP1-3], [C], [Hu], [Ma], [Mi1-4], 
[MP1-2] and [P]). 
In fact, the construction of modules and the understanding of the
new algebraic structures underlying the construction using
$Z$-algebras amounts to
the vertex-operator-theoretic interpretation of 
generalized Rogers-Ramanujan partition identities [LW4-7].
%The $Z$-operators (generators of $Z$-algebras) 
%were constructed in terms of ordinary vertex operators
%and they centralize an underlying Heisenberg algebra.
To place the untwisted $Z$-algebras [LP1-3]
into an elegant axiomatic context, 
and in particular, to embed such algebras into larger, 
more natural algebras, Dong and Lepowsky in [DL2]
introduced ``relative (untwisted) vertex operators'' 
and ``quotient relative vertex operators''
and studied the algebraic structure of such vertex operators.
This study led in [DL2] to three levels of generalization of the
concept of vertex operator algebra (and module). 
One of these notions --- that of ``generalized vertex operator algebra''
--- enabled 
Dong and Lepowsky to clarify the essential equivalence between
the $Z$-algebras of [LP1-2] and the parafermion algebra of [ZF1],
among other things, providing a precise mathematical foundation for
``parafermion conformal field theory''.
In this sense, $Z$-algebras  or parafermion algebras 
are attached to the family of generalized 
vertex operator algebras.
Algebraic structures similar to generalized vertex operator algebras, 
and to generalized vertex algebras have been independently
introduced and studied in [FFR], with certain different 
motivations and examples (involving spinor constructions), 
and also in [Mo].

The theory of relative (untwisted) vertex operators in [DL2] was built on
a space $V_{L}$ associated to a rational lattice $L$.
For a certain abelian subgroup $A$
of a central extension of $L$, the space $\Omega_{V_{L}}^{A}$ and 
quotient relative vertex operators on $\Omega_{V_{L}}^{A}$
were explicitly constructed.
Note that on the one hand, when $A=0$, $\Omega_{V_{L}}^{A}$ specializes to 
$\Omega_{V_{L}}$, the vacuum space of a Heisenberg algebra in $V_{L}$,
and on the other hand, 
$\Omega_{V_{L}}^{A}$ is a quotient space of $\Omega_{V_{L}}$.
Among the results in [DL2], it was proved that 
$\Omega_{V_{L}}^{A}$ is a 
generalized vertex algebra for an even lattice $L$.
Furthermore, for a simple Lie algebra $\fg$ of type $A,D,E$ 
and for a positive integer level $\ell$,
by using the technique embedding the higher-level standard modules 
into tensor products of level-one standard modules,
it was proved that for a certain $A$, $\Omega_{L(\ell,0)}^{A}$
defined similarly is 
a generalized vertex operator algebra and $\Omega_{L(\ell,\lambda)}^{A}$ 
is an irreducible $\Omega_{L(\ell,0)}^{A}$-module for any standard 
$\hat{\fg}$-module $L(\ell,\lambda)$ of level $\ell$.

In [GL], as an application of a general result on
generalized vertex algebras generated by parafermion-like 
vertex operators, it was proved that
the vacuum space $\Omega_{L(\ell,0)}$ of the homogeneous
Heisenberg algebra of $\hat{\fg}$ has a canonical
generalized vertex algebra structure for any nonzero level $\ell$
(not just positive integer levels).

As discussed in [DL2], the construction of the Virasoro
algebra of generalized vertex operator algebras $\Omega_{L(\ell,0)}^{A}$
is a special case of ``coset construction'' [GKO1-2] (cf. [FZ]).
Coset constructions and dual pairs for a general pair 
of vertex operator algebras have been studied in [DM]
and a beautiful duality of Schur-Weyl type was obtained.
The study of [DM] is in the context of vertex operator algebras 
and modules, i.e., the vacuum space was studied as a module for
the commutant vertex operator subalgebra of the Heisenberg algebra.

In this paper, we revisit the abelian coset theory from 
a different point of view with a different approach and
our naive purpose is to give a solution to 
the following problems that naturally arise from [DL2]:

a) Explicitly construct the generalized vertex algebra
$\Omega_{L(\ell,0)}$ for {\em any} nonzero level $\ell$
by using the vertex operator algebra structure on $L(\ell,0)$. 

b) Classify all irreducible $\Omega_{L(\ell,0)}^{A}$-modules
for a positive integer level $\ell$.

c) Prove the conjecture that for a positive integer level $\ell$,
every $\Omega_{L(\ell,0)}^{A}$-module is completely reducible.
(It was known (cf. [FZ], [Li2]) that the category of $L(\ell,0)$-modules 
is semisimple.)

d) Prove the conjecture that $\Omega_{L(\ell,0)}^{A}$ is 
the vacuum space of a vertex operator subalgebra in $L(\ell,0)$.
(Note that $\Omega_{L(\ell,0)}$ by definition is the vacuum space of
the vertex operator subalgebra generated by the Heisenberg algebra.)

These are rather basic and natural questions in vertex 
operator algebra theory and a solution to these questions, 
in particular, b) and c), is also important in the study of
abelian coset conformal field theory. In this paper,
not only do we completely solve 
all the four problems, but this is done with great generality 
with a general vertex operator algebra $V$
in place of $V_{L}$ and $L(\ell,0)$ also. 
Consequently, this provides a mathematical foundation for
general abelian coset conformal field theory.

Now, we explain our main results.
Let $V$ be a vertex operator algebra with  a subspace ${\bf h}$ of 
$V_{(1)}$ such that components $h(m)$ of $Y(h,z)$ for $h\in {\bf h}$
satisfy the Heisenberg algebra relation with a level $\ell$
and such that $h(0)$ acts semisimply on $V$.
Let $\Omega_{V}$
be the vacuum space of the Heisenberg algebra $\hat{\bf h}$ in $V$.
We define a vertex operator map $Y_{\Omega}$ on $\Omega_{V}$ by 
using the construction of $Z$-operators ([LW4-7], [LP1-3]) and
we prove that $\Omega_{V}$ equipped with $Y_{\Omega}$ is 
a generalized vertex algebra in the sense of [DL2].
We also prove that for any $V$-module $W$, 
$\Omega_{W}$ is a natural $\Omega_{V}$-module.
We define a notion of $\Omega_{V}$-${\bf h}$-module,
analogous to a notion in [LW4] and [LP2], and then we prove that
the category of weak $V$-modules of a certain type is naturally 
equivalent to the category of $\Omega_{V}$-${\bf h}$-modules.
As the first step in classifying irreducible $\Omega_{V}$-modules, 
we prove that for two irreducible $V$-modules $W_{1},W_{2}$, 
$\Omega_{W_{1}}$ and $\Omega_{W_{2}}$ are isomorphic $\Omega_{V}$-modules
if and only if $W_{2}$ is isomorphic to $W_{1}^{(\a)}$ for some 
$\a\in {\bf h}$, where $W_{1}^{(\a)}$ is a $V$-module constructed
in [Li3] (see also [Li4-6]) as a deformation of $W_{1}$.
When $\ell$ is irrational, it is proved that the 
generalized vertex algebra $\Omega_{V}$ is simple and 
all irreducible $\Omega_{V}$-modules are classified
as those $\Omega_{W}$ with $W$ being irreducible $V$-modules.

In general, generalized vertex algebra $\Omega_{V}$ may have 
nontrivial ideals.
To explicitly construct quotient algebras, motivated by
our previously mentioned results on $\Omega_{V}$ and by [DL2] 
we study the automorphism group $\Aut_{\Omega_{V}}\Omega_{V}$ 
of the adjoint $\Omega_{V}$-module.
Assuming that $V$ is simple, we show that the group 
$\Aut_{\Omega_{V}}\Omega_{V}$ 
is a central extension of a subgroup $K$ 
in ${\bf h}$ by $\C^{\times}$, where $K$ consists of those 
$\a\in {\bf h}$ such that $V^{(\a)}\simeq V$ as a $V$-module.
We then construct quotient generalized vertex algebras 
$\Omega_{V}^{A}$ of $\Omega_{V}$ for certain central 
subgroups $A$ of $\Aut_{\Omega_{V}}\Omega_{V}$. 
We also construct certain functors between the categories
of $\Omega_{V}$-${\bf h}$-modules and $\Omega_{V}^{A}$-modules
(with a certain grading group) are constructed.
With these functors,
all irreducible $\Omega_{V}^{A}$-modules are classified.
Furthermore, if the category of $V$-modules is semisimple, 
we show that the category of $\Omega_{V}^{A}$-modules is
semisimple. As the last result of our study on the general case,
we identify $\Omega_{V}^{A}$ as the vacuum 
space of a lattice vertex operator subalgebra of $V$.
In the last section, we apply our results 
for $V=L(\ell,0)$. In particular, we recover 
the corresponding results of [DL2].

General abelian coset theory has been also studied 
in physics such as in [Br1-3] and [Gep]. 
As usual, physics papers on 
abelian coset theory, or parafermion algebras 
focused on the conformal-field-theoretic 
structure, including the role of the Virasoro algebra.
On the other hand, the present paper like [DL1-3], [LW4-7] 
and [LP1-3] focuses on the vertex-operator-algebraic structure.

A natural continuation of the present study is to 
study the twisted case as a counterpart of [DL3].
This will be done in a future publication.

This paper is organized as follows: Section 2 is preliminary;
We recall certain basic definitions and results from [DL2].
In Section 3, we prove that the vacuum space $\Omega_{V}$ 
is a generalized vertex algebra with $\Omega_{W}$ as a module
for any $V$-module $W$. In Section 4, we determine the group
$\Aut_{\Omega_{V}}\Omega_{V}$.
In Section 5, we study quotient generalized 
vertex algebras $\Omega_{V}^{A}$ and classify their 
irreducible modules. In Section 6, as an application 
we study $\Omega_{L(\ell,0)}^{A}$ and its modules.

\section{Definitions of generalized vertex algebra and module}
This section is preliminary.
In this section, we recall from [DL2] the notions of 
generalized vertex (operator) algebra and module, and 
for our need we define and discuss
the notions of submodule, homomorphism and ideal.

First, we briefly review some formal variable notations.
Throughout this paper, $z, z_{0},z_{1},z_{2}$ and $x, y$ will be
mutually commuting (independent) formal variables. 
We shall use $\N$ for the nonnegative integers,
$\Z_{+}$ for the positive integers, ${\bf Q}$
for the rational numbers, $\C$ for the complex numbers
and $\C^{\times}$ for the nonzero complex numbers.

We shall use standard formal variable notations as defined 
in [FLM] and [FHL]. For example, for a vector space $U$, 
\begin{eqnarray}
U\{z\}=\{ \sum_{n\in \C}u(n)z^{n}\;|\; u(n)\in U
\;\;\;\mbox{ for }n\in \C\}.
\end{eqnarray}
The following are useful subspaces of $U\{z\}$:
\begin{eqnarray}
U[[z,z^{-1}]]&=&\{ \sum_{n\in \Z}u(n)z^{n}\;|\; u(n)\in U
\;\;\;\mbox{ for }n\in \Z\},\\
U((z))&=&\{ \sum_{n\in \Z}u(n)z^{n}\in U[[z,z^{-1}]]\;|\; u(n)=0
\;\;\;\mbox{ for }n\;\;\mbox{sufficiently small}\},\\
U[[z]]&=&\{ \sum_{n\in \Z}u(n)z^{n}\in U[[z,z^{-1}]]\;|\; u(n)=0
\;\;\;\mbox{ for }n<0\}.
\end{eqnarray}
A typical element of ${\bf C}[[z,z^{-1}]]$ is the formal Fourier expansion
of the delta-function at $0$:
\begin{eqnarray}
\delta(z)=\sum_{n\in \Z}z^{n}.
\end{eqnarray}
Its fundamental property is:
\begin{eqnarray}
f(z)\delta(z)=f(1)\delta(z)\;\;\;\mbox{ for }\; f(z)\in {\bf C}[z,z^{-1}].
\end{eqnarray}
For $\a\in \C$, by definition,
\begin{eqnarray}
(z_{1}-z_{2})^{\alpha}=\sum_{i\ge 0}{\a\choose i}(-1)^{i}z_{1}^{\a-i}z_{2}^{i}.
\end{eqnarray}
Then as a formal series,
\begin{eqnarray}
\delta\left(\frac{z_{1}-z_{2}}{z_{0}}\right)
=\sum_{n\in \Z}\left(\frac{z_{1}-z_{2}}{z_{0}}\right)^{n}
=\sum_{n\in \Z}\sum_{i\ge 0}{n\choose i}(-1)^{i}z_{0}^{-n}z_{1}^{n-i}z_{2}^{i}.
\end{eqnarray}
We have the following fundamental properties of delta function
(see [FLM], [FHL], [Le], [Zhu]):

\bl{ldeltajacobi}
For $\a\in \C$,
\begin{eqnarray}
z_{0}^{-1}\left(\frac{z_{1}-z_{2}}{z_{0}}\right)^{\a}
\delta\left(\frac{z_{1}-z_{2}}{z_{0}}\right)
=z_{1}^{-1}\left(\frac{z_{0}+z_{2}}{z_{1}}\right)^{-\a}
\delta\left(\frac{z_{0}+z_{2}}{z_{2}}\right);
\end{eqnarray}
For $r,s,k\in \Z$ and for $p(z_{1},z_{2})\in {\bf C}[[z_{1},z_{2}]]$,
\begin{eqnarray}
& &z_{0}^{-1}\delta\left(\frac{z_{1}-z_{2}}{z_{0}}\right)
z_{1}^{r}z_{2}^{s}(z_{1}-z_{2})^{k}p(z_{1},z_{2})
-z_{0}^{-1}\delta\left(\frac{z_{2}-z_{0}}{-z_{0}}\right)
z_{1}^{r}z_{2}^{s}(-z_{2}+z_{1})^{k}p(z_{1},z_{2})\nonumber\\
&=&z_{2}^{-1}\delta\left(\frac{z_{1}-z_{0}}{z_{2}}\right)
(z_{2}+z_{0})^{r}z_{2}^{s}z_{0}^{k}p(z_{2}+z_{0},z_{2}).
\end{eqnarray}
In particular,
\begin{eqnarray}
z_{0}^{-1}\delta\left(\frac{z_{1}-z_{2}}{z_{0}}\right)-
z_{0}^{-1}\delta\left(\frac{z_{2}-z_{0}}{-z_{0}}\right)
=z_{2}^{-1}\delta\left(\frac{z_{1}-z_{0}}{z_{2}}\right).
\end{eqnarray}
\el

A generalized vertex algebra as defined in [DL2] is associated to
an abelian group $G$, a symmetric 
$\C /2{\bf Z}$-valued ${\bf Z}$-bilinear form 
(not necessarily nondegenerate) on $G$:
\begin{eqnarray}
(g,h)\in {\bf C}/2{\bf Z}\;\;\;\mbox{ for }g,h\in G
\end{eqnarray}
and $c(\cdot,\cdot)$ is a ${\bf C}^{\times}$-valued
alternating ${\bf Z}$-bilinear form on $G$.

\bd{dgva}
{\em [DL2] {\em A generalized vertex algebra} associated
with the group $G$ and the forms $(\cdot,\cdot)$ and $c(\cdot,\cdot)$
is a vector space with two gradations:
\begin{eqnarray}
V=\coprod_{n\in \C}V_{n}=\coprod_{g\in G}V^{g};
\;\;\mbox{ for }v\in V_{n},\;\; n=\wt v;
\end{eqnarray}
such that
\begin{eqnarray}
V^{g}=\coprod_{n\in \C}V_{n}^{g}\;\;\;\;
\mbox{(where $V_{n}^{g}=V_{n}\cap V^{g}$)}\;\;\mbox{ for }g\in G,
\end{eqnarray}
equipped with a linear map
\begin{eqnarray}
Y: & &V \rightarrow (\End V)\{z\}\nonumber\\ 
& &v\mapsto Y(v,z)=\sum_{n\in \C}v_{n}z^{-n-1}\;\;\;\mbox{($v_{n}\in \End V$)}
\end{eqnarray}
and with two distinguished vectors 
${\bf 1}\in V_{0}^{0},\; \omega\in V_{2}^{0},$
satisfying the following conditions for $g,h\in G,\; u,v\in V$ and $l\in \C$:
\begin{eqnarray}
& &u_{l}V^{h}\subset V^{g+h}\;\;\;\mbox{ if }\;\; u\in V^{g};\\
& &u_{l}v=0\;\;\;\mbox{ if the real part of }\;l\;\;
\mbox{ is sufficiently large};\\
& &Y({\bf 1},z)=1;\\
& &Y(v,z){\bf 1}\in V[[z]]\;\;\mbox{ and }
\;\;\lim_{z\rightarrow 0}Y(v,z){\bf 1}=v;\\
& &Y(v,z)|_{V^{h}}=\sum_{n\equiv (g,h)\mod \bZ}v_{n}z^{-n-1}
\;\;\;\mbox{ if }v\in V^{g}
\end{eqnarray}
(i.e., $n+2\bZ\equiv (g,h)\;\mod \bZ/2\bZ$);
\begin{eqnarray}
& &z_{0}^{-1}\left(\frac{z_{1}-z_{2}}{z_{0}}\right)^{(g,h)}
\delta\left(\frac{z_{1}-z_{2}}{z_{0}}\right)Y(u,z_{1})Y(v,z_{2})\nonumber\\
& &-c(g,h)z_{0}^{-1}\left(\frac{z_{2}-z_{1}}{z_{0}}\right)^{(g,h)}
\delta\left(\frac{z_{2}-z_{1}}{-z_{0}}\right)Y(v,z_{2})Y(u,z_{1})\nonumber\\
&=&z_{2}^{-1}\delta\left(\frac{z_{1}-z_{0}}{z_{2}}\right)Y(Y(u,z_{0})v,z_{2})
\left(\frac{z_{1}-z_{0}}{z_{2}}\right)^{-g}
\end{eqnarray}
(the {\em generalized Jacobi identity}) if $u\in V^{g},\; v\in V^{h}$, where
\begin{eqnarray}
\delta\left(\frac{z_{1}-z_{0}}{z_{2}}\right)
\left(\frac{z_{1}-z_{0}}{z_{2}}\right)^{-g}\cdot w=
\left(\frac{z_{1}-z_{0}}{z_{2}}\right)^{-(g,k)}
\delta\left(\frac{z_{1}-z_{0}}{z_{2}}\right)w
\end{eqnarray}
for $w\in V^{k},\; k\in G$; 
\begin{eqnarray}\label{e2.25}
[L(m),L(n)]=(m-n)L(m+n)+{1\over 12}(m^{3}-m)\delta_{m+n,0}(\mbox{rank }V)
\end{eqnarray}
for $m,n\in \Z$, where
\begin{eqnarray}
L(n)=\omega_{n+1}\;\;\mbox{ for }n\in {\bf Z}, \;\mbox{ i.e., } 
Y(\omega,z)=\sum_{n\in {\bf Z}}L(n)z^{-n-2}
\end{eqnarray}
and
\begin{eqnarray}
& &\mbox{ rank }V\in \C;\\
& &L(0)v=nv=(\wt v)v\;\;\;\mbox{ for }n\in \C,\; v\in V_{n};\\
& &{d\over dz}Y(v,z)=Y(L(-1)v,z).\label{e2.29}
\end{eqnarray}}
\ed

This completes the definition. This generalized vertex algebra is denoted by 
\begin{eqnarray}
(V,Y, {\bf 1},\omega, G, (\cdot,\cdot), c(\cdot,\cdot)),
\end{eqnarray}
or simply by $V$.

Note: We here allow the form $(\cdot,\cdot)$ to be $\C/2\Z$-valued 
instead of $({1\over T}\Z)/2\Z$-valued,
so that the parameter $T$ in the original definition given in [DL2]
disappeared.

We recall the following remark from [DL2]:

\br{rvoa}
{\em If $G=0$, the notion of generalized vertex algebra reduces 
to the notion of vertex algebra.
If $G=\Z/2\Z$ with $(m+2\Z,n+2\Z)=mn+2\Z$ for $m,n\in \Z$, 
the notion of generalized vertex algebra reduces 
to the notion of vertex superalgebra, noting that $c(\cdot,\cdot)=1$.}
\er

\bd{dgvoa}
{\em [DL2] A generalized vertex algebra 
$(V,Y,{\bf 1},\omega, G, c(\cdot,\cdot),(\cdot,\cdot))$ 
is called a {\em generalized vertex operator algebra}
if $G$ is finite, $c(\cdot,\cdot)=1$, $(\cdot,\cdot)$ is 
nondegenerate and if $V$ 
satisfies the following two {\em grading restrictions}:
\begin{eqnarray}
& &\dim V_{n}<\infty \;\;\;\mbox{ for }n\in \C,\\
& &V_{n}=0\;\;\;\mbox{ for $n$ whose real part is sufficiently small}.
\end{eqnarray}}
\ed

\br{ronG}
{\em Note that it is possible that $V^{g}=0$ for some $g\in G$.
For example, $V=V^{0}$, i.e., $V^{g}=0$ for $0\ne g\in G$. 
Denote by $H$ the subgroup of $G$ generated by $g$ with $V^{g}\ne 0$.
(For those studied in [DL2], $H=\{ g\in G\;|\; V^{g}\ne 0\}\ne G$.)
Then $V$ is an $H$-graded space.
However, the form $(\cdot,\cdot)$ restricted to $H$ may be 
degenerate.}
\er

\br{rnondegenerate}
{\em  We here show that on the other hand, the nondegeneracy
of $(\cdot,\cdot)$ can always be achieved 
by enlarging the grading group $G$.
Let $G$ be a finite abelian group and $(\cdot,\cdot)$ 
a degenerate symmetric $\C/2\Z$-valued $\Z$-bilinear form on $G$.
Let $F$ be a free covering group of $G$ of finite rank with 
covering map $\pi$. Through $\pi$, $(\cdot,\cdot)$ lifts to 
a $\C/2\Z$-valued $\Z$-bilinear form $(\cdot,\cdot)$ on $F$.
Furthermore, by using a basis of $F$ we can lift $(\cdot,\cdot)$ 
to a nondegenerate symmetric $\C$-valued $\Z$-bilinear form 
$\<\cdot,\cdot\>$ on $F$. Since $G$ is finite, it is easy to see that
$\<\cdot,\cdot\>$ is ${\bf Q}$-valued. Set 
\begin{eqnarray}
{\bf h}=\C\otimes_{\Z}F,
\end{eqnarray}
and then extend $\<\cdot,\cdot\>$ to be a symmetric $\C$-bilinear form on 
${\bf h}$.
Let $F_{0}$ be the kernel of $\pi$ and denote by $F_{0}^{o}$ the dual
lattice of $F_{0}$ in ${\bf h}$. Then $F\subset 2F_{0}^{o}$.
Set
\begin{eqnarray}
\tilde{G}=2F_{0}^{o}/F_{0}.
\end{eqnarray}
Since $|G|F\subset F_{0}$,
$F_{0}$ must also span ${\bf h}$ over $\C$. 
It follows that $\tilde{G}$ is finite.
Then $G$ is a natural subgroup of $\tilde{G}$ and
$\<\cdot,\cdot\>$ gives rise to a natural
 nondegenerate symmetric $\C/2\Z$-valued $\Z$-bilinear form
on $\tilde{G}$  whose restriction on $G$ gives $(\cdot,\cdot)$.}
\er

\bp{pdlduality}
[DL2] In the presence of all the axioms
except the generalized Jacobi identity in defining the notion 
of generalized vertex algebra,
the generalized Jacobi identity is equivalent to the following
generalized weak commutativity and associativity:

(A) For $g_{1},g_{2}\in G$ and $v_{1}\in V^{g_{1}},\; v_{2}\in V^{g_{2}}$,
there exists $k\in \N$ such that
\begin{eqnarray}
& &(z_{1}-z_{2})^{k+(g_{1},g_{2})}Y(v_{1},z_{1})Y(v_{2},z_{2})\nonumber\\
&=&(-1)^{k}c(g_{1},g_{2})(z_{2}-z_{1})^{k+(g_{1},g_{2})}
Y(v_{2},z_{2})Y(v_{1},z_{1}).
\end{eqnarray}
(B) For $g_{1},g_{2},h\in G$ and 
$v_{1}\in V^{g_{1}},\; v_{2}\in V^{g_{2}},\; w\in V^{h}$,
there exists $l\in \N$ such that
\begin{eqnarray}
(z_{0}+z_{2})^{l+(g_{1},h)}Y(v_{1},z_{0}+z_{2})Y(v_{2},z_{2})w
=(z_{2}+z_{0})^{l+(g_{1},h)}Y(Y(v_{1},z_{0})v_{2},z_{2})w,
\end{eqnarray}
where $l$ is independent of $v_{2}$.
\ep

\bd{dmodule1}
{\em [DL2] A {\em $V$-module} is a vector space 
$W=\coprod_{s\in S}W^{s}$, where $S$ is a $G$-set with action of $G$ 
on $S$ denoted by $+$,
equipped with a $\C/\Z$-valued function $(\cdot,\cdot)$
on $G\times S$ such that in $\C/\Z$,
\begin{eqnarray}
(g_{1}+g_{2},g_{3}+s)=(g_{1},g_{3})+(g_{2},g_{3})+(g_{1},s)+(g_{2},s)
\end{eqnarray}
for $g_{1},g_{2},g_{3}\in G,\;s\in S$, such that all the axioms
defining the notion of generalized vertex algebra that make sense
hold for $W$ with obvious modifications.
Denote this module by $(W,Y,S,(\cdot,\cdot))$.}
\ed

When $V$ is a generalized vertex operator algebra, a {\em $V$-module}
is a module for $V$ as a generalized vertex algebra
such that the two grading restrictions on the $\C$-gradation 
hold.
A {\em weak module} for a generalized vertex operator algebra $V$
is a module for $V$ as a generalized vertex algebra.

\br{rconsequence}
{\em In the definition of a $V$-module, the axioms (\ref{e2.25}) 
and (\ref{e2.29}) 
are consequences of the others (cf. [FHL], [DLM2]).}
\er

A {\em homomorphism} from $(W_{1},Y,S_{1},(\cdot,\cdot)_{1})$
to $(W_{2},Y,S_{2},(\cdot,\cdot)_{2})$ consists of a $G$-set map 
$\bar{f}$ from  $S_{1}$ to $S_{2}$ and
a linear map
$f$ from $W_{1}$ to $W_{2}$ such that 
\begin{eqnarray}
& &(g,s)_{1}=(g,\bar{f}(s))_{2}\;\;\;\mbox{
for }g\in G,\; s\in S_{1},\\
& &f(W_{1}^{s})\subset W_{2}^{\bar{f}(s)}\;\;\;\mbox{
for }s\in S_{1},\\
& &f(Y(v,z)w)=Y(v,z)f(w)\;\;\;\mbox{ for }v\in V,\; w\in W_{1}.
\end{eqnarray}
It is called an {\em isomorphism} if both $f$ and $\bar{f}$ 
are one-to-one and onto.

Next we shall introduce the notion of a submodule of $(W,Y,S,(\cdot,\cdot))$.
A {\em submodule} of $(W,Y,S,(\cdot,\cdot))$
is an $S$-graded subspace which is stable
under the action of $v_{n}$ for $v\in V,\; n\in \C$. 
Later, we shall need the following a little more general notion.

Set
\begin{eqnarray}
& &G_{0}=\{ g\in G\;|\; c(g,h)=1,\; (g,h)=0+2\Z\;\;\;\mbox{ for }h\in G\},\\
& &G_{1}=\{ g\in G\;|\; (g,h)=0+2\Z\;\;\;\mbox{ for }h\in G\}.
\end{eqnarray}
Clearly, $G_{0}$ and $G_{1}$ are subgroups of $G$.

Let $A\subset G_{1}$ be a subgroup. Then 
for any $V$-module $(W,Y,S,(\cdot,\cdot))$,
$W$ can be naturally $S/A$-graded as
\begin{eqnarray}
W=\coprod_{A+s\in S/A}W^{A+s},\;\;\;\mbox{ where }
\;\;W^{A+s}=\coprod_{\a\in A}W^{\a+s}.
\end{eqnarray}
Furthermore, the form $(\cdot,\cdot)$ on $G\times S$ reduces 
to $G\times (S/A)$.
Thus, $(W,Y,S/A,(\cdot,\cdot))$ is a $V$-module.
We shall use the term ``an $S/A$-graded submodule of $W$'' 
for a submodule of $(W,Y,S/A,(\cdot,\cdot))$.

%By analogy a submodule of $(W,Y,S,(\cdot,\cdot))$
%should be a subspace that is stable of 
%the action of $v_{n}$ for $v\in V,\; n\in \C$. 
%Note that in the generalized Jacobi identity, the third term
%$$z_{2}^{-1}\delta\left(\frac{z_{1}-z_{0}}{z_{2}}\right)
%\left(\frac{z_{1}-z_{0}}{z_{2}}\right)^{-g}$$
%needs special attention.

We define an {\em ideal} of $V$ to be a $G/A$-graded submodule of $V$
for some subgroup $A$ of $G_{0}$.
Recall from [GL] the following skew-symmetry:
\begin{eqnarray}
Y(u,z)v=c(g,h)e^{\pi i (g,h)}e^{zL(-1)}Y(v,e^{\pi i}z)u
\end{eqnarray}
for $u\in V^{g},\; v\in V^{h},\; g,h\in G$.
Then for an ideal $I$ of $V$ we have
\begin{eqnarray}
Y(u,z)v,\;\;Y(v,z)u\in I\{z\}\;\;\;\mbox{ for }u\in V,\; v\in I,
\end{eqnarray}
and both the forms $c(\cdot,\cdot)$ and $(\cdot,\cdot)$ reduce to
$G/A$. Therefore, $V/I$ is a generalized vertex algebra with
$G/A$ as the grading group.

%In ideal cases, we can choose a larger grading group $G$ such that 
%every $V$-module is $G$-graded.

%Conversely, let $S$ be a $G$-set equipped with a 
%$\C/2\Z$-valued $\Z$-linear function on $G\times S$ 
%such that.. Then consider the product group $G\times \Z S$.
%Consider the relation: $(g,s)=(0,g+s)$ for $g\in G,\; s\in S$.

\section{Abelian coset generalized vertex algebras $\Omega_{V}$}
In this section, we shall establish a general
abelian coset theory in the context of generalized vertex algebras.
Let $V$ be a general vertex operator algebra with a vertex operator 
subalgebra $M(\ell)$ associated to an affine Lie algebra $\hat{\bf h}$
of level $\ell$ where ${\bf h}$ is a finite-dimensional 
abelian Lie algebra.
($M(\ell)$ and $V$ possibly have different Virasoro vectors.)
We assume that $V$ is a semisimple ${\bf h}$-module with $h\in {\bf h}$
being represented by $h(0)$.
By extending the construction of $Z$-operators in 
[LW4-6] (or $U$-operators in [LP2])
we define a vertex operator map $Y_{\Omega}$ on $V$
and we show that the vacuum space $\Omega_{V}$ of $\hat{\bf h}$, namely
the space of highest weight vectors for $\hat{\bf h}$, 
equipped with $Y_{\Omega}$
has a natural generalized vertex algebra structure with $V$ as a module.
Similarly, we show that for a $V$-module $W$, 
$W$ is a natural $\Omega_{V}$-module
with $\Omega_{W}$ as a submodule. Furthermore, we show that for 
irreducible $V$-modules $W_{1}$ and $W_{2}$,
$\Omega_{W_{1}}\simeq \Omega_{W_{2}}$ if and only if 
$W_{1}^{(\a)}\simeq W_{2}$ for some $\a\in {\bf h}$, 
where $W_{1}^{(\a)}$ is a $V$-module which was constructed in [Li3]
as a deformation of $W_{1}$.
Motivated by the $Z$-algebra theory developed in [LW4-6] and [LP1-2],
we define a notion of $\Omega_{V}$-${\bf h}$-module
and we prove that the category of $V$-modules is naturally
equivalent to the category of $\Omega_{V}$-${\bf h}$-modules.

Now we start to establish our basic notations.
Let ${\bf h}$ be a $d$-dimensional vector space 
equipped with a nondegenerate
symmetric bilinear form $\<\cdot,\cdot\>$. 
Viewing ${\bf h}$ as an abelian Lie algebra, we 
have the corresponding affine Lie algebra
\begin{eqnarray}
\hat{\bf h}={\bf h}\otimes \C[t,t^{-1}]\oplus\C c,
\end{eqnarray}
where
\begin{eqnarray}
& &[c, \hat{\bf h}]=0,\nonumber\\
& &[a\otimes t^{m},b\otimes t^{n}]=m\<a,b\>\delta_{m+n,0}c
\;\;\;\mbox{ for }a,b\in {\bf h},\;m,n\in \Z.
\end{eqnarray}
Following the tradition we also use $h(m)$ for $h\otimes t^{m}$.

Set
\begin{eqnarray}
\hat{\bf h}^{+}={\bf h}\otimes t\C[t],\;\;\;\;\; 
\hat{\bf h}^{-}={\bf h}\otimes t^{-1}\C[t^{-1}].
\end{eqnarray}
The subalgebra
\begin{eqnarray}
\hat{\bf h}_{\Z}=\hat{\bf h}^{+}\oplus \hat{\bf h}^{-}\oplus\C c
\end{eqnarray}
of $\hat{\bf h}$ is a Heisenberg algebra.

Let $\ell$ be a nonzero complex number. Consider
the induced irreducible $\hat{\bf h}$-module, 
irreducible even under $\hat{\bf h}_{\Z}$,
\begin{eqnarray}
M(\ell)
=U(\hat{\bf h})\otimes _{U({\bf h}\otimes \C[t]\oplus \C c)}\C
\simeq S(\hat{\bf h}^{-}) \;\;\;\mbox{(linearly)},
\end{eqnarray}
${\bf h}\otimes \C[t]$ acting trivially on $\C$ and $c$ acting as $\ell$.
An $\hat{\bf h}$-module on which $c$ acts as $\ell$ is
called a {\em level } $\ell$ module.
Then $M(\ell)$ is of level $\ell$.

\br{rleveloneell}
{\em Note that a level-$\ell$ module for the affine Lie algebra 
$\hat{\bf h}$ with respect to $\<\cdot,\cdot\>$ amounts
to a level-$1$ module for the affine Lie algebra $\hat{\bf h}$ 
with respect to $\ell\<\cdot,\cdot\>$.}
\er

Let $\{\b_{1},\dots,\b_{d}\}$ be an orthonormal basis 
of ${\bf h}$ with respect to the bilinear form $\<\cdot,\cdot\>$.
(Then $\{{1\over\sqrt{\ell}}\b_{1},\dots,{1\over\sqrt{\ell}}\b_{d}\}$ 
is an orthonormal basis with respect to $\ell\<\cdot,\cdot\>$, where 
$\sqrt{\ell}$ is a square root of $\ell$.)
Set
\begin{eqnarray}
\omega_{\bf h}
={1\over 2\ell}\sum_{i=1}^{d}\b_{i}(-1)\b_{i}(-1)1\in M(\ell).
\end{eqnarray}
Then from [FLM] $M(\ell)$ is a vertex operator algebra with 
the Virasoro vector
$\omega_{\bf h}$ of rank $d$. Furthermore, a weak $M(\ell)$-module
structure on a vector space $W$ amounts to a level-$\ell$ 
$\hat{\bf h}$-module structure 
on $W$ such that for every $w\in W$, $({\bf h}\otimes t^{n})w=0$ for $n$ 
sufficiently large.

Set
\begin{eqnarray}
Y(\omega_{\bf h},z)=\sum_{n\in \Z}L_{\bf h}(n)z^{-n-2}.
\end{eqnarray}
Then (cf. [FLM])
\begin{eqnarray}\label{eL(m)h(n)1}
[L_{\bf h}(m),h(n)]=-nh(m+n)\;\;\;\mbox{ for }h\in {\bf h},\; m,n\in \Z.
\end{eqnarray}
The following results are well known ([FLM], (1.9.56), (8.7.23), (8.7.25)):

\bl{lvirasoro}
Let $\ell\in \C^{\times}$ and let $W$ be a weak $M(\ell)$-module, 
or what is equivalent to, a level-$\ell$ $\hat{\bf h}$-module
such that for every $w\in W$, $({\bf h}\otimes t^{n})w=0$ 
for $n$ sufficiently large.
Let 
$a\in W$ be a highest weight vector of $W$ of weight $\a\in {\bf h}$, i.e.,
\begin{eqnarray}
h(n)a=\<h,\a\>\delta_{n,0}a\;\;\;\mbox{ for }h\in {\bf h},\; n\ge 0.
\end{eqnarray}
Then
\begin{eqnarray}
& &L_{\bf h}(0)a={\<\a,\a\>\over 2\ell}a,\label{el0}\\
& &L_{\bf h}(-1)a={1\over \ell}\a(-1)a.\label{el-1}
\end{eqnarray}
\el

%\pf (\ref{el0}) follows from [FLM] (formula (1.9.56)).
%The formula was observed in [FLM]. It can also be proved as follows:
%Set $v=L_{\bf h}(-1)a-{1\over \ell}\a(-1)a$. Clearly,
%\begin{eqnarray}
%h(n)v=0\;\;\;\mbox{ for }h\in {\bf h},\; n\ge 2.
%\end{eqnarray}
%Furthermore, using the properties
%\begin{eqnarray}
%[L_{\bf h}(-1),h(n)]=-nh(n-1),\;\;\;\;\; [h(1),\a(-1)]=\ell \<h,\a\>,
%\end{eqnarray}
%we get
%\begin{eqnarray}
%h(1)v=h(1)L_{\bf h}(-1)a-{1\over \ell}h(1)\a(-1)a
%=h(0)a-\<h,\a\>a=0.
%\end{eqnarray}
%This shows that $v$ is a lowest weight vector
%of $M_{\a}$ if $v\ne 0$. Because $M_{\a}$ is irreducible, 
%we must have $v=0$. This proves (\ref{el-1}). $\;\;\;\;\Box$

{\bf Basic Assumption 1:} Throughout Sections 3-5,  we assume that
$V$ is a vertex operator algebra, $\ell \in \C^{\times}$ and 
that  $M(\ell)$ is a vertex operator subalgebra 
generated by the subspace ${\bf h}$ of $V_{(1)}$
with the Virasoro vector $\omega_{\bf h}$ which is possibly 
different from the Virasoro vector $\omega$ of $V$. We also assume that 
$V$ is a semisimple ${\bf h}$-module with $h\in {\bf h}$ 
being represented by $h(0)$ and that
\begin{eqnarray}
L(n){\bf h}=0\;\;\;\mbox{ for }n\ge 1.
\end{eqnarray}
It is known (cf. [FHL]) that the last condition together 
with ${\bf h}\subset V_{(1)}$ is equivalent to
\begin{eqnarray}\label{eL(m)h(n)}
[L(m),h(n)]=-nh(m+n)\;\;\;\mbox{ for }h\in {\bf h},\; m,n\in \Z.
\end{eqnarray}

For $\a\in {\bf h}$, set
\begin{eqnarray}
V^{\a}=\{v\in V\;|\; h(0)v=\<\a,h\>v\;\;\;\mbox{ for }h\in {\bf h}\}.
\end{eqnarray}
Using the general notation $\<\cdot\>$ for ``group generated by,'' we 
define
\begin{eqnarray}
L=\< \a\in {\bf h}\;|\; V^{\a}\ne 0\>,
\end{eqnarray}
a subgroup of the additive group ${\bf h}$.
Then
\begin{eqnarray}
V=\coprod_{\a\in L}V^{\a}.
\end{eqnarray}
When $V$ is simple, one can show  that 
$L=\{\a\in {\bf h}\;|\; V^{\a}\ne 0\}$ (cf. [LX]).

Define 
\begin{eqnarray}
\Omega_{V}=\{v\in V\;|\; h(n)v=0\;\;\;\mbox{ for }h\in {\bf h}, \; n\ge 1\}.
\end{eqnarray}
Since $[h(0), h'(m)]=0$ for $h,h'\in {\bf h},\;m\in \Z$, 
$h(0)$ preserves $\Omega_{V}$.
Then
\begin{eqnarray}
\Omega_{V}=\coprod_{\a\in L}\Omega_{V}^{\a},
\;\;\;\mbox{ where }\;\;\Omega_{V}^{\a}=\Omega_{V}\cap V^{\a}.
\end{eqnarray}
Since ${\bf h}$ commutes with $\hat{\bf h}$, $\Omega_{V}^{\a}$ is an
$\hat{\bf h}$-submodule of $V$ for $\a\in L$.
For $h\in {\bf h},\; v\in \Omega_{V}^{\a}$, it follows from 
the commutator formula (see [B], [FLM]) that
\begin{eqnarray}\label{e3.17}
[h(n), Y(v,z)]=\<h,\a\>z^{n}Y(v,z)\;\;\;\mbox{ for }n\in \Z.
\end{eqnarray}

\bd{dVcategory}
{\em Define ${\cal{C}}$ to be the category of weak 
$V$-modules $W$ on which ${\bf h}$ acts semisimply and 
$\hat{\bf h}^{+}$ acts locally nilpotently.}
\ed

It follows from [LW3] and [K] that each $W$ from ${\cal{C}}$ 
is a completely reducible $\hat{\bf h}$-module. Then
\begin{eqnarray}
W=U(\hat{\bf h}^{-})\Omega_{W}=U(\hat{\bf h}^{-})\otimes \Omega_{W},
\end{eqnarray}
where 
\begin{eqnarray}
\Omega_{W}=\{w\in W\;|\; h(n)w=0\;\;\;\mbox{ for }h\in {\bf h},\; n\ge 1\}.
\end{eqnarray}
As $V$ satisfies the two grading restrictions ([FLM], [FHL]), 
the adjoint module $V$ is in the category ${\cal{C}}$.

For $v\in V^{\a},\;\a\in L$, following [LW2] we
set
\begin{eqnarray}
Z(v,z)=E^{-}({1\over \ell}\a,z)Y(u,z)
E^{+}({1\over \ell}\a,z),
\end{eqnarray}
where for $h\in {\bf h}$,
\begin{eqnarray}
E^{\pm}(h,z)=\exp \left(\sum_{n\in \pm \Z_{+}}{h(n)\over n}z^{-n}\right).
\end{eqnarray}
Then for any $W\in {\cal{C}}$, $Z(v,z)$ is a canonical element
of $(\End W)[[z,z^{-1}]]$.
Furthermore, with (\ref{e3.17}), from [LP1] and [LW5] we immediately have:

\bl{llw1}
Let $W\in {\cal{C}},\;h\in {\bf h},\;v\in \Omega_{V}^{\a}$ for $\a\in L$. 
On $W$,
\begin{eqnarray}
[h(n), Z(v,z)]=\delta_{n,0}\<h,\a\>Z(v,z)\;\;\;\mbox{ for }n\in \Z.
\;\;\;\;\Box
\end{eqnarray}
\el

For a $V$-module $W\in {\cal{C}}$, let
\begin{eqnarray}
Y_{\Omega}(\cdot,z): \Omega_{V}\mapsto (\End W)\{z\}
\end{eqnarray}
be the linear map given by
\begin{eqnarray}
Y_{\Omega}(u,z)=Z(u,z)z^{-{1\over \ell}\a(0)}=
E^{-}({1\over \ell}\a,z)Y(u,z)
E^{+}({1\over \ell}\a,z)z^{-{1\over \ell}\a(0)}
\end{eqnarray}
for $u\in V^{\a},\;\a\in L$. In particular, we have
\begin{eqnarray}
Y_{\Omega}(u,z)=Y(u,z)\;\;\;\mbox{ for }u\in \Omega_{V}^{0}.
\end{eqnarray}
Clearly, the following lower truncation condition holds:
\begin{eqnarray}\label{elowertruncation}
Y_{\Omega}(u,z)z^{{1\over\ell}\a(0)}w\in W((z))\;\;\;\mbox{ for } 
u\in \Omega_{V}^{\a},\; \a\in L,\; w\in W.
\end{eqnarray}

Let $W\in {\cal{C}}$. For $\lambda\in {\bf h}$, we define
\begin{eqnarray}
W^{\lambda}=\{w\in W\;|\; h(0)w=\<\lambda,h\>w
\;\;\;\mbox{ for }h\in {\bf h}\}.
\end{eqnarray}
Set 
\begin{eqnarray}
P(W)=\cup_{\lambda\in {\bf h},\; W^{\lambda}\ne 0}(\lambda+L),
\end{eqnarray}
i.e., $P(W)$ is the $L$-subset of ${\bf h}$ generated by $\lambda$ 
with $W^{\lambda}\ne 0$.
If $W$ is irreducible, one can easily show that
$P(W)=\lambda+L$ for any $\lambda\in {\bf h}$ with $W^{\lambda}\ne 0$.

In view of Lemma \ref{llw1}, we immediately have:

\bl{lextlw1}
Let $W\in {\cal{C}}$ be a $V$-module. Then
for $h\in {\bf h},\; n\in \Z,\; u\in \Omega_{V}^{\a},\;\a\in L,\;
\lambda\in P(W)$,
\begin{eqnarray}
& &[h(n), Y_{\Omega}(u,z)]=\delta_{n,0}\<h,\a\>Y_{\Omega}(u,z),
\label{ehmyomega}\\
& &Y_{\Omega}(u,z)\Omega_{W}^{\lambda}\subset \Omega_{W}^{\a+\lambda}\{z\}.
\;\;\;\;\Box\label{e3.29}
\end{eqnarray}
\el

Set
\begin{eqnarray}
\omega_{\Omega}=\omega-\omega_{\bf h}\in V_{(2)}
\end{eqnarray}
and
\begin{eqnarray}
Y(\omega_{\Omega},z)=\sum_{n\in \Z}L_{\Omega}(n)z^{-n-2}.
\end{eqnarray}
{}From (\ref{eL(m)h(n)1}) and (\ref{eL(m)h(n)}) we immediately have
\begin{eqnarray}\label{elomegamhn}
[L_{\Omega}(m), h(n)]=0\;\;\;\mbox{ for }h\in {\bf h},\; m,n\in \Z.
\end{eqnarray}
Using (\ref{elomegamhn}) we get
\begin{eqnarray}
h(n)\omega_{\Omega}=h(n)L_{\Omega}(-2){\bf 1}=L_{\Omega}(-2)h(n){\bf 1}=0
\;\;\;\mbox{ for }h\in {\bf h},\; n\ge 0.
\end{eqnarray}
Then $\omega_{\Omega}\in \Omega_{V}^{0}$. Hence
$Y_{\Omega}(\omega_{\Omega},z)=Y(\omega_{\Omega},z)$.
Furthermore, from (\ref{elomegamhn}) we have
\begin{eqnarray}
[L_{\Omega}(m), L_{\bf h}(n)]=0\;\;\;\mbox{ for }m,n\in \Z.
\end{eqnarray}
Then it follows (cf. [FZ], [DL2], [GKO1-2]) that component operators 
$L_{\Omega}(n)$ satisfy the Virasoro relations
with central charge ${\rm rank }V-d$. 
To summarize we have:

\bl{lcosetvirasoro}
The relation (\ref{elomegamhn}) holds and
\begin{eqnarray}
\omega_{\Omega}=\omega-\omega_{\bf h}\in \Omega_{V}^{0}.
\end{eqnarray}
Furthermore, 
\begin{eqnarray}
[L_{\Omega}(m),L_{\Omega}(n)] =(m-n)L_{\Omega}(m+n)
+{1\over 12}(m^{3}-m)\delta_{m+n,0}
({\rm rank }V-d)
\end{eqnarray}
for $m,n\in \Z$, where $d=\dim {\bf h}$. $\;\;\;\;\Box$
\el

\bp{pderivative}
Let $W\in {\cal{C}}$ be a $V$-module. On $W$, we have
\begin{eqnarray}
[L_{\Omega}(-1),Y_{\Omega}(a,z)]=Y_{\Omega}(L_{\Omega}(-1)a,z)
={d\over dz}Y_{\Omega}(a,z)
\;\;\;\mbox{ for }a\in \Omega_{V}.
\end{eqnarray}
\ep

\pf For $u\in V$, set
\begin{eqnarray}
Y(u,z)^{+}=\sum_{n\ge 0}u_{n}z^{-n-1},\;\;\;\;
Y(u,z)^{-}=\sum_{n\ge 0}u_{-n-1}z^{n}.
\end{eqnarray}
Then for $u,v\in V$,
\begin{eqnarray}\label{enormalorderproduct}
Y(u_{-1}v,z)=Y(u,z)^{-}Y(v,z)+Y(v,z)Y(u,z)^{+}
\end{eqnarray}
(cf. formula (13.24) of [DL2]). Note that
\begin{eqnarray}\label{ebracketderivative1}
[L_{\Omega}(-1),Y(v,z)]
=Y(L_{\Omega}(-1)v,z)\;\;\;\mbox{ for }v\in V
\end{eqnarray}
because $Y_{\Omega}(\omega_{\Omega},z)=Y(\omega_{\Omega},z)$
and $L_{\Omega}(-1)=(\omega_{\Omega})_{0}$.

Let $a\in \Omega_{V}^{\a},\;\a\in L$. Then
using (\ref{elomegamhn}), (\ref{ebracketderivative1}),
Lemma \ref{lvirasoro} and 
(\ref{enormalorderproduct}) we obtain
\begin{eqnarray}
& &[L_{\Omega}(-1),Y_{\Omega}(a,z)]\nonumber\\
&=&E^{-}({1\over \ell}\a,z)[L_{\Omega}(-1), Y(a,z)]
E^{+}({1\over \ell}\a,z)z^{-{1\over\ell} \a(0)}
\nonumber\\
&=&E^{-}({1\over\ell}\a,z)Y(L_{\Omega}(-1)a,z)
E^{+}({1\over \ell}\a,z)z^{-{1\over\ell} \a(0)}
\nonumber\\
&(=&Y_{\Omega}(L_{\Omega}(-1)a,z),\;\;\mbox{ noting that from (\ref{e3.29}) } 
L_{\Omega}(-1)a\in \Omega_{V}^{\a}.)\nonumber\\
&=&E^{-}({1\over\ell}\a,z)(Y(L(-1)a,z)-Y(L_{\bf h}(-1)a,z))
E^{+}({1\over \ell}\a,z)
z^{-{1\over\ell} \a(0)}\nonumber\\
&=&E^{-}({1\over\ell}\a,z)\left({d\over dz}Y(a,z)\right)
E^{+}({1\over \ell}\a,z)
z^{-{1\over\ell} \a(0)}\nonumber\\
& &-{1\over\ell} E^{-}({1\over\ell}\a,z)Y(\a (-1)a,z)
E^{+}({1\over \ell}\a,z)z^{-{1\over\ell} \a(0)}\nonumber\\
&=&E^{-}({1\over\ell}\a,z)\left({d\over dz}Y(a,z)\right)
E^{+}({1\over \ell}\a,z)
z^{-{1\over\ell} \a(0)}\nonumber\\
& &-{1\over\ell} \a(z)^{-}E^{-}({1\over\ell}\a,z)Y(a,z)
E^{+}({1\over \ell}\a,z)
z^{-{1\over\ell} \a(0)}\nonumber\\
& &-{1\over\ell} E^{-}({1\over\ell}\a,z)Y(a,z)\a(z)^{+}
E^{+}({1\over \ell}\a,z)z^{-{1\over\ell} \a(0)}\nonumber\\
&=&{d\over dz}\left(E^{-}({1\over\ell}\a,z)Y(a,z)E^{+}({1\over \ell}\a,z)
z^{-{1\over\ell} \a(0)}\right)\nonumber\\
&=&{d\over dz}Y_{\Omega}(a,z).
\end{eqnarray}
This completes the proof.$\;\;\;\;\Box$

Now we present our first main result of this paper.
%Then $P(W)$ is an $L$-set and ${1\over\ell}\<\cdot,\cdot\>$ is a 
%$\C$-valued function on $L\times P(W)$. 

\bt{tcosetalgebra}
Let $W\in {\cal{C}}$ be a $V$-module and let
$u\in \Omega_{V}^{\a_{1}},\; v\in \Omega_{V}^{\a_{2}},\;w\in W^{\lambda}$
with $\a_{1},\a_{2}\in L,\; \lambda\in P(W)$. Then
\begin{eqnarray}\label{egjacobithem}
& &z_{0}^{-1}
\left(\frac{z_{1}-z_{2}}{z_{0}}\right)^{{1\over\ell}\<\a_{1},\a_{2}\>}
\delta\left(\frac{z_{1}-z_{2}}{z_{0}}\right)
Y_{\Omega}(u,z_{1})Y_{\Omega}(v,z_{2})w\nonumber\\
&-&z_{0}^{-1}
\left(\frac{z_{2}-z_{1}}{z_{0}}\right)^{{1\over\ell}\<\a_{1},\a_{2}\>}
\delta\left(\frac{z_{2}-z_{1}}{-z_{0}}\right)
Y_{\Omega}(v,z_{2})Y_{\Omega}(u,z_{1})w\nonumber\\
&=&z_{2}^{-1}
\left(\frac{z_{1}-z_{0}}{z_{2}}\right)^{-{1\over\ell}\<\a_{1},\lambda\>}
\delta\left(\frac{z_{1}-z_{0}}{z_{2}}\right)
Y_{\Omega}(Y_{\Omega}(u,z_{0})v,z_{2})w.
\end{eqnarray}
\et

\pf Since $W=U(\hat{\bf h}^{-})\Omega_{W}$ and 
$Y_{\Omega}(a,z)$ for $a\in \Omega_{V}$
commutes with the action of $U(\hat{\bf h}^{-})$
by (\ref{ehmyomega}), it is enough to prove 
(\ref{egjacobithem}) for $w\in \Omega_{W}^{\lambda}$.

We first prove a conjugation formula. 
Let $h\in {\bf h},\; u\in \Omega_{V}^{\a}$. 
%By [Li5] (Lemma 3.6, or more specifically, formula (3.35)) we have
%\begin{eqnarray}
%& &E^{-}(h,z_{2})Y(u,z_{1})E^{-}(-h,z_{2})\nonumber\\
%&=&Y\left((1-z_{2}/z_{1})^{h(0)}E^{+}(-h,-z_{1}+z_{2})E^{+}(h,-z_{1})u,z_{2}
%\right).
%\end{eqnarray}
%In particular, if $u\in \Omega_{V}^{\a}$, we have
Using (\ref{e3.17}) we get
\begin{eqnarray}
\left[\sum_{n\ge 1}{h(-n)\over -n}z_{2}^{n}, Y(u,z_{1})\right]
= \sum_{n\ge 1}{\<h,\a\>\over -n}z_{1}^{-n}z_{2}^{n}Y(u,z_{1})
=\left(\log (1-z_{2}/z_{1})^{\<h,\a\>}\right)Y(u,z_{1}).\;\;\;
\end{eqnarray}
Then
\begin{eqnarray}\label{econjugation2}
E^{-}(h,z_{2})Y(u,z_{1})E^{-}(-h,z_{2})
=(1-z_{2}/z_{1})^{\<\a,h\>}Y(u,z_{2}).
\end{eqnarray}

Let $u\in \Omega_{V}^{\a_{1}},\; v\in \Omega_{V}^{\a_{2}},\; 
w\in \Omega_{V}^{\a_{3}}$. Using (\ref{econjugation2}) and the fact that
\begin{eqnarray}
z_{1}^{h(0)}Y(v,z_{2})=Y(v,z_{2})z_{1}^{\<\a_{2},h\>+h(0)}
\;\;\;\mbox{ for }h\in {\bf h},
\end{eqnarray}
we get
\begin{eqnarray}
& &Y_{\Omega}(u,z_{1})Y_{\Omega}(v,z_{2})w\nonumber\\
&=&E^{-}({1\over\ell}\a_{1},z_{1})Y(u,z_{1})z_{1}^{-{1\over\ell}\a_{1}(0)}
E^{-}({1\over\ell}\a_{2},z_{2})Y(v,z_{2})z_{2}^{-{1\over\ell}\a_{2}(0)}w
\nonumber\\
&=&(1-z_{2}/z_{1})^{-{1\over\ell}\<\a_{1},\a_{2}\>}
E^{-}({1\over\ell}\a_{1},z_{1})E^{-}({1\over\ell}\a_{2},z_{2})\nonumber\\
& &\times Y(u,z_{1})Y(v,z_{2})
z_{1}^{-{1\over\ell}\<\a_{1},\a_{2}\>-{1\over\ell}\a_{1}(0)}
z_{2}^{-{1\over\ell}\a_{2}(0)}w
\nonumber\\
&=&(z_{1}-z_{2})^{-{1\over\ell}\<\a_{1},\a_{2}\>}
E^{-}({1\over\ell}\a_{1},z_{1})E^{-}({1\over\ell}\a_{2},z_{2})
Y(u,z_{1})Y(v,z_{2})
z_{1}^{-{1\over\ell}\a_{1}(0)}z_{2}^{-{1\over\ell}\a_{2}(0)}w.\;\;\;\;
\end{eqnarray}
Then
\begin{eqnarray}\label{ej1}
& &z_{0}^{-1}
\left(\frac{z_{1}-z_{2}}{z_{0}}\right)^{{1\over\ell}\<\a_{1},\a_{2}\>}
\delta\left(\frac{z_{1}-z_{2}}{z_{0}}\right)
Y_{\Omega}(u,z_{1})Y_{\Omega}(v,z_{2})w\nonumber\\
&=&z_{0}^{-{1\over\ell}\<\a_{1},\a_{2}\>}
E^{-}({1\over\ell}\a_{1},z_{1})E^{-}({1\over\ell}\a_{2},z_{2})\nonumber\\
& &\times z_{0}^{-1}\delta\left(\frac{z_{1}-z_{2}}{z_{0}}\right)
Y(u,z_{1})Y(v,z_{2})
z_{1}^{-{1\over\ell}\a_{1}(0)}z_{2}^{-{1\over\ell}\a_{2}(0)}w.\;\;\;\;\;
\end{eqnarray}
Symmetrically, we have
\begin{eqnarray}\label{ej2}
& &z_{0}^{-1}
\left(\frac{z_{2}-z_{1}}{z_{0}}\right)^{{1\over\ell}\<\a_{1},\a_{2}\>}
\delta\left(\frac{z_{2}-z_{1}}{-z_{0}}\right)
Y_{\Omega}(v,z_{2})Y_{\Omega}(u,z_{1})w\nonumber\\
&=&z_{0}^{-{1\over\ell}\<\a_{1},\a_{2}\>}
E^{-}({1\over\ell}\a_{1},z_{1})E^{-}({1\over\ell}\a_{2},z_{2})\nonumber\\
& &\times z_{0}^{-1}\delta\left(\frac{z_{2}-z_{1}}{-z_{0}}\right)
Y(v,z_{2})Y(u,z_{1})
z_{1}^{-{1\over\ell}\a_{1}(0)}z_{2}^{-{1\over\ell}\a_{2}(0)}w.\;\;\;\;\;
\end{eqnarray}
Recall from [Li5] (Lemma 3.5) (see also (8.6.9) of [FLM]) that
\begin{eqnarray}\label{egeneralvertexoperator}
& &Y(E^{-}(h,z_{0})a,z_{2})=E^{-}(h,z_{2}+z_{0})E^{-}(-h,z_{2})Y(a,z_{2})
\nonumber\\
& &\;\;\;\; \times E^{+}(-h,z_{2})E^{+}(h,z_{2}+z_{0})
z_{2}^{h(0)}(z_{2}+z_{0})^{-h(0)}
\end{eqnarray}
for $h\in {\bf h},\; a\in V$.
Noticing that $Y_{\Omega}(u,z_{0})v\in \Omega_{V}^{\a_{1}+\a_{2}}\{z_{0}\}$, 
then using the fundamental properties of delta function we obtain
\begin{eqnarray}
& &z_{2}^{-1}
\left(\frac{z_{1}-z_{0}}{z_{2}}\right)^{-{1\over\ell}\<\a_{1},\lambda\>}
\delta\left(\frac{z_{1}-z_{0}}{z_{2}}\right)
Y_{\Omega}(Y_{\Omega}(u,z_{0})v,z_{2})w\nonumber\\
&=&z_{2}^{-1}
\left(\frac{z_{1}-z_{0}}{z_{2}}\right)^{-{1\over\ell}\<\a_{1},\lambda\>}
\delta\left(\frac{z_{1}-z_{0}}{z_{2}}\right)\nonumber\\
& &\times E^{-}({1\over\ell}(\a_{1}+\a_{2}),z_{2})
Y(E^{-}({1\over\ell}\a_{1},z_{0})
Y(u,z_{0})z_{0}^{-{1\over\ell}\a_{1}(0)}v,z_{2})
z_{2}^{-{1\over\ell}(\a_{1}(0)+\a_{2}(0))}w \nonumber\\
&=&z_{1}^{-1}
\left(\frac{z_{2}+z_{0}}{z_{1}}\right)^{{1\over\ell}\<\a_{1},\lambda\>}
\delta\left(\frac{z_{2}+z_{0}}{z_{1}}\right)
z_{0}^{-{1\over\ell}\<\a_{1},\a_{2}\>}
\nonumber\\
& &\times E^{-}({1\over\ell}\a_{2},z_{2})E^{-}({1\over\ell}\a_{1},z_{2}+z_{0})
Y(Y(u,z_{0})v,z_{2})\nonumber\\
& &\times z_{2}^{{1\over\ell}\a_{1}(0)}
(z_{2}+z_{0})^{-{1\over\ell}\a_{1}(0)}
z_{2}^{-{1\over\ell}(\a_{1}(0)+\a_{2}(0))}w \nonumber\\
&=&z_{1}^{-1}
\left(\frac{z_{2}+z_{0}}{z_{1}}\right)^{{1\over\ell}\<\a_{1},\lambda\>}
\delta\left(\frac{z_{2}+z_{0}}{z_{1}}\right)
z_{0}^{-{1\over\ell}\<\a_{1},\a_{2}\>}
\nonumber\\
& &\times E^{-}({1\over\ell}\a_{2},z_{2})E^{-}({1\over\ell}\a_{1},z_{1})
Y(Y(u,z_{0})v,z_{2})(z_{2}+z_{0})^{-{1\over\ell}\<\a_{1},\lambda\>}
z_{2}^{-{1\over\ell}\a_{2}(0)}w\nonumber\\
&=&z_{0}^{-{1\over\ell}\<\a_{1},\a_{2}\>}
E^{-}({1\over\ell}\a_{2},z_{2})E^{-}({1\over\ell}\a_{1},z_{1})
\nonumber\\
& &\times z_{1}^{-1}\delta\left(\frac{z_{2}+z_{0}}{z_{1}}\right)
Y(Y(u,z_{0})v,z_{2})
z_{1}^{-{1\over\ell}\a_{1}(0)}z_{2}^{-{1\over\ell}\a_{2}(0)}w.\label{ej3}
\end{eqnarray}
We are using the fact that
$$E^{-}(h,z)E^{-}(h',z)=E^{-}(h+h',z),\;\;\;E^{+}(h,z)b=b
\;\;\;\mbox{ for }h,h'\in {\bf h},\; b\in \Omega_{V}.$$
It follows from (\ref{ej1}), (\ref{ej2}) and (\ref{ej3}) that
the generalized Jacobi identity of $Y_{\Omega}$ for $(u,v,w)$
is equivalent to the Jacobi identity of $Y$ for $(u,v,w)$.
The proof is complete.$\;\;\;\;\Box$

With Proposition \ref{pderivative} and Theorem \ref{tcosetalgebra} we have:

\bt{tcosetmodule}
The tuple
$(\Omega_{V}, Y_{\Omega}, {\bf 1}, \omega_{\Omega}, L, (\cdot,\cdot))$
carries the structure of a generalized vertex algebra with 
$c(\cdot,\cdot)=1$ and
\begin{eqnarray}
(\a,\b)={1\over\ell}\<\a,\b\>+2\Z\;\;\;\mbox{ for }\a,\b\in L.
\end{eqnarray}
Furthermore, for $W\in {\cal{C}}$,
$(W,Y_{\Omega})$ is an $\Omega_{V}$-module associated with
the free $L$-set $P(W)$ equipped with the $\C$-valued function 
${1\over\ell}\<\cdot,\cdot\>$ on $L\times P(W)$
with $\Omega_{W}$ as a submodule. 
\et

\pf From definition we have $Y_{\Omega}({\bf 1},z)=1$ and
\begin{eqnarray}
Y_{\Omega}(u,z){\bf 1}=E^{-}({1\over\ell}\a,z)Y(u,z){\bf 1}\in \Omega_{V}[[z]],
\end{eqnarray}
which implies that $\lim_{z\rightarrow 0}Y_{\Omega}(u,z){\bf 1}=u$.
Then it follows immediately from Proposition \ref{pderivative}
and Theorem \ref{tcosetalgebra}.$\;\;\;\;\Box$

For our need in a sequel, we prove the following:

\bp{pautoreduction}
Let $\tau$ be an automorphism of the vertex operator algebra $V$
such that $\tau({\bf h})={\bf h}$ and $\tau$ preserves the bilinear form
$\<\cdot,\cdot\>$ on ${\bf h}$. Then $\tau(V_{\Omega})=V_{\Omega}$
and $\tau$ restricted to $V_{\Omega}$ is an 
automorphism of the generalized vertex algebra $\Omega_{V}$.
\ep

\pf For $h\in {\bf h},\; n\in \Z$, we have
\begin{eqnarray}
h(n)\tau=\tau (\tau^{-1}(h))(n).
\end{eqnarray}
Then it follows that $\tau(V_{\Omega})=V_{\Omega}$
and that there is a map $\tau_{L}$ from $L$ to $L$ such that
\begin{eqnarray}
\tau (\Omega_{V}^{\a})=\Omega_{V}^{\tau_{L}(\a)}
\;\;\;\mbox{ for }\a\in L.
\end{eqnarray}
Let $v\in \Omega_{V}^{\a},\; \a\in L$. Then
\begin{eqnarray}
\tau Y_{\Omega}(v,z)\tau^{-1}
&=&\tau E^{-}({1\over \ell}\a,z)Y(v,z)E^{+}({1\over \ell}\a,z)z^{{1\over\ell}\a(0)}
\tau^{-1}\nonumber\\
&=&E^{-}({1\over \ell}\tau_{L}(\a),z)Y(\tau(v),z)
E^{+}({1\over \ell}\tau_{L}(\a),z)z^{{1\over\ell}\tau(\a)(0)}\nonumber\\
&=&Y_{\Omega}(\tau(v),z).
\end{eqnarray}
By definition, we have $\tau({\bf 1})={\bf 1}$ and $\tau(\omega)=\omega$.
Since $\tau({\bf h})={\bf h}$ and $\tau$ preserves
the form $\<\cdot,\cdot\>$, $\tau(\omega_{{\bf h}})=\omega_{{\bf h}}$.
Consequently, $\tau(\omega_{\Omega})=\omega_{\Omega}$
Then $\tau$ is an automorphism of the generalized vertex algebra $\Omega_{V}$.
$\;\;\;\;\Box$

Sometimes, we are interested in knowing how to get a set of generators of 
the generalized vertex algebra $\Omega_{V}$.
The following result is a simple generalization of 
Proposition 14.9 of [DL2]:

\bp{pgenerators}
Let $U$ be an ${\bf h}$-submodule of $\Omega_{V}$ such that $U+{\bf h}$ 
generates $V$ as a vertex algebra. Then
$U$ generates $\Omega_{V}$ in the sense that $\Omega_{V}$ is 
linearly spanned by 
\begin{eqnarray}\label{espanningset}
Z(u_{1},m_{1})\cdots Z(u_{k},m_{k}){\bf 1}
\end{eqnarray}
for $k\in \N,\; u_{i}\in U^{\a_{i}},\;\a_{i}\in L,\;m_{i}\in \Z$, where
\begin{eqnarray}
Z(u_{i},z)=Y_{\Omega}(u_{i},z)z^{{1\over\ell}\a_{i}(0)}
=\sum_{m\in \Z}Z(u_{i},m)z^{m}.
\end{eqnarray}
\ep

\pf Let $\bar{U}$ be the subspace spanned by elements of the form
(\ref{espanningset}). With $U$ being an ${\bf h}$-submodule of $\Omega_{V}$,
$\bar{U}$ is an ${\bf h}$-submodule of $\Omega_{V}$, hence
$\bar{U}$ is stable under
the action of components of $Y_{\Omega}(u,z)$ for $u\in U$.
Then 
$U(\hat{\bf h}^{-})\bar{U}=M(\ell)\otimes \bar{U}$
is stable under the actions of $\hat{\bf h}$ and the
components of $Y(u,z)$ for $u\in U$. Because by assumption
$U+{\bf h}$ generates $V$, we must have $M(\ell)\otimes \bar{U}=V$, Thus
$\bar{U}=\Omega_{V}$.
This completes the proof. $\;\;\;\;\Box$

The following definition was motivated by [LW4-6] and [LP1-2]:

\bd{dcategoryD}
{\em An {\em $\Omega_{V}$-${\bf h}$-module} $U$ is 
a semisimple ${\bf h}$-module and an $\Omega_{V}$-module 
associated to the $L$-set $S=P(U)\subset {\bf h}$ with 
$(\cdot,\cdot)={1\over\ell}\<\cdot,\cdot\>$ such that
\begin{eqnarray}
[h,Y_{\Omega}(v,z)]=\<h,\a\>Y_{\Omega}(v,z)
\;\;\;\mbox{ for }h\in {\bf h},\; v\in \Omega_{V}^{\a},\; \a\in L.
\end{eqnarray}
We denote the category of
$\Omega_{V}$-${\bf h}$-modules by ${\cal{D}}$.}
\ed

As an immediate consequence we have (cf. [LW5], Theorem 5.7):

\bc{conetoone}
Let $W\in {\cal{C}}$. The correspondences
$$M\mapsto M\cap \Omega_{W},\;\;\;\; U\mapsto U(\hat{\bf h})\cdot U $$
define mutually inverse bijections between the set
of all $V$-submodules of $W$ and the set of all 
$\Omega_{V}$-${\bf h}$-submodules of $\Omega_{V}$. In particular, 
$W$ is $V$-irreducible if and only if $\Omega_{V}$ is 
$\Omega_{V}$-${\bf h}$-irreducible.$\;\;\;\;\Box$
\ec

In view of Theorem \ref{tcosetmodule}, we have a functor
$\Omega$ from ${\cal{C}}$ to ${\cal{D}}$ which maps $W$ to $\Omega_{W}$.
Conversely, given an $\Omega_{V}$-${\bf h}$-module $U$
we shall construct a $V$-module in ${\cal{C}}$ (cf. [LW3], [LP2]).

Set
\begin{eqnarray}
E(U)=M(\ell)\otimes U\;\;\;\mbox{ as a vector space}.
\end{eqnarray}
First, view $E(U)$ as an $\hat{\bf h}$-module with
\begin{eqnarray}
& &h(n)\cdot (v\otimes w)=h(n)v\otimes w,\\
& &h(0)\cdot (v\otimes w)=h\cdot (u\otimes v)=v\otimes hw
\end{eqnarray}
for $h\in {\bf h},\;n\in \Z-\{0\},\;v\in M(\ell),\;w\in U$.
For $h\in {\bf h}$, set
\begin{eqnarray}
h(z)=\sum_{m\in\Z}h(m)z^{-m-1}\in (\End W)[[z,z^{-1}]].
\end{eqnarray}
For $a\in \Omega_{V}^{\a},\; \a\in L$, we define
\begin{eqnarray}
Y(a,z)
=E^{-}(-{1\over\ell}\a,z)E^{+}(-{1\over\ell}\a,z)\otimes
Y_{\Omega}(a,z)z^{{1\over\ell}\a(0)}\in (\End E(U))[[z,z^{-1}]].
\end{eqnarray}
(Recall that $Y_{\Omega}(a,z)z^{{1\over\ell}\a(0)}\in (\End U)[[z,z^{-1}]]$.)
Then we inductively define
\begin{eqnarray}\label{edefinitionY}
Y(h(-n-1)v,z)
={1\over n!}\left(h^{(n)}(z)^{-}Y(v,z)+Y(v,z)h^{(n)}(z)^{+}\right)
\end{eqnarray}
for $h\in {\bf h},\; n\ge 0,\; v\in V$,
where $h^{(n)}(z)$ is the $n$-th derivative of $h(z)$.
With $V=U(\hat{\bf h}^{-})\otimes \Omega_{V}$, it follows from 
induction that this really defines a
vertex operator map $Y$ from $V$ to $(\End E(U))[[z,z^{-1}]]$.

\bp{pback}
The space $E(U)$ equipped with the defined vertex operator map $Y$
is a weak $V$-module in ${\cal{C}}$ such that
$\Omega_{E(U)}$ is naturally
isomorphic to $U$ as an $\Omega_{V}$-module.
On the other hand, for any $W\in {\cal{C}}$, $E(\Omega_{W})$ 
is naturally isomorphic to $W$ as a $V$-module.
\ep

\pf Clearly, $Y(a,z)w\in E(U)((z))$ for 
$a\in \Omega_{V}^{\a},\; \a\in L,\; w\in E(U)$.
Then it follows from induction that $Y(v,z)w\in E(U)((z))$ 
for all $v\in V,\; w\in E(U)$. From the definition, we have
$Y({\bf 1},z)=1$. 
Then for proving $E(U)$ is a weak $V$-module,
it remains to prove the Jacobi identity.

First, we list  certain properties that $Y(\cdot,z)$ satisfies.

{}From definition, we have
\begin{eqnarray}
[h(m), Y_{\Omega}(a,z)]=\delta_{m,0}\<h,\a\>Y_{\Omega}(a,z)
\;\;\;\mbox{ for }h\in {\bf h},\; m\in \Z,\; 
a\in \Omega_{V}^{\a},\; \a\in L.
\end{eqnarray}
Recall from [FLM] that for $h, h'\in {\bf h}$,
\begin{eqnarray}
& &[h(0), E^{\pm}(h',z)]=0\label{e3.51}\\
& &[h(\pm m),E^{\pm}(h',z)]=0\label{e3.52}\\
& &[h(\pm m),E^{\mp}(h',z)]=-\ell\<h,h'\>z^{m}E^{\mp}(h',z)
\;\;\;\mbox{ for }m\in \Z_{+\label{e3.53}}\\
& &E^{+}(h,z_{1})E^{-}(h',z_{2})
=(1-z_{2}/z_{1})^{\ell\<h,h'\>}E^{-}(h',z_{2})E^{+}(h,z_{1}).\label{eE+E-}
\end{eqnarray}

It follows from the definition of $Y(u,z)$ and (\ref{eE+E-})
that (\ref{econjugation2}) holds. 

Using (\ref{e3.51})-(\ref{e3.53})  we get
\begin{eqnarray}\label{ecommproof}
[h(m),Y(u,z)]=\<h,\a\>z^{m}Y(u,z)\;\;\;\mbox{ for }
h\in {\bf h},\; m\in \Z,\; u\in \Omega_{V}^{\a}
\end{eqnarray}
(cf. (\ref{e3.17})). Furthermore, using (\ref{ecommproof}) we obtain
\begin{eqnarray}
& &[h(z)^{-},Y(u,z_{1})]=\<h,\a\>(z_{1}-z)^{-1}Y(u,z_{1}),\label{eh-yu}\\
& &[h(z)^{+},Y(u,z_{1})]=\<h,\a\>(z-z_{1})^{-1}Y(u,z_{1}).\label{eh+yu}
\end{eqnarray}

It follows from (\ref{edefinitionY}) that
(\ref{egeneralvertexoperator}) holds. 
(Cf. (8.6.9) of [FLM], Lemma 3.5 of [Li4].)

Writing (\ref{edefinitionY}) in terms of generating series we have
\begin{eqnarray}\label{ehiterate}
Y(h(z_{0})^{-}v,z)
&=&\left(e^{z_{0}{\partial\over\partial z}}h(z)^{-}\right)
Y(v,z)+Y(v,z)\left(e^{z_{0}{\partial\over\partial z}}h(z)^{+}\right)
\nonumber\\
&=&h(z+z_{0})^{-}Y(v,z)+Y(v,z)h(z+z_{0})^{+}.
\end{eqnarray}

Claim 1: The Jacobi identity of $Y$ holds for $(u,v,w)$ with 
$u\in \Omega_{V}^{\a},\;v\in \Omega_{V}^{\b},\;w\in U$.
With (\ref{econjugation2}) and (\ref{egeneralvertexoperator}),
it was shown in the proof of Theorem \ref{tcosetalgebra} that
the Jacobi identity of $Y$ for $(u,v,w)$ is equivalent
to the generalized Jacobi identity of $Y_{\Omega}$ for $(u,v,w)$.

Claim 2: The Jacobi identity of $Y$ holds for $(u,v,w)$ 
with $u\in \Omega_{V}^{\a},\;v\in \Omega_{V}^{\b},\;w\in E(U)$.
Assume that the Jacobi identity of $Y$ holds for $(u,v,w)$ 
for some $w\in E(U)$.
Let $h\in {\bf h},\; m\in -\Z_{+}$.
Then using (\ref{ecommproof}) we get
\begin{eqnarray}
& &[h(m), Y(u,z_{1})Y(v,z_{2})]
=(\<h,\a\>z_{1}^{m}+\<h,\b\>z_{2}^{m})Y(u,z_{1})Y(v,z_{2})
\label{e3.59}\\
& &[h(m), Y(v,z_{2})Y(u,z_{1})]
=(\<h,\a\>z_{1}^{m}+\<h,\b\>z_{2}^{m})Y(v,z_{2})Y(u,z_{1})\label{e3.60}
\end{eqnarray}
and using (\ref{egeneralvertexoperator}) and (\ref{ecommproof}) we get
\begin{eqnarray}\label{e3.61}
& &[h(m),Y(Y(u,z_{0})v,z_{2})]\nonumber\\
&=&\left[h(m), Y(E^{-}(-{1\over\ell}\a,z_{0})Y_{\Omega}(u,z_{0})v,z_{2})\right]
z_{0}^{{1\over\ell}\<\a,\b\>}\nonumber\\
&=&\left[h(m), E^{-}(-{1\over\ell}\a,z_{2}+z_{0})E^{-}({1\over\ell}\a,z_{2})
Y(Y_{\Omega}(u,z_{0})v,z_{2})
E^{+}(-{1\over\ell}\a,z_{2}+z_{0})E^{+}({1\over\ell}\a,z_{2})
\right]\nonumber\\
& &\hspace{2cm} \times z_{0}^{{1\over\ell}\<\a,\b\>}
z_{2}^{-{1\over\ell}\a(0)}(z_{2}+z_{0})^{{1\over\ell}\a(0)}
\nonumber\\
&=&(\<h,\a+\b\>z_{2}^{m}+
\<h,\a\>(z_{2}+z_{0})^{m}-\<h,\a\>z_{2}^{m})z_{0}^{{1\over\ell}\<\a,\b\>}
\nonumber\\
& &\times E^{-}(-{1\over\ell}\a,z_{2}+z_{0})E^{-}({1\over\ell}\a,z_{2})
Y(Y_{\Omega}(u,z_{0})v,z_{2})
E^{+}(-{1\over\ell}\a,z_{2}+z_{0})E^{+}({1\over\ell}\a,z_{2})
\nonumber\\
& &\times z_{2}^{-{1\over\ell}\a(0)}(z_{2}+z_{0})^{{1\over\ell}\a(0)}
\nonumber\\
&=&(\<h,\a\>(z_{2}+z_{0})^{m}+\<h,\b\>z_{2}^{m})Y(Y(u,z_{0})v,z_{2}).
\end{eqnarray}
Note that
\begin{eqnarray}
z_{2}^{-1}\delta\left(\frac{z_{1}-z_{0}}{z_{2}}\right)(z_{2}+z_{0})^{m}
=z_{2}^{-1}\delta\left(\frac{z_{1}-z_{0}}{z_{2}}\right)z_{1}^{m}.
\end{eqnarray}
Then we see that the Jacobi identity of $Y$ for $(u,v, h(m)w)$ follows from
the Jacobi identity for $(u,v,w)$.
It follows from induction that
the Jacobi identity of $Y$ holds for $(u,v,w)$ 
with $w\in E(U)$ being arbitrary.

Claim 3: The Jacobi identity holds
for $(u,v,w)$ with $u\in \Omega_{V}^{\a},\; v\in V,\;w\in E(U)$.
Assume that the Jacobi identity holds
for $(u,v,w)$ for some $v\in V$ and for all $w\in E(U)$ with $u$ being fixed.
Let $h\in {\bf h}$. Using (\ref{ehiterate}) and (\ref{eh-yu}) we get
\begin{eqnarray}
& &Y(u,z_{1})Y(h(z)^{-}v,z_{2})\nonumber\\
&=&Y(u,z_{1})h(z_{2}+z)^{-}Y(v,z_{2})+Y(u,z_{1})Y(v,z_{2})h(z_{2}+z)^{+}
\nonumber\\
&=&h(z_{2}+z)^{-}Y(u,z_{1})Y(v,z_{2})
+Y(u,z_{1})Y(v,z_{2})h(z_{2}+z)^{+}\nonumber\\
& &-\<h,\a\>(z_{1}-z_{2}-z)^{-1})Y(u,z_{1})Y(v,z_{2})
\end{eqnarray}
and using (\ref{ehiterate}) and (\ref{eh+yu}) we get
\begin{eqnarray}
& &Y(h(z)^{-}v,z_{2})Y(u,z_{1})\nonumber\\
&=&h(z_{2}+z)^{-}Y(v,z_{2})Y(u,z_{1})
+Y(v,z_{2})Y(u,z_{1})h(z_{2}+z)^{+}\nonumber\\
& &+\<h,\a\>(z_{2}+z-z_{1})^{-1}Y(v,z_{2})Y(u,z_{1}).
\end{eqnarray}
On the other hand, using (\ref{eh-yu}) on $V$ and (\ref{ehiterate}) we get
\begin{eqnarray}
& &Y(Y(u,z_{0})h(z)^{-}v,z_{2})\nonumber\\
&=&Y(h(z)^{-}Y(u,z_{0})v,z_{2})
-\<h,\a\>(z_{0}-z)^{-1}Y(Y(u,z_{0})v,z_{2})\nonumber\\
&=&h(z_{2}+z)^{-}Y(Y(u,z_{0})v,z_{2})
+Y(Y(u,z_{0})v,z_{2})h(z_{2}+z)^{+}\nonumber\\
& &-\<h,\a\>(z_{0}-z)^{-1}Y(Y(u,z_{0})v,z_{2}).
\end{eqnarray}
Note that
\begin{eqnarray}
& &z_{0}^{-1}\delta\left(\frac{z_{1}-z_{2}}{z_{0}}\right)
(z_{1}-z_{2}-z)^{-1}=z_{0}^{-1}\delta\left(\frac{z_{1}-z_{2}}{z_{0}}\right)
(z_{0}-z)^{-1},\\
& &z_{0}^{-1}\delta\left(\frac{z_{2}-z_{1}}{-z_{0}}\right)
(z_{2}+z-z_{1})^{-1}=-z_{0}^{-1}\delta\left(\frac{z_{2}-z_{1}}{-z_{0}}\right)
(z_{0}-z)^{-1}.
\end{eqnarray}
Then the Jacobi identity for $(u,h(z)^{-}v)$ on $E(U)$ follows from
the Jacobi identity for $(u,v)$ on $E(U)$.
Now, it follows from induction that the Jacobi identity
holds for $(u,v,w)$ with $u\in \Omega_{V},\; v\in V,\;w\in E(U)$.

Claim 4: The Jacobi identity of $Y$ on $E(U)$ holds for $(h,u)$ with 
$h\in {\bf h},\;u\in V$.

First, consider $u\in \Omega_{V}^{\a},\; \a\in L$.
Using (\ref{eh-yu}) and (\ref{eh+yu}) we get
\begin{eqnarray}
& &z_{0}^{-1}\delta\left(\frac{z_{1}-z_{2}}{z_{0}}\right)h(z_{1})Y(u,z_{2})
\nonumber\\
&=&z_{0}^{-1}\delta\left(\frac{z_{1}-z_{2}}{z_{0}}\right)
\left(h(z_{1})^{-}Y(u,z_{2})+Y(u,z_{2})h(z_{1})^{+}+\<h,\a\>(z_{1}-z_{2})^{-1}
Y(u,z_{2})\right)\nonumber\\
&=&z_{0}^{-1}\delta\left(\frac{z_{1}-z_{2}}{z_{0}}\right)
\left(h(z_{1})^{-}Y(u,z_{2})+Y(u,z_{2})h(z_{1})^{+}+\<h,\a\>z_{0}^{-1}
Y(u,z_{2})\right)
\end{eqnarray}
and
\begin{eqnarray}
& &z_{0}^{-1}\delta\left(\frac{z_{2}-z_{1}}{-z_{0}}\right)Y(u,z_{2})h(z_{1})
\nonumber\\
&=&z_{0}^{-1}\delta\left(\frac{z_{2}-z_{1}}{-z_{0}}\right)
\left(Y(u,z_{2})h(z_{1})^{+}+h(z_{1})^{-}Y(u,z_{2})-\<h,\a\>(z_{2}-z_{1})^{-1}
Y(u,z_{2})\right)\nonumber\\
&=&z_{0}^{-1}\delta\left(\frac{z_{2}-z_{1}}{-z_{0}}\right)
\left(h(z_{1})^{-}Y(u,z_{2})+Y(u,z_{2})h(z_{1})^{+}+\<h,\a\>z_{0}^{-1}
Y(u,z_{2})\right).
\end{eqnarray}
Then using the fundamental properties of delta function and (\ref{ehiterate})
we get
\begin{eqnarray}
& &z_{0}^{-1}\delta\left(\frac{z_{1}-z_{2}}{z_{0}}\right)h(z_{1})Y(u,z_{2})
-z_{0}^{-1}\delta\left(\frac{z_{2}-z_{1}}{-z_{0}}\right)Y(u,z_{2})h(z_{1})
\nonumber\\
&=&z_{2}^{-1}\delta\left(\frac{z_{1}-z_{0}}{z_{2}}\right)
\left(h(z_{1})^{-}Y(u,z_{2})+Y(u,z_{2})h(z_{1})^{+}
+\<h,\a\>z_{0}^{-1}Y(u,z_{2})\right)\nonumber\\
&=&z_{2}^{-1}\delta\left(\frac{z_{1}-z_{0}}{z_{2}}\right)
\left(h(z_{2}+z_{0})^{-}Y(u,z_{2})+Y(u,z_{2})h(z_{2}+z_{0})^{+}
+\<h,\a\>z_{0}^{-1}Y(u,z_{2})\right)\nonumber\\
&=&z_{2}^{-1}\delta\left(\frac{z_{1}-z_{0}}{z_{2}}\right)
Y((h(z_{0})^{-}+h(0)z_{0}^{-1})u,z_{2})\nonumber\\
&=&z_{2}^{-1}\delta\left(\frac{z_{1}-z_{0}}{z_{2}}\right)Y(Y(h,z_{0})u,z_{2}).
\end{eqnarray}

For $h_{1},h_{2}\in {\bf h}$, we have
\begin{eqnarray}\label{eh1h2-}
[h_{1}(z_{1}),h_{2}(z_{2})^{-}]
&=&\sum_{m\in \Z}\sum_{n\ge 1}[h_{1}(m),h_{2}(-n)]
z_{1}^{-m-1}z_{2}^{n-1}\nonumber\\
&=&\sum_{n\ge 1}n\ell \<h_{1},h_{2}\>
z_{1}^{-n-1}z_{2}^{n-1}\nonumber\\
&=&\ell \<h_{1},h_{2}\> (z_{1}-z_{2})^{-2}
\end{eqnarray}
and similarly,
\begin{eqnarray}\label{eh1h2+}
[h_{1}(z_{1}),h_{2}(z_{2})^{+}]=-\ell \<h_{1},h_{2}\> (z_{2}-z_{1})^{-2}.
\end{eqnarray}
Following the proof of Claim 3,  using (\ref{eh1h2-}) and (\ref{eh1h2+})
we easily see the Jacobi identity holds on $E(U)$
for $(h,v)$ with all $v\in V$.

Now, the Jacobi identity holds on $E(U)$ for $(u,v)$ with 
$u\in \Omega_{V}\cup {\bf h},\; v\in V$. 
Note that $\Omega_{V}\cup {\bf h}$ generates $V$ as 
a vertex algebra. Immediately after this theorem 
we shall prove a simple general result (Lemma \ref{lspecial}) from which
it follows that $E(U)$ is a $V$-module.

Clearly, $E(U)$ is in ${\cal{C}}$ and as an $\Omega_{V}$-module, 
$\Omega_{E(U)}=U$.

For $W\in {\cal{C}}$, let $\eta$ be the linear map from
$E(\Omega_{W})$ to $W$ defined by $\eta(a\otimes u)=a\cdot u$
for $a\in U(\hat{\bf h}^{-}),\; u\in U$. Then $\eta$ is 
an $\hat{\bf h}$-isomorphism from [LW3]. 

Furthermore, 
for $v\in \Omega_{V}^{\a},\;\a\in L,\;a\in U(\hat{\bf h}^{-1}),\; u\in 
\Omega_{W}$, we have
\begin{eqnarray}
& &\eta (Y(v,z)(a\otimes u))\nonumber\\
&=&\eta \left(E^{-}(-{1\over\ell}\a,z)E^{+}(-{1\over\ell}\a,z)a
\otimes Y_{\Omega}(v,z)z^{{1\over\ell}\a(0)}u\right)\nonumber\\
&=&E^{-}(-{1\over\ell}\a,z)E^{+}(-{1\over\ell}\a,z)
\eta (a\otimes Y_{\Omega}(v,z)z^{{1\over\ell}\a(0)}u)\nonumber\\
&=&E^{-}(-{1\over\ell}\a,z)E^{+}(-{1\over\ell}\a,z)
(a\cdot Y_{\Omega}(v,z)z^{{1\over\ell}\a(0)}u)\nonumber\\
&=&E^{-}(-{1\over\ell}\a,z)E^{+}(-{1\over\ell}\a,z)
Y_{\Omega}(v,z)z^{{1\over\ell}\a(0)}(a\cdot u)\nonumber\\
&=&Y(v,z)\eta(a\otimes u).
\end{eqnarray}
Since ${\bf h}$ and $\Omega_{V}^{\a},\;\a\in L$ generate $V$ 
as a vertex algebra, $\eta$ must be a $V$-homomorphism.
Thus $\eta$ is a $V$-isomorphism from $E(\Omega_{W})$ to $W$. 
$\;\;\;\;\Box$

\bl{lspecial}
Let $V$ be a vertex algebra for now and $(W,Y)$ be a pair
satisfying all the axioms defining the notion of a $V$-module
except the Jacobi identity. In addition, assume that 
the Jacobi identity holds for $(Y_{W};a,v)$ with $a\in A,\; v\in V$, where
$A$ is a set of generators of $V$ as a vertex algebra. Then
$(W,Y_{W})$ is a $V$-module.
\el

\pf One can in principle use induction to prove that
the Jacobi identity holds for any pair of elements $u,v$ of $V$.
Here we use a result of [Li2] to prove this.

Set
\begin{eqnarray}
\tilde{A}=\{ Y_{W}(a,z)\in (\End W)[[z,z^{-1}]]\;|\; a\in A\}.
\end{eqnarray}
For $a\in A,\; v\in V$,
due to the Jacobi identity of $Y_{W}$ for $(a,v)$,
$Y_{W}(a,z)$ and $Y_{W}(v,z)$ are mutually local, i.e.,
they satisfy the generalized weak commutativity with
$c(\a,\b)=1$ and $(\a,\b)=0$.
In particular, $\tilde{A}$ is a local set of 
vertex operators on $W$.
{}From [Li2], $\tilde{A}$ generates a canonical vertex algebra $\tilde{V}$
inside $(\End W)[[z,z^{-1}]]$ with $W$ as a natural module.

For $a\in A,\; u\in V$, using the Jacobi identity of $Y_{W}$
for $(a,u)$ and the definition of the vertex operator map of $\tilde{V}$
[Li2] we have
\begin{eqnarray}\label{eiterateformula}
Y_{W}(Y(a,z_{0})u,z)
&=&\Res_{z_{1}}z_{0}^{-1}\delta\left(\frac{z_{1}-z}{z_{0}}\right)
Y_{W}(a,z_{1})Y_{W}(u,z)\nonumber\\
& &-\Res_{z_{1}}z_{0}^{-1}\delta\left(\frac{z-z_{1}}{-z_{0}}\right)
Y_{W}(u,z)Y_{W}(a,z_{1})\nonumber\\
&=&Y_{\tilde{V}}(Y_{W}(a,z),z_{0})Y_{W}(u,z).
\end{eqnarray}
Since $A$ generates $V$ and $Y_{W}(a,z)\in \tilde{V}$ for $a\in A$,
it follows from (\ref{eiterateformula}) and induction that
$Y_{W}(u,z)\in \tilde{V}$ for $u\in V$.
Then we have a linear map from $V$ to $\tilde{V}$:
\begin{eqnarray}
\rho : v\mapsto Y_{W}(v,z)\;\;\;\mbox{ for }v\in V.
\end{eqnarray}
In view of (\ref{eiterateformula}), we have
\begin{eqnarray}
\rho (Y(a,z_{0})v)=Y(\rho(a),z_{0})\rho (v)\;\;\;
\mbox{ for }a\in A,\; v\in V.
\end{eqnarray}
Since $A$ generates $V$, it follows from induction that
\begin{eqnarray}
\rho (Y(u,z_{0})v)=Y(\rho(u),z_{0})\rho (v)\;\;\;\mbox{ for }u,v\in V.
\end{eqnarray}
Thus $\rho$ is a homomorphism of vertex algebras, noting that
$\rho({\bf 1})=1$ follows from the assumption.
Consequently, $W$ is a $V$-module.
$\;\;\;\;\Box$

With Theorem \ref {tcosetmodule} and Proposition \ref{pback} 
we immediately have (cf. [LW5], Theorem 5.5):

\bt{tequivalencecategory}
The functors $\Omega$ from ${\cal{C}}$ to ${\cal{D}}$ and 
$E$ from ${\cal{D}}$ to ${\cal{C}}$ define equivalences of categories.
$\;\;\;\;\Box$
\et

Note that an $\Omega_{V}$-${\bf h}$-module amounts to an
$\Omega_{V}$-module associated with a {\em free} $L$-set
$S$ and with a {\em $\C$-valued} function on $L\times S$.
Since different $\C$-valued functions on $L\times S$
may give rise to the same $\C/2\Z$-valued (or $\C/\Z$-valued) 
function on $L\times S$, nonisomorphic $V$-modules may give rise to
isomorphic $\Omega_{V}$-modules. In the following, we shall give
a criterion to determine when two $V$-modules give rise to 
isomorphic $\Omega_{V}$-modules.

Define
\begin{eqnarray}
L^{o}=\{\b\in {\bf h}\;|\; \<\a,\b\>\in \Z\;\;\;
\mbox{ for }\a\in L\}.
\end{eqnarray}
Let $\b\in L^{o}$ and let $W$ be a weak $V$-module. Recall from [Li3]
(see also [Li6]) the weak $V$-module $W^{(\b)}$:
\begin{eqnarray}\label{edvbeta1}
W^{(\b)}=\C e^{\b}\otimes W\;\;\mbox{ as a vector space},
\end{eqnarray}
where $\C e^{\b}$ is a one-dimensional vector space with a distinguished
basis element $e^{\b}$, with action
\begin{eqnarray}\label{edeformation}
Y(v,z)(e^{\b}\otimes w)=e^{\b}\otimes Y(E^{+}(-\b,-z)z^{\b(0)}v,z)w
\;\;\;\mbox{ for }v\in V,\; w\in W.
\end{eqnarray}
Furthermore, for weak $V$-modules $W_{1},W_{2}$, 
\begin{eqnarray}
e^{\b}\otimes f\in \Hom_{V}(W_{1}^{(\b)},W_{2}^{(\b)})
\;\;\;\mbox{ for }f\in \Hom_{V}(W_{1},W_{2}).
\end{eqnarray}
Thus, we have a functor $F_{\b}$ of the category of weak $V$-modules
defined in the obvious way.
The following are some straightforward consequences:

\bl{lsimplecurrents}
Let $\b\in L^{o}$ and let $W$ be a weak $V$-module. Then
the map $e^{\b}\otimes\cdot$ gives rise to a one-to-one
correspondence between the set of $V$-submodules of $W$
and the set of $V$-submodules of $W^{(\b)}$.
In particular, $W^{(\b)}$ is irreducible if and only if
$W$ is irreducible. 
Furthermore, as weak $V$-modules,
\begin{eqnarray}
& &W^{(0)}\simeq W,\label{e3.98}\\
& &(W^{(\b_{1})})^{(\b_{2})}\simeq W^{(\b_{1}+\b_{2})}
\;\;\;\mbox{ for }\b_{1},\b_{2}\in L^{o},\label{eassociativitydef}
\end{eqnarray}
where the map 
\begin{eqnarray}\label{eassociativitydefisom}
e^{\b_{2}}\otimes (e^{\b_{1}}\otimes 
w)\mapsto e^{\b_{1}+\b_{2}}\otimes w
\end{eqnarray}
is a $V$-isomorphism from $(W^{(\b_{1})})^{(\b_{2})}$ to $W^{(\b_{1}+\b_{2})}$.
$\;\;\;\;\Box$
\el

Next we show that each functor $F_{\b}$ preserves the subcategory ${\cal{C}}$.
%In fact, $F_{\b}$ is an isomorphism [Li6].

\bp{pidenfitication}
Let $\b\in L^{o},\; W\in {\cal{C}}$. 
Then $W^{(\b)}\in {\cal{C}}$ and the linear map 
\begin{eqnarray}
e^{\b}\otimes \cdot :\;\; W\rightarrow W^{(\b)}=\C e^{\b}\otimes W; \;\;\;\;
 w\mapsto e^{\b}\otimes w
\end{eqnarray}
is an $\Omega_{V}$-isomorphism such that
\begin{eqnarray}
& &\C e^{\b}\otimes \Omega_{W}=\Omega_{W^{(\b)}},\label{etwoomegas}\\
& &e^{\b}\otimes \Omega_{W}^{s}=\Omega_{W^{(\b)}}^{s+\ell\b}\;\;\;
\mbox{ for }s\in P(W).\label{etwoomegas1}
\end{eqnarray}
\ep

\pf Note that for $h\in {\bf h}$,
\begin{eqnarray}
E^{+}(-\b,-z)z^{\b(0)}h=h+\ell \<\b,h\>{\bf 1}z^{-1}.
\end{eqnarray}
Then for $w\in W$,
\begin{eqnarray}\label{e3.97}
Y(h,z)(e^{\b}\otimes w)=e^{\b}\otimes (Y(h,z)+\ell \<\b,h\>z^{-1})w.
\end{eqnarray}
In terms of components, we have
\begin{eqnarray}\label{ehpsicomm}
h(n)(e^{\b}\otimes w)=e^{\b}\otimes (h(n)+\ell \<\b,h\>\delta_{n,0})w 
\;\;\;\mbox{ for }n\in \Z.
\end{eqnarray}
Then it follows immediately that $W^{(\b)}\in {\cal{C}}$ and 
(\ref{etwoomegas}), (\ref{etwoomegas1}) hold.
%We also have
%\begin{eqnarray}
%P(W^{(\b)})=P(W)+\ell \b.
%\end{eqnarray}

For $v\in \Omega_{V}^{\a},\;\a\in L$, we have
\begin{eqnarray}
E^{+}(-\b,-z)z^{\b(0)}v=z^{\<\b,\a\>}v=z^{\<\a,\b\>}v.
\end{eqnarray}
Then for $w\in W$, using (\ref{ehpsicomm}) 
and the definition (\ref{edeformation}) we get
\begin{eqnarray}\label{e3.64}
& &e^{\b}\otimes Y_{\Omega}(v,z)w\nonumber\\
&=&e^{\b}\otimes E^{-}({1\over\ell}\a,z)Y(v,z)
E^{+}({1\over\ell}\a,z)z^{-{1\over \ell}\a(0)}w
\nonumber\\
&=&E^{-}({1\over\ell}\a,z)(e^{\b}\otimes 
Y(v,z)E^{+}({1\over\ell}\a,z)z^{-{1\over \ell}\a(0)}w)
\nonumber\\
&=&E^{-}({1\over\ell}\a,z)Y(v,z)(e^{\b}
\otimes E^{+}({1\over\ell}\a,z)z^{-\<\a,\b\>-{1\over \ell}\a(0)}w)
\nonumber\\
&=&E^{-}({1\over\ell}\a,z)Y(v,z)E^{+}({1\over\ell}\a,z)
z^{-{1\over \ell}\a(0)}(e^{\b}\otimes w)
\nonumber\\
&=&Y_{\Omega}(v,z)(e^{\b}\otimes w).
\end{eqnarray}
This proves that $e^{\b}\otimes \cdot$ is an 
$\Omega_{V}$-isomorphism from $W$ to $W^{(\b)}$. $\;\;\;\;\Box$

Let $W_{1}, W_{2}\in {\cal{C}}$ be weak $V$-modules
and let $f\in \Hom_{V}(W_{1}^{(\b)},W_{2})$ for some
$\b\in L^{o}$. It follows from Proposition \ref{pidenfitication} that
the restriction on $\Omega_{W_{1}}$ of the map 
$ f\circ (e^{\b}\otimes\cdot)$ is an $\Omega_{V}$-homomorphism.
Furthermore, we have:

\bp{pisomorphismrelation}
Let $W_{1}, W_{2}\in {\cal{C}}$ be weak $V$-modules
and let $\b\in L^{o}$. 
Denote by $\Hom_{\Omega_{V}}^{\b}(\Omega_{W_{1}},\Omega_{W_{2}})$
the space of $\Omega_{V}$-homomorphisms $f$ from $\Omega_{W_{1}}$
to $\Omega_{W_{2}}$ such that
\begin{eqnarray}
f(\Omega_{W_{1}}^{\lambda})\subset \Omega_{W_{2}}^{\lambda+\ell \b}
\;\;\;\mbox{ for }\lambda\in P(W_{1}).
\end{eqnarray}
Then the map
\begin{eqnarray}\label{emap}
\Hom_{V}(W_{1}^{(\b)},W_{2})&\rightarrow&
\Hom_{\Omega_{V}}^{\b}(\Omega_{W_{1}},\Omega_{W_{2}})\nonumber\\
&f\mapsto& f\circ (e^{\b}\otimes\cdot)_{\Omega},
\end{eqnarray}
where $(e^{\b}\otimes\cdot)_{\Omega}$ 
denotes the restriction of $(e^{\b}\otimes\cdot)$ onto $\Omega_{W_{1}}$,
is a linear isomorphism which sends $V$-isomorphisms to 
$\Omega_{V}$-isomorphisms.
\ep

\pf It is easy to see that the linear map (\ref{emap}) is injective.

On the other hand, let 
$f\in \Hom_{\Omega_{V}}^{\b}(\Omega_{W_{1}},\Omega_{W_{2}})$.
Define a linear map
\begin{eqnarray}
\tilde{f}: \;W_{1}^{(\b)}=\C e^{\b}\otimes M(\ell)\otimes \Omega_{W_{1}}
&\rightarrow& W_{2}=M(\ell)\otimes \Omega_{W_{2}}\nonumber\\
e^{\b}\otimes a\otimes w&\mapsto& a\otimes f(w)
\end{eqnarray}
for $a\in M(\ell),\; w\in \Omega_{W_{1}}$. 

For $u\in \Omega_{V}^{\a},\; \a\in L,\; a\in M(\ell),\;
w\in \Omega_{W_{1}}^{s},\; s\in P(W_{1})$, we have
\begin{eqnarray}
& &Y(u,z)\tilde{f}(e^{\b}\otimes a\otimes w)\nonumber\\
&=&z^{{1\over\ell}\<\a,s+\ell\b\>}
E^{-}(-{1\over\ell}\a,z)E^{+}(-{1\over\ell}\a,z)a
\otimes Y_{\Omega}(u,z)f(w)\nonumber\\
&=&z^{{1\over\ell}\<\a,s+\ell\b\>}
E^{-}(-{1\over\ell}\a,z)E^{+}(-{1\over\ell}\a,z)a
\otimes f(Y_{\Omega}(u,z)w)\nonumber\\
&=&z^{{1\over \ell}\<\a,s+\ell\b\>}
\tilde{f}\left(e^{\b}\otimes E^{-}(-{1\over\ell}\a,z)E^{+}(-{1\over\ell}\a,z)a
\otimes Y_{\Omega}(u,z)w\right)\nonumber\\
&=&\tilde{f}(e^{\b}\otimes z^{\<\a,\b\>}Y(u,z)(a\otimes w))\nonumber\\
&=&\tilde{f}(e^{\b}\otimes Y(E^{+}(-\b,-z)z^{\b(0)}u,z)(a\otimes w))
\nonumber\\
&=&\tilde{f}(Y(u,z)(e^{\b}\otimes a\otimes w)).
\end{eqnarray}
On the other hand, for $h\in {\bf h}$, using (\ref{e3.97}) we get
\begin{eqnarray}
& &Y(h,z)\tilde{f}(e^{\b}\otimes a\otimes w)\nonumber\\
&=&Y(h,z)(a\otimes f(w))\nonumber\\
&=&(Y(h,z)-h(0)z^{-1})a\otimes f(w)+a\otimes h(0)f(w)z^{-1}\nonumber\\
&=&(Y(h,z)-h(0)z^{-1})a\otimes f(w)+\ell\<h,\b\> z^{-1})(a\otimes f(w))
\nonumber\\
&=&\tilde{f}(e^{\b}\otimes (Y(h,z)-h(0)z^{-1})a\otimes w)
+\<h,s+\ell\b\> \tilde{f}(e^{\b}\otimes a\otimes w)z^{-1}\nonumber\\
&=&\tilde{f}(e^{\b}\otimes (Y(h,z)+\ell\<h,\b\> z^{-1})(a\otimes w))
\nonumber\\
&=&\tilde{f}(e^{\b}\otimes Y(E^{+}(-\b,-z)z^{\b (0)}h,z)(a\otimes w))
\nonumber\\
&=&\tilde{f}(Y(h,z)(e^{\b}\otimes a\otimes w)).
\end{eqnarray}
That is,  $\tilde{f}$ is a linear map from $W_{1}^{(\b)}$
to $W_{2}$ such that
\begin{eqnarray}\label{elast1}
Y(v,z)\tilde{f}(e^{\b}\otimes w)=\tilde{f}(Y(v,z)(e^{\b}\otimes w))
\end{eqnarray}
for $v\in \Omega_{V}\cup {\bf h},\; w\in W_{1}$.
Since ${\bf h}$ and $\Omega_{V}$ generate $V$ as a vertex algebra, 
it follows from induction and the Jacobi identity of the vertex algebra $V$
that (\ref{elast1}) holds for all $v\in V,\; w\in W_{1}$. That is,
$\tilde{f}$ is a $V$-homomorphism. Clearly, $f$ is the restriction of
$\tilde{f}\circ (e^{\b}\otimes\cdot)$. This proves 
that the linear map (\ref{emap}) is also onto.
Clearly, $\tilde{f}$ is a linear isomorphism if and only if 
$f$ is a linear isomorphism. This completes the proof.
$\;\;\;\Box$

Now, we present our second main result of this section.

\bt{tequivalence}
Let $W_{1}, W_{2}$ be irreducible weak $V$-modules in ${\cal{C}}$.
Then $\Omega_{W_{1}}$ and $\Omega_{W_{2}}$ are isomorphic $\Omega_{V}$-modules
if and only if $W_{2}\simeq W_{1}^{(\b)}$ for some $\b\in L^{o}$.
\et

\pf The ``if'' part follows immediately from 
Proposition \ref{pidenfitication}.
Now, suppose that $\Omega_{W_{1}}$ and $\Omega_{W_{2}}$ are 
isomorphic $\Omega_{V}$-modules.
Let $f$ be an $\Omega_{V}$-isomorphism from $\Omega_{W_{1}}$ onto 
$\Omega_{W_{2}}$ with an $L$-set map
$\bar{f}$  from $P(W_{1})$ to $P(W_{2})$.
Let $\lambda\in {\bf h}$ with $W_{1}^{\lambda}\ne 0$. 
Since $W_{1},W_{2}$ are irreducible, 
$P(W_{1})=\lambda+L$ and
$P(W_{2})=\bar{f}(\lambda)+L$.

Let $0\ne u\in \Omega_{V}^{\a}$ for $\a\in L$ and 
$0\ne w\in \Omega_{W_{1}}^{\lambda}$.
Then
\begin{eqnarray}
z^{{1\over\ell }\<\a,\lambda\>}Y_{\Omega}(u,z)w\in W_{1}[[z,z^{-1}]].
\end{eqnarray}
With $f$ being an $\Omega_{V}$-isomorphism, we have
\begin{eqnarray}
z^{{1\over\ell }\<\a,\lambda\>}Y_{\Omega}(u,z)f(w)\in W_{2}[[z,z^{-1}]].
\end{eqnarray}
On the other hand, with $f(w)\in \Omega_{W_{2}}^{\bar{f}(\lambda)}$ we have
\begin{eqnarray}
z^{{1\over\ell }\<\a,\bar{f}(\lambda)\>}Y_{\Omega}(u,z)f(w)
\in W_{2}[[z,z^{-1}]].
\end{eqnarray}
Since $W_{2}$ is irreducible, by Proposition 11.9 of [DL2],
$Y(u,z)f(w)\ne 0$. Then we have $Y_{\Omega}(u,z)f(w)\ne 0$. 
Consequently,
\begin{eqnarray}\label{e3.118}
{1\over\ell}\<\a,\bar{f}(\lambda)-\lambda\>\in \Z.
\end{eqnarray}
Since $L$ by definition is generated by $\a$ with $\Omega_{V}^{\a}\ne 0$,
(\ref{e3.118}) holds for all $\a\in L$. Thus
${1\over\ell}(\bar{f}(\lambda)-\lambda)\in L^{o}$.
With $\bar{f}$ being an $L$-set map, clearly
${1\over\ell}(\bar{f}(\lambda)-\lambda)$ does not depend on $\lambda$.
Set
\begin{eqnarray}
\b={1\over\ell}(\bar{f}(\lambda)-\lambda)\in L^{o}.
\end{eqnarray}
Then
\begin{eqnarray}
\bar{f}(s)=s+\ell \b\;\;\;\mbox{ for }s\in P(W_{1}).
\end{eqnarray}
Thus $f\in \Hom_{\Omega_{V}}^{\b}(\Omega_{W_{1}},\Omega_{W_{2}})$.
It follows immediately from Proposition \ref{pisomorphismrelation}
that $\tilde{f}$ is a $V$-isomorphism from
$W_{1}^{(\b)}$ to $W_{2}$. The proof is complete. $\;\;\;\;\Box$

\br{rsimple}
{\em In view of Lemma \ref{lsimplecurrents} and 
Proposition \ref{pidenfitication}, the abelian group 
$L^{o}$ naturally acts on the set $Irr({\cal{C}})$ of equivalence classes 
of irreducible weak $V$-modules in ${\cal{C}}$.  
Then Theorem \ref{tequivalence} states that
for irreducible weak $V$-modules $W_{1},W_{2}\in {\cal{C}}$, 
$\Omega_{W_{1}}\simeq \Omega_{W_{2}}$ as an $\Omega_{V}$-module if and only if
the equivalence classes $[W_{1}]$ and $[W_{2}]$ are in the same $L^{o}$-orbit.}
\er

We conclude this section with the special case where $\ell$ is not rational.

\bl{lirrational}
Suppose that $L$ equipped with $\<\cdot,\cdot\>$ is a
nondegenerate rational lattice and that $\ell$ is not rational. 
(1) Let $S$ be an $L$-set
equipped with a $\C/\Z$-valued function $(\cdot,\cdot)$ 
on $L\times S$ such that in $\C/\Z$,
\begin{eqnarray}
(g_{1}+g_{2},g_{3}+s)={1\over\ell}\<g_{1},g_{3}\>+
{1\over\ell}\<g_{2},g_{3}\>+(g_{1},s)+(g_{2},s)
\end{eqnarray}
for $g_{1},g_{2},g_{3}\in L,\;s\in S$. 
Then $S$ is a free $L$-set.  (2) The $\C/2\Z$-valued
$\Z$-bilinear form $(\cdot,\cdot)$ on $L$ defined by
\begin{eqnarray}
(\a,\b)={1\over\ell}\<\a,\b\>+2\Z\;\;\;\mbox{ for }\a,\b\in L
\end{eqnarray}
is nondegenerate.
\el

\pf (1) If $S$ is not free, there exist 
$\a_{1},\a_{2}\in L,\; s\in S$ 
such that $\a_{1}\ne \a_{2}$ and $\a_{1}+s=\a_{2}+s$.
Then for any $\a\in L$ from 
$(\a,\a_{1}+s)=(\a,\a_{2}+ s)$ we get
\begin{eqnarray}
{1\over \ell}\<\a,\a_{1}-\a_{2}\>\in \Z.
\end{eqnarray}
Since $\<\a,\a_{1}-\a_{2}\>$ is rational and 
$\ell$ is not rational, $\<\a,\a_{1}-\a_{2}\>$ must be zero.
This is impossible because $\<\cdot,\cdot\>$ by assumption
is nondegenerate on $L$. 

(2) Let $\b\in L$ such that $(\a,\b)=0$ in $\C/2\Z$ for all $\a\in L$. Then
${1\over\ell}\<\a,\b\>\in 2\Z$ for all $\a\in L$.
With $\ell$ being not rational and $\<\cdot,\cdot\>$ being rational,
we have $\<\a,\b\>=0$. Since $\<\cdot,\cdot\>$ is nondegenerate, $\b=0$.
Thus $(\cdot,\cdot)$ is nondegenerate.
$\;\;\;\;\Box$

\bp{pirrational}
Suppose that $L$ equipped with $\<\cdot,\cdot\>$ is a
nondegenerate rational lattice and that $\ell$ is not rational. Then
for any irreducible weak $V$-module $W\in {\cal{C}}$, $\Omega_{W}$ is an
irreducible $\Omega_{V}$-module. Furthermore, the map $W\mapsto \Omega_{W}$
gives rise to a one-to-one map between the set of $L^{o}$-orbits of
the set of equivalence classes of irreducible weak $V$-modules in ${\cal{C}}$ 
and the set of equivalence classes of irreducible $\Omega_{V}$-modules.
\ep

\pf Let $W\in {\cal{C}}$ be an irreducible weak $V$-module.
Then $P(W)=\lambda+L$ for any $\lambda\in {\bf h}$ with $W^{\lambda}\ne 0$.
It follows from the definition of a submodule and
Lemma \ref{lirrational} that
any submodule of $\Omega_{W}$ must be $P(W)$-graded, hence
an $\Omega_{V}$-${\bf h}$-submodule.
Then $\Omega_{W}$ is an irreducible $\Omega_{V}$-module
because $\Omega_{W}$ is an irreducible $\Omega_{V}$-${\bf h}$-module
(Corollary \ref{conetoone}).

On the other hand, let $(U,Y_{\Omega},S, (\cdot,\cdot))$ be 
an irreducible $\Omega_{V}$-module.
Then by Lemma \ref{lirrational} $S$ is a free $L$-set.
Since $U$ is irreducible, $S=L+s$ for any $s\in S$ with $U^{s}\ne 0$.
Consequently, $S$ is a transitive $L$-set.
Let $\{\a_{1},\dots,\a_{n}\}$ be a basis of $L$ and fix an 
element $s_{0}$ of $S$.
Choose a representative $(\a_{i},s_{0})$ in $\C$ for each $i$.
Define a $\C$-valued function $f$ on $L\times S$ by
\begin{eqnarray}
f(\a,\b+s_{0})
={1\over\ell}\<\a,\b\>+\sum_{i=1}^{n}m_{i}(\a_{i},s_{0})
\end{eqnarray}
for $\a=\sum_{i=1}^{n}m_{i}\a_{i},\;\b\in L$.
Set
\begin{eqnarray}
{\bf h}'=\C L,\;\;\; {\bf h}''=({\bf h}')^{\perp}.
\end{eqnarray}
Since $\<\cdot,\cdot\>$ is nondegenerate on $L$, we have
\begin{eqnarray}
{\bf h}={\bf h}'\oplus {\bf h}''.
\end{eqnarray}
We then lift $f(\cdot,\cdot)$ to a $\C$-valued function 
on ${\bf h}\times S$ which is $\C$-linear in the first variable
with $f({\bf h}'',S)=0$. We may view each $s\in S$ as an element of ${\bf h}$
by
\begin{eqnarray}
\<h,s\>=\ell f(h,s)\;\;\;\mbox{ for }h\in {\bf h}.
\end{eqnarray}
Using the fact that $S=L+s_{0}$ and that
$\<\cdot,\cdot\>$ is nondegenerate on ${\bf h}'$ we easily see that
$S$ is embedded as a subset of ${\bf h}$. Furthermore, we have
\begin{eqnarray}
{1\over\ell}\<\a,s\>+\Z=f(\a,s)+\Z=(\a,s)\in \C/\Z\;\;\;
\mbox{ for }\a\in L,\; s\in S.
\end{eqnarray}
We make $U$ an ${\bf h}$-module by defining
\begin{eqnarray}
hw=\<h,s\>w\;\;\;\mbox{ for }h\in {\bf h},\; w\in U^{s},\; s\in S.
\end{eqnarray}
It is clear that $U$ becomes an $\Omega_{V}$-${\bf h}$-module.
By Proposition  \ref{pback}, $E(U)\in {\cal{C}}$ with $\Omega_{E(U)}=U$.
It follows from Corollary \ref{conetoone} that
$E(U)$ is $V$-irreducible.
Now the second assertion follows from Theorem \ref{tequivalence}.
$\;\;\;\;\Box$

\section{Automorphism group $\Aut_{\Omega_{V}}\Omega_{V}$}
In this section, we shall determine the automorphism group
$\Aut_{\Omega_{V}}\Omega_{V}$ of the adjoint $\Omega_{V}$-module
for the purpose of constructing quotient generalized 
vertex (operator) algebras $\Omega_{V}^{A}$ in Section 5.
We prove that $\Aut_{\Omega_{V}}\Omega_{V}$ is a central extension of 
a free group $K$ (a subgroup of the additive group ${\bf h}$) by
$\C^{\times}$ and we furthermore determine the commutator map.

Set
\begin{eqnarray}
K=\{ \a\in L^{o}\;|\; V^{(\a)}\simeq V \;\;\;\mbox{ as a $V$-module}\}.
\end{eqnarray}
It follows from (\ref{e3.98}) and (\ref{eassociativitydef}) that
$K$ is a subgroup of $L^{o}$ (cf. [DLM1]). 

\bl{lfapsia}
Let $\a\in K$ and let $\pi_{\a}$ be a $V$-isomorphism from $V^{(\a)}$
to $V$. Set 
\begin{eqnarray}
\psi_{\a}=\pi_{\a}\circ (e^{\a}\otimes \cdot)\in \End V,
\end{eqnarray}
i.e., 
\begin{eqnarray}
\psi_{\a}(v)=\pi_{\a}(e^{\a}\otimes v)\;\;\;\mbox{ for }v\in V.
\end{eqnarray}
Then $\psi_{\a}$ is an $\Omega_{V}$-automorphism 
of $V$ such that
\begin{eqnarray}\label{ebrackethnpsia}
[h(n),\psi_{\a}]=\delta_{n,0}\ell\<\a,h\> \psi_{\a}\;\;\;
\mbox{ for }h\in {\bf h},\; n\in \Z,
\end{eqnarray}
and such that $\psi_{\a}(\Omega_{V})=\Omega_{V}$ and
\begin{eqnarray}\label{epsiagrading}
\psi_{\a}(\Omega_{V}^{g})=\Omega_{V}^{g+\ell\a}\;\;\;
\mbox{ for }g\in L.
\end{eqnarray}
Conversely,
let $\phi$ be an $\Omega_{V}$-automorphism of $\Omega_{V}$. 
Then there exist a unique $\a\in K$ and a unique $V$-isomorphism 
$\tilde{\phi}$ from
$V^{(\a)}$ to $V$ such that $\phi=\tilde{\phi}\circ (e^{\a}\otimes \cdot)$.
\el

\pf Clearly, $\pi_{\a}$ is also an $\Omega_{V}$-isomorphism from $V^{(\a)}$
to $V$. With $e^{\a}\otimes\cdot$ being an $\Omega_{V}$-isomorphism
(Proposition \ref{pidenfitication}),
$\psi_{\a}$ is an $\Omega_{V}$-automorphism of $V$. 
The relation (\ref{ebrackethnpsia}) follows immediately from 
(\ref{ehpsicomm}). Consequently,
$\psi_{\a}(\Omega_{V})=\Omega_{V}$ and (\ref{epsiagrading}) holds.
This proves the first part.

For the converse, from the definition of $\Omega_{V}$-isomorphism,
associated to $\phi$ there is an $L$-set map $\bar{\phi}$ from $L$ to $L$
such that 
\begin{eqnarray}
& &{1\over\ell}\<g,g'\>+\Z={1\over\ell}\<g,\bar{\phi}(g')\>+\Z,
\label{elsetmap4.6}\\
& & \phi(\Omega_{V}^{g})=\Omega_{V}^{\bar{\phi}(g)}
\;\;\;\mbox{ for }g,g'\in L.
\end{eqnarray}
Set $\b=\bar{\phi}(0)\in L$. Then $\b\in \ell L^{o}$ 
(from (\ref{elsetmap4.6}) with $g'=0$) and
\begin{eqnarray}
\phi(\Omega_{V}^{g})=\Omega_{V}^{g+\b}\;\;\;\mbox{ for }g\in L.
\end{eqnarray}
By Proposition \ref{pisomorphismrelation}, there is 
a $V$-isomorphism $\tilde{\phi}$ from 
$V^{({1\over\ell}\b)}$ to $V$ such that 
$\phi=\tilde{\phi}\circ (e^{{1\over\ell}\b}\otimes\cdot)$
Consequently, ${1\over\ell}\b\in K$. Clearly, ${1\over\ell}\b$ is uniquely
determined by $\phi$. The uniqueness of $\tilde{\phi}$ 
follows from Proposition \ref{pisomorphismrelation}.
This completes the proof.$\;\;\;\;\Box$

In view of Lemma \ref{lfapsia} we have a map $\gamma$ from 
$\Aut_{\Omega_{V}}\Omega_{V}$ onto $K$ such that
\begin{eqnarray}
\sigma (\Omega_{V}^{g})=\Omega_{V}^{g+\ell \gamma(\sigma)}
\;\;\;\mbox{ for }g\in L.
\end{eqnarray}
It is clear that $\gamma$ is a group homomorphism.
Then we have a short sequence of groups:
\begin{eqnarray}\label{eshortsequence1}
1\rightarrow \C^{\times}\rightarrow \Aut_{\Omega_{V}}\Omega_{V}
\rightarrow K\rightarrow 1
\end{eqnarray}
with $\C^{\times}$ being naturally embedded into $\Aut_{\Omega_{V}}\Omega_{V}$.

\bp{pautostructure}
If $V$ is a simple vertex operator algebra, the short sequence 
(\ref{eshortsequence1}) is exact. In particular,
$\Aut_{\Omega_{V}}\Omega_{V}$ is a central extension of $K$ 
by $\C^{\times}$. Furthermore, any collection of $V$-isomorphisms
$\pi_{\a}$ from $V^{(\a)}$ to $V$ for $\a\in K$ gives rise to a section
$\psi$ as defined in Lemma \ref{lfapsia}.
\ep

\pf We only need to show $\Ker \gamma= \C^{\times}$.
{}From definition we have
\begin{eqnarray}
\Ker \gamma = \Aut_{\Omega_{V}}^{0}\Omega_{V},
\end{eqnarray}
which consists of grading-preserving automorphisms.
Since $V$ is simple, it follows from Schur lemma that 
$$\Hom_{V}(V^{(0)},V)={\rm End}_{V}(V)=\C.$$
Then by Proposition \ref{pisomorphismrelation},
\begin{eqnarray}
\Aut_{\Omega_{V}}^{0}\Omega_{V}=\C^{\times}.
\end{eqnarray}
This completes the proof.$\;\;\;\;\Box$

The following is an immediate consequence of (\ref{epsiagrading})
and the definition of $L$:

\bc{crelationAQ}
We have
\begin{eqnarray}
\ell K\subset L.\;\;\;\;\Box
\end{eqnarray}
\ec

\br{rKzero}
{\em Suppose that $L$ equipped with the form 
$\<\cdot,\cdot\>$ is a nondegenerate rational lattice 
and that $\ell$ is a non-rational complex number. 
We now show that $K=0$.
Let $\a\in \ell K,\; \b\in L$.
By Corollary \ref{crelationAQ} we have
$\a\in \ell K\subset L$. Then
$\<\a,\b\>$ is rational.
On the other hand, with $\b\in L$ and
$\a\in \ell K\subset \ell L^{o}$, we have
$\<\a,\b\>\in \ell \Z$. 
Because $\ell$ is not rational, we must have $\<\a,\b\>=0$. 
Since by assumption $\<\cdot,\cdot\>$ is nondegenerate on $L$,
$\a=0$. Then $\ell K=0$. That is, $K=0$.}
\er

{\bf Basic Assumption 2:} From now on we assume that 
$V$ is simple, $\ell$ is rational and that $L$ 
equipped with $\<\cdot,\cdot\>$
is a nondegenerate rational lattice of finite rank.

We shall explicitly determine $\Aut_{\Omega_{V}}\Omega_{V}$ as 
an extension of $K$ by determining the commutator map of the extension 
(\ref{eshortsequence1}) (cf. [FLM], Chapter 5).
We shall first relate $\psi_{\a}$
with the vertex operator $Y_{\Omega}(\psi_{\a}({\bf 1}),z)$, and
then use the generalized weak commutativity to determine 
the commutator map.

\bp{pvacuumlike}
Let $\psi$ be a section of the extension (\ref{eshortsequence1})
obtained through Proposition \ref{pautostructure}.
Let $\a\in K$ and $W\in {\cal{C}}$. Then
\begin{eqnarray}
\<\a,s\>\in \Z\;\;\;\mbox{ for }s\in P(W)
\end{eqnarray}
and
\begin{eqnarray}
Y_{\Omega}(\psi_{\a}({\bf 1}),z)\in \End W,
\end{eqnarray}
i.e., $Y_{\Omega}(\psi_{\a}({\bf 1}),z)$ is independent of $z$.
Furthermore, 
\begin{eqnarray}
& &Y_{\Omega}(v,z)Y_{\Omega}(\psi_{\a}({\bf 1}),0)w
=(-1)^{\<g,\a\>}Y_{\Omega}(\psi_{\a}({\bf 1}),0)Y_{\Omega}(v,z)w
\label{ealternating}\\
& &Y_{\Omega}(v,z)Y_{\Omega}(\psi_{\a}({\bf 1}),0)w
=Y_{\Omega}(\psi_{\a}(v),z)w
\;\;\;\mbox{ for }v\in \Omega_{V}^{g},\; w\in W.\label{epsiarelation}
\end{eqnarray}
In particular, for $\a\in K\cap 2L^{o}$, 
$Y_{\Omega}(\psi_{\a}({\bf 1}),0)$ is an $\Omega_{V}$-automorphism of $W$.
\ep

\pf Since $W$ is an $\Omega_{V}$-module and $\psi_{\a}$ is 
an $\Omega_{V}$-automorphism of $V$, we have
\begin{eqnarray}
{d\over dz}Y_{\Omega}(\psi_{\a}({\bf 1}),z)
=Y_{\Omega}(L_{\Omega}(-1)\psi_{\a}({\bf 1}),z)
=Y_{\Omega}(\psi_{\a}L_{\Omega}(-1){\bf 1},z)=0.
\end{eqnarray}
Then $Y_{\Omega}(\psi_{\a}({\bf 1}),z)\in \End W$.

Let $s\in P(W)$ with $W^{s}\ne 0$. Then $\Omega_{W}^{s}\ne 0$.
Let $0\ne w\in \Omega_{W}^{s}$. 
Since by assumption $V$ is simple,
$Y(\psi_{\a}({\bf 1}),z)w\ne 0$ ([DL2], Proposition 11.9).
Then 
\begin{eqnarray}
0\ne z^{{1\over\ell}\<\ell\a,s\>}
Y_{\Omega}(\psi_{\a}({\bf 1}),z)w\in W[[z,z^{-1}]]
\;\;\;\mbox{ and }\;\; Y_{\Omega}(\psi_{\a}({\bf 1}),z)w\in W,
\end{eqnarray}
noting that $\psi_{\a}({\bf 1})\in \Omega_{V}^{\ell\a}$ from 
(\ref{epsiagrading}).
Consequently, $\<\a,s\>\in \Z$. Since 
$\a\in K\subset L^{o}$ and 
$P(W)$ by definition is an $L$-subset of ${\bf h}$ generated
by $s$ with $W^{s}\ne 0$, we have $\<\a,s\>\in \Z$ for $s\in P(W)$.

Let $v\in \Omega_{V}^{g},\; w\in W^{s}$ with $g\in L,\; s\in P(W)$.
By the generalized weak commutativity, there is a nonnegative integer $k$
such that
\begin{eqnarray}\label{e4gwcomm}
& &(z_{1}-z_{2})^{k+{1\over\ell}\<g,\ell\a\>}
Y_{\Omega}(v,z_{1})Y_{\Omega}(\psi_{\a}({\bf 1}),z_{2})w\nonumber\\
&=&(-1)^{k}(z_{2}-z_{1})^{k+{1\over\ell}\<g,\ell\a\>}
Y_{\Omega}(\psi_{\a}({\bf 1}),z_{2})Y_{\Omega}(v,z_{1})w.
\end{eqnarray}
Since $g\in L$ and $\a\in K\subset L^{o}$, 
we have $\<g,\a\>\in \Z$. Let $k'\in \N$ be such that
$$r=k'+k+\<g,\a\>\ge 0.$$
Then multiplying (\ref{e4gwcomm}) by $(z_{1}-z_{2})^{k'}$ 
we get
\begin{eqnarray}
(z_{1}-z_{2})^{r}Y_{\Omega}(v,z_{1})Y_{\Omega}(\psi_{\a}({\bf 1}),z_{2})w
=(-1)^{\<g,\a\>}(z_{1}-z_{2})^{r}
Y_{\Omega}(\psi_{\a}({\bf 1}),z_{2})Y_{\Omega}(v,z_{1})w.
\end{eqnarray}
Setting $z_{2}=0$, we obtain
\begin{eqnarray}
z_{1}^{r}Y_{\Omega}(v,z_{1})Y_{\Omega}(\psi_{\a}({\bf 1}),0)w
=(-1)^{\<g,\a\>}z_{1}^{r}Y_{\Omega}(\psi_{\a}({\bf 1}),0)Y_{\Omega}(v,z_{1})w,
\end{eqnarray}
which immediately gives (\ref{ealternating}).

For $v\in \Omega_{V}^{g},\; w\in W^{s}$ with $g\in L,\; s\in P(W)$,
by the generalized weak associativity, there is a nonnegative 
integer $l$ such that
\begin{eqnarray}
& &(z_{0}+z_{2})^{l+{1\over\ell}\<g,s\>}Y_{\Omega}(v,z_{0}+z_{2})
Y_{\Omega}(\psi_{\a}({\bf 1}),z_{2})w\nonumber\\
&=&(z_{2}+z_{0})^{l+{1\over\ell}\<g,s\>}
Y_{\Omega}(Y_{\Omega}(v,z_{0})\psi_{\a}({\bf 1}),z_{2})w.
\end{eqnarray}
Since
\begin{eqnarray}
Y_{\Omega}(Y_{\Omega}(v,z_{0})\psi_{\a}({\bf 1}),z_{2})w
&=&Y_{\Omega}(\psi_{\a}Y_{\Omega}(v,z_{0}){\bf 1},z_{2})w\nonumber\\
&=&Y_{\Omega}(\psi_{\a}e^{z_{0}L_{\Omega}(-1)}v,z_{2})w\nonumber\\
&=&Y_{\Omega}(e^{z_{0}L_{\Omega}(-1)}\psi_{\a}(v),z_{2})w\nonumber\\
&=&Y_{\Omega}(\psi_{\a}(v),z_{2}+z_{0})w,
\end{eqnarray}
we get
\begin{eqnarray}\label{etemp1}
(z_{0}+z_{2})^{l+{1\over\ell}\<g,s\>}Y_{\Omega}(v,z_{0}+z_{2})
Y_{\Omega}(\psi_{\a}({\bf 1}),z_{2})w
=(z_{2}+z_{0})^{l+{1\over\ell}\<g,s\>}Y_{\Omega}(\psi_{\a}(v),z_{2}+z_{0})w.
\end{eqnarray}
Note that
$$Y_{\Omega}(\psi_{\a}({\bf 1}),z_{2})w\in W^{s+\ell\a},\;\;\;\;
\psi_{\a}(v)\in \Omega_{V}^{g+\ell\a}$$
and 
$${1\over\ell}\<g,s+\ell\a\>,\;\;\;\;{1\over\ell}\<g+\ell\a,s\>
\in {1\over\ell}\<g,s\>+\Z.$$
Then we may update $l$ by a larger positive integer so that
$$z^{l+{1\over\ell}\<g,s\>}Y_{\Omega}(\psi_{\a}(v),z)w\in W[[z]].$$
Now we may set $z_{2}=0$ in (\ref{etemp1}) with $\ell$ being updated
to get
\begin{eqnarray}
z_{0}^{l+{1\over\ell}\<g,s\>}Y_{\Omega}(v,z_{0})
Y_{\Omega}(\psi_{\a}({\bf 1}),0)w
=z_{0}^{l+{1\over\ell}\<g,s\>}Y_{\Omega}(\psi_{\a}(v),z_{0})w.
\end{eqnarray}
Then (\ref{epsiarelation}) follows immediately.
$\;\;\;\;\Box$

\br{rvacuumlike}
{\em Let $0\ne u\in \Omega_{V}^{\ell \a}$, where $\a\in {1\over\ell}L$,
be such that $L_{\Omega}(-1)u=0$.
The same proof of Proposition \ref{pvacuumlike}
with $u$ in place of $\psi_{\a}({\bf 1})$ shows that 
all the assertions of
Proposition \ref{pvacuumlike} hold with $u$ in place of $\psi_{\a}({\bf 1})$.
Furthermore, set
\begin{eqnarray}
f=Y_{\Omega}(u,0)(-1)^{\a(0)}\in \End W.
\end{eqnarray}
Then for $v\in \Omega_{V}^{g},\;g\in L, \; w\in W^{s},\; s\in P(W)$,
using (\ref{ealternating}) we get
\begin{eqnarray}
f(Y_{\Omega}(v,z)w)&=&(-1)^{\<\a,g+s\>}Y_{\Omega}(u,0)Y_{\Omega}(v,z)w
=(-1)^{\<\a,s\>}Y_{\Omega}(v,z)Y_{\Omega}(u,0)w\nonumber\\
& =&Y_{\Omega}(v,z)f(w).
\end{eqnarray}
That is, $f$ is an $\Omega_{V}$-endomorphism of $W$. Clearly, $f$ 
preserves $\Omega_{W}$.}
\er

%\br{rvacuumlike}
%{\em For a $V$-module $W$, a vector $e$ in $W$ is said to 
%be vacuum-like [Li1] if $Y(v,z)e\in W[[z]]$ for all $v\in V$, 
%which was proved to be equivalent to $L(-1)e=0$. 
%It was proved in [Li1] that the linear map
%$v\mapsto v_{-1}e=Y(v,0)e$ from $V$ to $W$ is 
%a $V$-homomorphism. When $W=V$, we have $Y(v,0)e=Y(e,0)v$.
%This fact motivated Proposition \ref{pvacuumlike}.}
%\er

Note that the $\Z$-bilinear form $\ell \<\cdot,\cdot\>$ on $K$
is $\Z$-valued because
$K\subset L^{o}$ (by definition) and $K\subset {1\over \ell}L$
(by Corollary \ref{crelationAQ}).

\bp{pcommutatormap}
Let $\psi$ be a section of the extension (\ref{eshortsequence1})
as in Proposition \ref{pvacuumlike}. Then
\begin{eqnarray}\label{epsivertexoperator}
\psi_{\a}=Y_{\Omega}(\psi_{\a}({\bf 1}),0)(-1)^{\a(0)}
\;\;\;\mbox{ for }\a\in K.
\end{eqnarray}
Furthermore,
\begin{eqnarray}\label{ecommutatorrelation}
\psi_{\a}\psi_{\b}=(-1)^{\ell \<\a,\b\>}\psi_{\b}\psi_{\a}
\;\;\;\mbox{ for }\a,\b\in K.
\end{eqnarray}
\ep

\pf Combining (\ref{ealternating}) and 
(\ref{epsiarelation}) we get
\begin{eqnarray}\label{ecombination}
Y_{\Omega}(\psi_{\a}(v),z)w=(-1)^{\<g,\a\>}Y_{\Omega}(\psi_{\a}({\bf 1}),0)
Y_{\Omega}(v,z)w.
\end{eqnarray}
Setting $w={\bf 1}$ and taking limit $z\rightarrow 0$ we obtain
(\ref{epsivertexoperator}).

Setting $v=\psi_{\b}({\bf 1})$ in (\ref{ealternating}) we get
\begin{eqnarray}\label{eypasicomm}
Y_{\Omega}(\psi_{\b}({\bf 1}),z)Y_{\Omega}(\psi_{\a}({\bf 1}),0)
=(-1)^{\ell\<\a,\b\>}Y_{\Omega}(\psi_{\a}({\bf 1}),0)
Y_{\Omega}(\psi_{\b}({\bf 1}),z),
\end{eqnarray}
noting that $\psi_{\b}({\bf 1})\in \Omega_{V}^{\ell\b}$.
Then (\ref{ecommutatorrelation}) follows immediately.
$\;\;\;\;\Box$

Now the group $\Aut_{\Omega_{V}}\Omega_{V}$ has been completely 
determined in terms of $K$ and the $\Z$-bilinear form 
$\ell \<\cdot,\cdot\>$. Furthermore, since by assumption $L$ is a
free group of finite rank, $\ell K\;(\subset L)$ is a 
free group of finite rank, and so is $K$.
The following is an immediate consequence
of (\ref{ecommutatorrelation}) (cf. [DLM1], Lemma 3.8):

\bc{cellklattice}
The group $K$ equipped with the $\Z$-bilinear form 
$\ell \<\cdot,\cdot\>$ is an even lattice. $\;\;\;\;\Box$
\ec

{}From [FLM], there exists a $\<\pm 1\>$-valued function 
$\e_{K}(\cdot,\cdot)$ on $K\times K$ such that
\begin{eqnarray}
& &\e_{K}(\a+\b,\gamma)=\e_{K}(\a,\gamma)\e_{K}(\b,\gamma),\;\;\;
\e_{K}(\gamma,\a+\b)=\e_{K}(\gamma,\a)\e_{K}(\gamma,\b),\\
& &\e_{K}(\a,\b)\e_{K}(\b,\a)^{-1}
=(-1)^{\ell\<\a,\b\>}
\;\;\;\mbox{ for }\a,\b,\gamma\in K.
\end{eqnarray}

{\bf Basic Assumption 3:} From now on we assume that 
$\psi$ is a section of the extension (\ref{eshortsequence1})
obtained through Proposition \ref{pautostructure} such that
\begin{eqnarray}\label{epsiassociator}
\psi_{\a}\psi_{\b}=\e_{K}(\a,\b)\psi_{\a+\b}
\;\;\;\mbox{ for }\a,\b\in K.
\end{eqnarray}
(The existence of such a section follows from 
Proposition \ref{pcommutatormap} and [FLM], Chapter 5.)

Note that for the generalized vertex algebra $\Omega_{V}$,
the $\C/2\Z$-valued form $(\cdot,\cdot)$ is given by
the form ${1\over\ell}\<\cdot,\cdot\>$ on $L$.
For a certain technical reason, instead of considering $K$ 
equipped with the form $\ell\<\cdot,\cdot\>$ we shall consider
$\ell K$ equipped with the form ${1\over\ell}\<\cdot,\cdot\>$.
By identifying $K$ with $\ell K$ in the obvious way we have 
an extension of $\ell K$:
\begin{eqnarray}\label{eshortsequence2}
1\rightarrow \C^{\times}\rightarrow \Aut_{\Omega_{V}}\Omega_{V}
\rightarrow \ell K\rightarrow 1,
\end{eqnarray}
with the commutator map $(-1)^{{1\over\ell}\<\cdot,\cdot\>}$.

Let $\e(\cdot,\cdot)$ be
the $\<\pm 1\>$-valued function on $\ell K$ defined by
\begin{eqnarray}
\e(\a,\b)=(-1)^{{1\over\ell} \<\a,\b\>}\e_{K}({1\over\ell}\a,{1\over\ell}\b)
=\e_{K}({1\over\ell}\b,{1\over\ell}\a)
\;\;\;\mbox{ for }\a,\b\in \ell K.
\end{eqnarray}
Then
\begin{eqnarray}
& &\e(\a+\b,\gamma)=\e(\a,\gamma)\e(\b,\gamma),\;\;\;
\e(\gamma,\a+\b)=\e(\gamma,\a)\e(\gamma,\b),\\
& &\e(\a,\b)\e(\b,\a)^{-1}
=(-1)^{{1\over\ell}\<\a,\b\>}
\;\;\;\mbox{ for }\a,\b,\gamma\in \ell K.
\end{eqnarray}
Let $\tilde{\ell K}$ be the central 
extension of $\ell K$ by the multiplicative group $\C^{\times}$ 
associated with the $2$-cocycle $\e$. That is,
\begin{eqnarray}
\tilde{\ell K}=\C^{\times}\times \ell K
\end{eqnarray}
with multiplication
\begin{eqnarray}
(a,\a)\cdot (b,\b)=(ab\e(\a,\b), \a+\b)
\;\;\;\mbox{ for }a,b\in \C^{\times},\; \a,\b\in \ell K.
\end{eqnarray}
Set
\begin{eqnarray}
\hat{\ell K}=\<\pm 1\>\times \ell K\subset \tilde{\ell K}.
\end{eqnarray}
Then $\hat{\ell K}$ is a subgroup and it is the central 
extension of $\ell K$ by $\<\pm 1\>$ associated with 
the $2$-cocycle $\e$.
We define the twisted group algebra $\C_{\e}[\ell K]$ to be 
the associative algebra with a basis $\{e^{\a}\;|\; \a\in \ell K\}$
and with multiplication
\begin{eqnarray}
e^{\a}\cdot e^{\b}=\e(\a,\b)e^{\a+\b}\;\;\;\mbox{ for }\a,\b\in \ell K.
\end{eqnarray}

\bd{dsigmawb}
{\em For $\a\in \ell K$, $W\in {\cal{C}}$, we define
\begin{eqnarray}
\sigma_{\a}^{W}=Y_{\Omega}(\psi_{{1\over\ell}\a}({\bf 1}),0)
\in \End W,
\end{eqnarray}
recalling Proposition \ref{pvacuumlike}.
When $W=V$, we simply use $\sigma_{\a}$.}
\ed

In view of Proposition \ref{pcommutatormap} we have
\begin{eqnarray}\label{esigmapsi}
\sigma_{\a}=\psi_{{1\over\ell}\a}(-1)^{{1\over\ell}\a(0)}
\;\;\;\mbox{ for }\a\in \ell K.
\end{eqnarray}
Note that $\psi_{\a}$ shifts the grading by $\ell \a$ while
$\sigma_{\a}$ shifts the grading by $\a$.

Let $\a\in \ell K, \; v\in V^{g},\; g\in L$. 
Using (\ref{ecombination}) and (\ref{esigmapsi}) we get
\begin{eqnarray}
\sigma_{\a}^{W}Y_{\Omega}(v,z)
=(-1)^{{1\over\ell}\<\a,g\>}Y_{\Omega}(\psi_{{1\over\ell}\a}(v),z)
=Y_{\Omega}(\sigma_{\a}(v),z).
\end{eqnarray}
In particular,
\begin{eqnarray}
\sigma_{\a}^{W}=Y_{\Omega}(\psi_{{1\over\ell}\a}({\bf 1}),z)
=Y_{\Omega}(\sigma_{\a}({\bf 1}),z).
\end{eqnarray}

We summarize the basic properties as follows:

\bp{psummary}
Let $W\in {\cal{C}}$. Then
\begin{eqnarray}\label{eshift}
\sigma_{\a}^{W}(\Omega_{W}^{s})=\Omega_{W}^{s+\a}\;\;\;
\mbox{ for }s\in P(W)
\end{eqnarray}
and
\begin{eqnarray}
& &Y_{\Omega}(v,z)\sigma_{\a}^{W}
=(-1)^{{1\over\ell}\<\a,g\>}\sigma_{\a}^{W}Y_{\Omega}(v,z),\\
& &\sigma_{\a}^{W}Y_{\Omega}(v,z)=
Y_{\Omega}(\sigma_{\a}(v),z)\;\;\;\;\mbox{ for }v\in V^{g},\; g\in L.
\label{e4.31}
\end{eqnarray}
In particular,
\begin{eqnarray}
\sigma_{\a}^{W}\in \Aut_{\Omega_{V}}W
\;\;\;\mbox{ for }\a\in \ell K\cap 2\ell L^{o}.
\end{eqnarray}
Furthermore,
\begin{eqnarray}\label{eprojsigma}
\sigma_{\a_{1}}^{W}\sigma_{\a_{2}}^{W}
=\e(\a_{1},\a_{2})\sigma_{\a_{1}+\a_{2}}^{W}
\;\;\;\mbox{ for }\a_{1},\a_{2}\in \ell K.
\end{eqnarray}
\ep

\pf (\ref{eprojsigma}) follows from (\ref{epsiarelation}) and 
(\ref{epsiassociator}),  and 
all the others follow from Proposition \ref{pvacuumlike}.
$\;\;\;\;\Box$

\br{rlatticeL}
{\em Set
\begin{eqnarray}
L_{0}=\ell K\cap 2\ell L^{o}\subset L.
\end{eqnarray}
For $\a\in L_{0},\; g\in L$, since by definition
$\a\in 2\ell L^{o}$, we have $\<\a,g\>\in 2\ell \Z.$ 
That is,
\begin{eqnarray}\label{e4.30}
{1\over\ell}\<\a,g\>\in 2\Z\;\;\;\mbox{ for }\a\in L_{0},\; g\in L.
\end{eqnarray}
In particular, $L_{0}$ equipped with the form ${1\over\ell}\<\cdot,\cdot\>$
is an even lattice. }
\er

Define
\begin{eqnarray}
K^{o}=\{ h\in {\bf h}\;|\;\<\a,h\>\in \Z\;\;\;\mbox{ for }\a\in K\}.
\end{eqnarray}
Note that $K^{o}$ is also the dual of $\ell K$ equipped with 
the form ${1\over\ell}\<\cdot,\cdot\>$.

\bl{lkernal}
For any $W\in {\cal{C}}$, we have
\begin{eqnarray}
P(W)\in K^{o}.
\end{eqnarray}
Furthermore, if $K$ spans ${\bf h}$ over $\C$, $2\ell K$ is
the kernel of
the $\C/2\Z$-valued $\Z$-bilinear form $(\cdot,\cdot)$
on $K^{o}$ defined by
\begin{eqnarray}
(\a,\b)={1\over\ell}\<\a,\b\>+2\Z\;\;\;\mbox{ for }\a,\b\in K^{o}.
\end{eqnarray}
\el

\pf It follows immediately from Proposition \ref{pvacuumlike} that
$P(W)\subset K^{o}$. In particular, $L=P(V)\subset K^{o}$.
Then using Corollary \ref{crelationAQ}, we get
\begin{eqnarray}
2\ell K\subset  2L\subset K^{o}.
\end{eqnarray}
If $K$ spans ${\bf h}$ over $\C$, 
we have $({1\over 2\ell}K^{o})^{o}=2\ell K$. Then
$2\ell K$ is the kernel of $(\cdot,\cdot)$.
$\;\;\;\;\Box$

\br{rGandform}
{\em Set
\begin{eqnarray}
G=K^{o}/2\ell K,
\end{eqnarray}
an abelian group. In view of Lemma \ref{lkernal}, 
if $K$ spans ${\bf h}$ over $\C$, we have a
nondegenerate $\C/2\Z$-valued symmetric $\Z$-bilinear 
form $(\cdot,\cdot)$ on $G$ defined by
\begin{eqnarray}
(\a+2\ell K,\b+2\ell K)={1\over \ell}\<\a,\b\>+2\Z\;\;\;
\mbox{ for }\a,\b\in K^{o}.
\end{eqnarray}}
\er

\section{Quotient generalized vertex (operator) algebras $\Omega_{V}^{A}$}
In this section, we shall construct quotient generalized 
vertex algebras $\Omega_{V}^{A}$ of $\Omega_{V}$
for central subgroups $A$ of $\Aut_{\Omega_{V}}\Omega_{V}$.
For certain $A$, the quotient algebras $\Omega_{V}^{A}$ are proved to be 
generalized vertex operator algebras. We shall
construct functors between the categories of $\Omega_{V}$-modules and 
$\Omega_{V}^{A}$-modules and
prove analogous assertions of Proposition \ref{pirrational} 
for $\Omega_{V}^{A}$, by which all
irreducible $\Omega_{V}^{A}$-modules are classified 
in terms of irreducible $V$-modules.

Let $A$ be a subgroup:
\begin{eqnarray}
A\subset L_{0}=\ell K\cap 2\ell L^{o}\subset L.
\end{eqnarray}
We set
\begin{eqnarray}
I_{V}^{A}=\mbox{linear span}\{ v-\sigma_{\a}(v)\;|\;v\in \Omega_{V},
\; \a\in A\}.
\end{eqnarray}

\bl{lideal}
The subspace $I_{V}^{A}$ of $\Omega_{V}$ is an ideal of $\Omega_{V}$.
\el

\pf By Proposition \ref{psummary}, for $\a\in A$
$\sigma_{\a}$ is an $\Omega_{V}$-homomorphism, so it follows that
$I_{V}^{A}$ is a left ideal of $\Omega_{V}$. 
For $\a\in A,\;u,v\in \Omega_{V}$, by Proposition \ref{psummary},
we have
\begin{eqnarray}
Y_{\Omega}(v-\sigma_{\a}(v),z)u
=Y_{\Omega}(v,z)u-Y_{\Omega}(\sigma_{\a}(v),z)u
=Y_{\Omega}(v,z)u-\sigma_{\a}
Y_{\Omega}(v,z)u\in I_{V}^{A}\{z\}.
\end{eqnarray}
This proves that $I_{V}^{A}$ is also a right ideal.
$\;\;\;\;\Box$

Because $A\subset L_{0}=\ell K\cap 2\ell L^{o}$,
${1\over\ell}\<\cdot,\cdot\>$ gives rise to
a $\C/2\Z$-valued $\Z$-bilinear form $(\cdot,\cdot)$ on $L/A$.
Then the quotient space of $\Omega_{V}$ modulo the two-sided ideal 
$I_{V}^{A}$ is a generalized vertex algebra with grading group $L/A$ 
and forms $c(\cdot,\cdot)=1$ and $(\cdot,\cdot)$.
We denote this quotient generalized vertex algebra by $\Omega_{V}^{A}$
(cf. [DL2]).
We shall still use ${\bf 1}$ and $\omega_{\Omega}$ for
the  vacuum vector and Virasoro vector for the 
quotient generalized vertex algebra $\Omega_{V}^{A}$.
 
\br{rellKquotient}
{\em  Clearly, one may define $I_{V}^{\ell K}$ and 
$\Omega_{V}^{\ell K}$ in the same way. 
However, $I_{V}^{\ell K}$ is a right ideal, but not a left ideal. 
Also, notice that ${1\over\ell}\<\cdot,\cdot\>$ does not
give rise to a $\C/2\Z$-valued form on $L/\ell K$.
It turns out that similar to a situation in [DLM1], one can make
the quotient space $\Omega_{V}^{\ell K}$ an abelian intertwining algebra 
in the sense of [DL2].}
\er

Note that $2\ell K\subset L_{0}$. An important case will be the one
with $A=2\ell K$.

\bl{lquotientgvoa}
Let $A$ be a subgroup of $L_{0}\;(=\ell K\cap 2\ell L^{o})$ 
such that $L/A$ is finite. Then
$(\Omega_{V}^{A}, Y_{\Omega}, {\bf 1}, \omega_{\Omega}, L/A,
(\cdot,\cdot))$ is a generalized vertex operator algebra
except that the form $(\cdot,\cdot)$ on $L/A$ may be degenerate.
\el

\pf We only need to prove that $\Omega_{V}^{A}$ graded by 
$L_{\Omega}(0)$-weight satisfies the two grading restrictions. 
{}From (\ref{eL(m)h(n)}), $L(0)$ preserves
$\Omega_{V}^{g}$ for $g\in L$.
For $g\in L$, because $V$ satisfies the two grading restrictions and
\begin{eqnarray}
& &L(0)=L_{\Omega}(0)+L_{\bf h}(0),\\
& &L_{\bf h}(0)={\<g,g\>\over 2\ell}
\;\;\;\mbox{ on }\Omega_{V}^{g},
\end{eqnarray}
recalling Lemma \ref{lvirasoro}, 
$\Omega_{V}^{g}$ graded by $L_{\Omega}(0)$-weight
satisfies the two grading restrictions. Let $\{g_{1},\dots, g_{n}\}$ 
be a complete set of representatives of cosets of $A$ in $L$.
Then the covering map from $\Omega_{V}$ to $\Omega_{V}^{A}$
restricted to 
$\sum_{i=1}^{n}\Omega_{V}^{g_{i}}$ is an $\Omega_{V}^{0}$-isomorphism.
Consequently, $\Omega_{V}^{A}$ graded by $L_{\Omega}(0)$-weight
satisfies the two grading restrictions. 
$\;\;\;\;\Box$.

Regarding generators of $\Omega_{V}^{A}$, we have
the following immediate consequence of 
Proposition \ref{pgenerators}:

\bp{pgeneratorquotientgva}
Let $A$ be a subgroup of $L_{0}\;(=\ell K\cap 2\ell L^{o})$ and let
$U$ be an ${\bf h}$-submodule of $\Omega_{V}$ such that
$U+{\bf h}$ generates $V$ as a vertex algebra. Then
the image of $U$ in $\Omega_{V}^{A}$ generates $\Omega_{V}^{A}$
as a generalized vertex algebra.
$\;\;\;\;\;\Box$
\ep

\br{rautoquotientgva}
{\em Let $\tau$ be an automorphism of the vertex operator algebra $V$
as in Proposition \ref{pautoreduction}. 
%Assume 
%\begin{eqnarray}
%\tau(\a)=\a\;\;\;\mbox{ for }\a\in A.
%\end{eqnarray}
%We also assume that $\Omega_{V}^{\a}\subset V^{\tau}$ if $\tau(\a)=\a$.
In addition we assume that 
$\sigma_{\a}\tau=\tau \sigma_{\a}$ for $\a\in A$.
Then $I_{V}^{A}$ is stable under $\tau$. Consequently,
$\tau$ gives rise to a canonical automorphism of $\Omega_{V}^{A}$.}
\er

Next, we shall construct an $\Omega_{V}^{A}$-module $\Omega_{W}^{A}$ from
a $V$-module $W$. Similarly, we define
\begin{eqnarray}
I_{W}^{A}
=\mbox{linear span} \{ w-\sigma_{\a}^{W}(w)\;|\; w\in \Omega_{W},\; \a\in A\}.
\end{eqnarray}
The same proof of Lemma \ref{lideal} gives:

\bl{lIWB}
The subspace  $I_{W}^{A}$ is a $(P(W)/A)$-graded $\Omega_{V}$-submodule 
of $\Omega_{W}$ and 
\begin{eqnarray}\label{eIWB=IVBW}
Y_{\Omega}(v-\sigma_{\a}(v),z)\Omega_{W}\subset I_{W}^{A}\{z\}
\;\;\;\mbox{ for }v\in \Omega_{V},\; \a\in A. \;\;\;\;\Box
\end{eqnarray}
\el

We define $\Omega_{W}^{A}$ to be the quotient $\Omega_{V}$-module 
of $\Omega_{W}$ modulo the submodule $I_{W}^{A}$.

\bp{pquotientmodule}
For any weak $V$-module $W\in {\cal{C}}$, $\Omega_{W}^{A}$ is 
a $(P(W)/A)$-graded 
$\Omega_{V}^{A}$-module. Furthermore, for $\b\in L^{o}$,
\begin{eqnarray}\label{eisomeb}
e^{\b}\otimes I_{W}^{A}=I_{W^{(\b)}}^{A}
\end{eqnarray}
and the linear map $e^{\b}\otimes \cdot$ gives rise to an 
$\Omega_{V}^{A}$-isomorphism from $\Omega_{W}^{A}$ 
onto $\Omega_{W^{(\b)}}^{A}$.
\ep

\pf It follows immediately from Lemma \ref{lIWB} that 
$\Omega_{W}^{A}$ is a $(P(W)/A)$-graded $\Omega_{V}^{A}$-module.

Recall from Proposition \ref{pidenfitication} that $e^{\b}\otimes\cdot$ is 
an $\Omega_{V}$-isomorphism from $\Omega_{W}$ onto $\Omega_{W^{(\b)}}$.
For $w\in W,\; \a\in A$, we have
\begin{eqnarray}
e^{\b}\otimes (w-\sigma_{\a}^{W}(w))&=&
e^{\b}\otimes Y_{\Omega}({\bf 1}-\sigma_{\a}({\bf 1}),z)w\nonumber\\
&=&Y_{\Omega}({\bf 1}-\sigma_{\a}({\bf 1}),z)(e^{\b}\otimes w)\nonumber\\
&=&e^{\b}\otimes w-\sigma_{\a}^{W^{(\b)}}(e^{\b}\otimes w).
\end{eqnarray}
Then (\ref{eisomeb}) follows immediately. Consequently,
$e^{\b}\otimes \cdot$ gives rise to an 
$\Omega_{V}^{A}$-isomorphism from $\Omega_{W}^{A}$ 
onto $\Omega_{W^{(\b)}}^{A}$.$\;\;\;\;\Box$

\br{rmorequotient}
{\em Let $U$ be an $\Omega_{V}$-${\bf h}$-module. 
By Proposition \ref{pback}, we have
a weak $V$-module $E(U)\in {\cal{C}}$ with $\Omega_{E(U)}=U$.
Then $\Omega_{E(U)}^{A}=U^{A}$, where $U^{A}$
is defined as a quotient space of $U$ in the obvious way.
By Proposition \ref{pquotientmodule}, $U^{A}$
is a $P(U)/A$-graded $\Omega_{V}^{A}$-module.}
\er

Next, we shall prove the irreducibility of $\Omega_{V}^{A}$-module
$\Omega_{W}^{A}$ with $W$ being irreducible. First, we 
prove the following result:

\bl{lfacts1}
Let $W\in {\cal{C}}$ be a weak $V$-module and 
let $\xi$ be the covering map from 
$\Omega_{W}$ to $\Omega_{W}^{A}$. Then for 
any $(P(W)/A)$-graded $\Omega_{V}^{A}$-submodule $U$ of
$\Omega_{W}^{A}$, there exists an $\Omega_{V}$-${\bf h}$-submodule
$\tilde{U}$ of $\Omega_{W}$ such that $\xi (\tilde{U})=U$, i.e.,
$\tilde{U}+I_{W}^{A}=\xi ^{-1}(U)$.
\el

\pf Let $\tilde{U}$ be the (unique) maximal 
$\Omega_{V}$-${\bf h}$-submodule
of $\Omega_{W}$ such that $\xi (\tilde{U})\subset U$.
Then we must prove that $U= \xi (\tilde{U})$, 
for which it suffices to prove
$U\subset \xi (\tilde{U})$.

In the following we shall use the proof of Theorem 14.20 of [DL2].
Set $S=P(W)/A$.
Since by hypothesis $U$ is $S$-graded, i.e.,
$U=\coprod_{s\in S}U^{s}$, it suffices to prove that
$U^{s}\subset \xi (\tilde{U})$ for $s\in S$.

Let $s\in S$ and let $e(s)\in P(W)$ be such that $s=e(s)+A$.
Let $u\in U^{s}$.  Consider an element of $\xi^{-1}(u)$:
\begin{eqnarray}
\sum_{\a\in A}u^{e(s)+\a}
\end{eqnarray}
(a finite sum), where $u^{e(s)+\a}\in \Omega_{W}^{e(s)+\a}$.
Set
\begin{eqnarray}
\tilde{u}=\sum_{\a\in A}\sigma_{-\a}^{W}(u^{e(s)+\a})\in \Omega_{W}^{e(s)}.
\end{eqnarray}
Then $\xi (\tilde{u})=u$.
Denote by $T(\tilde{u})$ the $\Omega_{V}$-submodule
of $\Omega_{W}$ generated by $\tilde{u}$. 
Since $\tilde{u}\in \Omega_{W}^{e(s)}$,
$T(\tilde{u})$ is clearly an $\Omega_{V}$-${\bf h}$-submodule.
With $\xi$ being an $\Omega_{V}$-homomorphism
and $\xi(\tilde{u})=u\in U$ we have $\xi (T(\tilde{u}))\subset U$.
By definition, $T(\tilde{u})\subset \tilde{U}$. Hence
$$u\in \xi(T(\tilde{u}))\subset \xi(\tilde{U}).$$
Thus, $U^{s}\subset  \xi(\tilde{U})$.
Therefore, $U\subset \xi(\tilde{U}).$
This completes the proof.$\;\;\;\;\Box$

As an immediate consequence of Lemma \ref{lfacts1} and
Theorem \ref{tequivalence} we have (cf. [DL2], Theorem 14.20):

\bp{pquotientmoduleirred}
For any irreducible weak $V$-module $W\in {\cal{C}}$, 
$\Omega_{W}^{A}$ is a $(P(W)/A)$-graded
irreducible  $\Omega_{V}^{A}$-module. 
In particular, $\Omega_{V}^{A}$ is 
an $(L/A)$-graded simple generalized vertex algebra.
$\;\;\;\;\Box$
\ep

\br{rmaybe}
{\em Note that from Lemma \ref{lkernal},
$\Omega_{W}^{A}$ is a natural $(K^{o}/L_{0})$-graded 
$\Omega_{V}^{A}$-module for any weak $V$-module $W\in {\cal{C}}$.
When $A\ne L_{0}$, one can show that $\Omega_{W}^{A}$ has nontrivial 
$(K^{o}/L_{0})$-graded $\Omega_{V}^{A}$-submodules.
(Cf. [DL2], Remark 14.21.)}
\er

Note that for any weak $V$-module $W\in {\cal{C}}$, 
$\Omega_{W}^{A}$ is a $(K^{o}/A)$-graded $\Omega_{V}^{A}$-module.
Next, we shall prove that $\Omega_{W}^{A}$ 
for irreducible $W\in {\cal{C}}$ exhaust 
all irreducible $(K^{o}/A)$-graded $\Omega_{V}^{A}$-modules.

{\em Fix a complete set $S$ of representatives of 
cosets of $A$ in $K^{o}$.} Then
for any $\lambda\in K^{o}$, there exists a unique
$t(\lambda)\in A$ such that $\lambda-t(\lambda)\in S$.
Denote by $t$ the map from $K^{o}$ to $A$ sending $\lambda$ 
to $t(\lambda)$.
Clearly, 
\begin{eqnarray}\label{epropertyt}
t(\a+\lambda)=\a+t(\lambda)\;\;\;\mbox{ for }\a\in A,\; \lambda\in K^{o}.
\end{eqnarray}
It is also clear that the map $\lambda\mapsto (t(\lambda),\lambda-t(\lambda))$
is a bijection from $K^{o}$ to $A\times S$.

Let $U$ be a $(K^{o}/A)$-graded $\Omega_{V}^{A}$-module. 
View $U$ as a natural $(K^{o}/A)$-graded $\Omega_{V}$-module.
Set
\begin{eqnarray}
U_{A}=\C[A]\otimes U.
\end{eqnarray}
For $\lambda\in K^{o}$, set
\begin{eqnarray}
U_{A}^{\lambda}=\C e^{t(\lambda)}\otimes U^{\lambda-t(\lambda)+A}.
\end{eqnarray}
Then $U_{A}=\coprod_{\lambda\in K^{o}}U_{A}^{\lambda}$.

For $v\in \Omega_{V}^{g},\;\a\in A,\;w\in U^{s+A},\;s\in S$, we define
\begin{eqnarray}\label{edef4.55}
Y_{\Omega}(v,z)(e^{\a}\otimes w)
=e^{\a+t(g+s)}\otimes Y_{\Omega}^{A}(v,z)w.
\end{eqnarray}
Then
\begin{eqnarray}
Y_{\Omega}(v,z)U_{A}^{\lambda}\subset U_{A}^{\lambda+g}\{z\}
\;\;\;\mbox{ for }\lambda\in K^{o}.
\end{eqnarray}
For $h\in {\bf h}$ we define
\begin{eqnarray}
h(e^{\a}\otimes w)=\<h,\a+s\>(e^{\a}\otimes w).
\end{eqnarray}

\bp{precover}
Let $U$ be a $(K^{o}/A)$-graded $\Omega_{V}^{A}$-module.
Then $U_{A}$ is a $K^{o}$-graded
$\Omega_{V}$-${\bf h}$-module such that $(U_{A})^{A}=U$.
Furthermore, if $U=\Omega_{W}^{A}$ for some $W\in {\cal{C}}$, then
$U_{A}\simeq \Omega_{W}$ as a $K^{o}$-graded $\Omega_{V}$-module.
\ep

\pf Let $u\in \Omega_{V}^{g_{1}},\;v\in \Omega_{V}^{g_{2}},
\; \a\in A,\; w\in U^{s+A},\; s\in S$.
By definition,
\begin{eqnarray}
Y_{\Omega}(v,z_{2})w\in U^{g_{2}+s-t(g_{2}+s)+A}\{z_{2}\}
\end{eqnarray}
with $g_{2}+s-t(g_{2}+s)\in S$.
Note that from (\ref{epropertyt}) we have
\begin{eqnarray}
t(g_{2}+s)+t(g_{1}+(g_{2}+s-t(g_{2}+s)))=t(g_{1}+g_{2}+s).
\end{eqnarray}
Then we have
\begin{eqnarray}
Y_{\Omega}(u,z_{1})Y_{\Omega}(v,z_{2})(e^{\a}\otimes w)
&=&Y_{\Omega}(u,z_{1})(e^{\a+t(g_{2}+s)}\otimes 
Y_{\Omega}^{A}(v,z_{2})w)\nonumber\\
&=&e^{\a+t(g_{1}+g_{2}+s)}\otimes 
Y_{\Omega}^{A}(u,z_{1})Y_{\Omega}^{A}(v,z_{2})w
\end{eqnarray}
and symmetrically,
\begin{eqnarray}
Y_{\Omega}(v,z_{2})Y_{\Omega}(u,z_{1})(e^{\a}\otimes w)
=e^{\a+t(g_{1}+g_{2}+s)}\otimes 
Y_{\Omega}^{A}(v,z_{2})Y_{\Omega}^{A}(u,z_{1})w.
\end{eqnarray}
With $Y_{\Omega}(u,z_{0})v\in \Omega_{V}^{g_{1}+g_{2}}\{z\}$,
we have
\begin{eqnarray}
Y_{\Omega}(Y_{\Omega}(u,z_{0})v,z_{2})(e^{\a}\otimes w)
=e^{\a+t(g_{1}+g_{2}+s)}\otimes 
Y_{\Omega}^{A}(Y_{\Omega}(u,z_{0})v,z_{2})w.
\end{eqnarray}
Then the generalized Jacobi identity of $Y_{\Omega}$ on $U_{A}$ 
follows immediately from the generalized Jacobi identity of 
$Y_{\Omega}^{A}$ on $U$. Clearly, $Y_{\Omega}({\bf 1},z)=1$.
Thus $U_{A}$ is a $K^{o}$-graded
$\Omega_{V}$-module.
It follows from the second part of the proof of 
Proposition \ref{pirrational} that
$U_{A}$ is an $\Omega_{V}$-${\bf h}$-module.

For $\a,\b\in A,\; w\in U^{s+A},\; s\in S$, we have
\begin{eqnarray}
\sigma_{\a}^{U_{A}}(e^{\b}\otimes w)
&=&Y_{\Omega}(\sigma_{\a}({\bf 1}),z)(e^{\b}\otimes w)\nonumber\\
&=&e^{\a+\b}\otimes Y_{\Omega}^{A}(\sigma_{\a}({\bf 1}),z)w\nonumber\\
&=&e^{\a+\b}\otimes Y_{\Omega}^{A}({\bf 1},z)w\nonumber\\
&=&e^{\a+\b}\otimes w.
\end{eqnarray}
Then $(U_{A})^{A}=U$. 

Let $U=\Omega_{W}^{A}$ and let $\zeta$ be the covering map from
$\Omega_{W}$ onto $\Omega_{W}^{A}$. Define a linear map
\begin{eqnarray}
\eta: \Omega_{W}&\rightarrow& 
(\Omega_{W}^{A})_{A}=\C[A]\otimes \Omega_{W}^{A}\nonumber\\
w&\mapsto& e^{t(\lambda)}\otimes \zeta(w)
\end{eqnarray}
for $w\in \Omega_{W}^{\lambda},\; \lambda\in P(W)\subset K^{o}$.
Clearly, $\eta$ is a linear isomorphism preserving the $K^{o}$-grading.
Let $v\in \Omega_{V}^{g},\; g\in L$. 
Recall that $\lambda=s+t(\lambda)$, where $s\in S,\;t(\lambda)\in A$.
Then 
$$t(\lambda+g)=t(t(\lambda)+g+s)=t(\lambda)+t(g+s).$$ 
Using the definition (\ref{edef4.55}) we get
\begin{eqnarray}
\eta(Y_{\Omega}(v,z)w)&=&e^{t(g+\lambda)}\otimes \zeta(Y_{\Omega}(v,z)w)
\nonumber\\
&=&e^{t(g+\lambda)}\otimes Y_{\Omega}^{A}(v,z)\zeta(w)\nonumber\\
&=&e^{t(\lambda)+t(g+s)}\otimes Y_{\Omega}^{A}(v,z)\zeta(w)\nonumber\\
&=&Y_{\Omega}(v,z)(e^{t(\lambda)}\otimes \zeta(w))\nonumber\\
&=&Y_{\Omega}(v,z)\eta(w).
\end{eqnarray}
This proves that $\eta$ is an $\Omega_{V}$-isomorphism
and completes the proof.$\;\;\;\;\Box$

Now, we present our main theorem of this section.

\bt{tmain}
For irreducible weak $V$-modules $W_{1},W_{2}\in {\cal{C}}$, 
$\Omega_{W_{1}}^{A}\simeq \Omega_{W_{2}}^{A}$
if and only if $W_{1}^{(\b)}\simeq W_{2}$ for some $\b\in L^{o}$.
Furthermore, the map $W\mapsto \Omega_{W}^{A}$ gives rise to a
one-to-one correspondence between the set of
$L^{o}$-orbits of the set of the equivalence classes of 
irreducible weak $V$-modules in ${\cal{C}}$ and the set
of equivalence classes of $(K^{o}/A)$-graded irreducible 
$\Omega_{V}^{A}$-modules.
Moreover, if ${\cal{C}}$ is semisimple, 
the category of $(K^{o}/A)$-graded $\Omega_{V}^{A}$-modules
is semisimple.
\et

\pf If $W_{1}^{(\b)}\simeq W_{2}$ for some $\b\in L^{o}$,
it follows from Proposition \ref{pquotientmodule} that 
$\Omega_{W_{1}}^{A}\simeq \Omega_{W_{2}}^{A}$.
Conversely, assume that $\Omega_{W_{1}}^{A}\simeq \Omega_{W_{2}}^{A}$.
By Proposition \ref{precover},
$$\Omega_{W_{1}}=(\Omega_{W_{1}}^{A})_{A}\simeq 
(\Omega_{W_{2}}^{A})_{A}=\Omega_{W_{2}}$$
as a $K^{o}$-graded $\Omega_{V}$-module. 
Then it follows from Theorem \ref{tequivalence} that 
$W_{1}^{(\b)}\simeq W_{2}$ for some $\b\in L^{o}$.
Now, the second assertion follows immediately.

Suppose that ${\cal{C}}$ is semisimple.
Let $U$ be any $(K^{o}/A)$-graded $\Omega_{V}^{A}$-module.
By Proposition \ref{precover}, we have an $\Omega_{V}$-${\bf h}$-module
$U_{A}$ with $U=(U_{A})^{A}$.
It follows from Theorem \ref{tequivalencecategory}
that $U_{A}$ is completely reducible.
Then it follows from Proposition \ref{pquotientmoduleirred} that
$U=(U_{A})^{A}$ (Proposition \ref{precover}) is completely reducible. 
$\;\;\;\;\Box$

\bp{pfinite}
Suppose that $\ell$ is a positive integer, 
$L$ spans ${\bf h}$ over $\C$ and that
$V$ has only finitely many irreducible weak modules in ${\cal{C}}$
up to equivalence. Then
$L/L_{0}$ is finite and $\Omega_{V}^{L_{0}}$ is a 
generalized vertex operator algebra with grading group $L/L_{0}$
except that the form 
$(\cdot,\cdot)$ may be degenerate. Furthermore, 
$\Omega_{V}^{2\ell K}$ is a generalized 
vertex operator algebra with $G=K^{o}/2\ell K$ (where the
form $(\cdot,\cdot)$ is nondegenerate on $G$) and 
for any $V$-module $W\in {\cal{C}}$, $\Omega_{W}^{2\ell K}$ is a $G$-graded
$\Omega_{V}^{2\ell K}$-module.
\ep

\pf For $\b\in L^{o}$, there exists a positive integer $n$ such that
$V^{(n\b)}\simeq V$. Hence $n\b \in K$. 
With $L^{o}$ being of finite rank, there is a positive integer 
$k_{1}$ such that $k_{1}L^{o}\subset K$.
Since $L$ is a rational lattice of finite rank and $L$ 
spans ${\bf h}$ over $\C$,
there is a positive integer $k_{2}$ such that $k_{2}L\subset L^{o}$.
Then $k_{1}k_{2}L\subset K$. Hence 
$2\ell k_{1}k_{2}L\subset 2\ell K\subset L_{0}$.
It follows immediately that $L/2\ell K$ and $L/L_{0}$ are finite.

In view of Lemma \ref{lquotientgvoa},
$\Omega_{V}^{A}$ is a generalized vertex operator algebra
except that the form $(\cdot,\cdot)$ may be degenerate.
However, if $A\subset 2\ell K$, in particular, $A=2\ell K$,
we may consider $\Omega_{V}^{A}$ as 
a $(K^{o}/2\ell K)$-graded space in the obvious way
and at the same time, by Lemma \ref{lkernal}
the form $(\cdot,\cdot)$ on $(K^{o}/2\ell K)$
is nondegenerate. Furthermore, by Lemma \ref{lkernal} again, for any
$V$-module $W\in {\cal{C}}$, $\Omega_{W}^{A}$ is a $G$-graded
$\Omega_{V}^{A}$-module.
$\;\;\;\;\Box$

For the rest of this section we shall show that $\Omega_{V}^{A}$
is equivalent to the vacuum space of a vertex subalgebra isomorphic
to the lattice vertex algebra $V_{A}$ in $V$. 
We assume that all the assumptions in Proposition \ref{pfinite} hold, 
i.e., $\ell$ is a positive integer, 
$L$ spans ${\bf h}$ over $\C$ and 
$V$ has only finitely many irreducible weak modules in ${\cal{C}}$
up to equivalence. From the proof of Proposition \ref{pfinite},
there are positive integers $k_{1},k_{2}$ such that
$k_{1}L^{o}\subset K$ and $k_{2}L\subset L^{o}$.
Consequently,
$$2\ell k_{1}k_{2}L\subset L_{0}=\ell K\cap 2\ell L^{o}.$$ 
Thus, both $\ell K$ and $L_{0}$ span ${\bf h}$ over $\C$.

Let $V_{\ell K}$ be the vertex algebra constructed 
in [FLM] (cf. [B]) from the even lattice $\ell K$ 
with the form ${1\over\ell}\<\cdot,\cdot\>$.
Note that $K^{o}$ is the dual lattice of $\ell K$ 
equipped with the form ${1\over\ell}\<\cdot,\cdot\>$.
Let $\{\b_{1},\dots, \b_{n}\}$ be a complete set of
representatives of cosets of $\ell K$ in $K^{o}$.
Then from [FLM],
\begin{eqnarray}
V_{K^{o}}=\coprod_{i=1}^{n}V_{\ell K+\b_{i}},
\end{eqnarray}
where $V_{\ell K+\b_{i}}$ and $V_{\ell K+\b_{j}}$ for $i\ne j$ are 
nonisomorphic irreducible $V_{\ell K}$-modules.
Furthermore, from [D]
any irreducible $V_{\ell K}$-module is isomorphic to 
$V_{\ell K+\b_{i}}$ for some $i$.

\bp{platticesubalgebra}
Let $V_{\ell K}$ be the vertex algebra 
associated with the even lattice $\ell K$ 
with the form ${1\over\ell}\<\cdot,\cdot\>$. 
Define a linear map $\Psi$ from $V_{\ell K}$ to $V$ by
\begin{eqnarray}
\Psi(e^{\a}\otimes h_{1}(-n_{1})\cdots h_{r}(-n_{r}))
=\ell^{-r}(\sigma_{\a}({\bf 1})\otimes h_{1}(-n_{1})\cdots h_{r}(-n_{r}))
\end{eqnarray}
for $\a\in \ell K,\;r\in \N,\; h_{i}\in {\bf h},\; n_{i}\ge 1$.
Then $\Psi$ is an injective vertex algebra homomorphism.
Furthermore, let $W\in {\cal{C}}$ and let
$0\ne w\in \Omega_{W}^{\b},\;\b \in P(W)$. Then
the $V_{\ell K}$-submodule of $W$ generated by $w$ 
through the embedding $\Psi$ is isomorphic to $V_{\ell K+\b}$.
\ep

\pf Because $0\ne \sigma_{\a}({\bf 1})\in \Omega_{V}^{\a}$ 
for $\a\in \ell K$ (Proposition \ref{psummary}), we have
\begin{eqnarray}
\Psi(V_{\ell K})
=\oplus_{\a\in L}\left(M(\ell)\otimes \C\sigma_{\a}({\bf 1})\right).
\end{eqnarray}
Then $f$ is a linear embedding of $V_{\ell K}$ into $V$.
Since ${\bf h}$ and $e^{\a}$ for $\a\in \ell K$ generates $V_{\ell K}$ as a 
vertex algebra, to prove that $\Psi$ is a vertex algebra homomorphism 
it suffices to prove
\begin{eqnarray}
\Psi(Y(v,z)u)=Y(\Psi(v),z)\Psi(u)\;\;\;\mbox{ for }
v\in {\bf h}\cup \{e^{\a}\;|\; \a\in \ell K\},
\; u\in V_{\ell K}.
\end{eqnarray}

By definition, we have
\begin{eqnarray}
& &\Psi(h)={1\over \ell}h\;\;\;\mbox{ for }h\in {\bf h},\\
& &\Psi(e^{\a})=\sigma_{\a}({\bf 1})\;\;\;\mbox{ for }\a\in \ell K.
\end{eqnarray}
Let $h_{1},h_{2}\in {\bf h},\;m,n\in \Z$. Then
\begin{eqnarray}
& &[h_{1}(m),h_{2}(n)]=m{1\over \ell}\<h_{1},h_{2}\>\delta_{m+n,0}
\;\;\;\mbox{ on }\;V_{\ell K},\\
& &\left[{1\over\ell}h_{1}(m),{1\over\ell}h_{2}(n)\right]
=m{1\over \ell}\<h_{1},h_{2}\>\delta_{m+n,0}\;\;\;\mbox{ on }\;V,
\end{eqnarray}
noting that for $h\in {\bf h}$,
$Y(h,z)=\sum_{m\in \Z}h(m)z^{-m-1}$ is a vertex operator on 
$V_{\ell K}$ and $V$. We also have
\begin{eqnarray}
& &h(n)e^{\a}=\delta_{n,0}{1\over \ell}\<h,\a\>e^{\a}
\;\;\;\mbox{ in }\;V_{\ell K},\\
& &{1\over\ell}h(n)\sigma_{\a}({\bf 1})=
\delta_{n,0}{1\over \ell}\<h,\a\>\sigma_{\a}({\bf 1})
\;\;\;\mbox{ in }\;V
\end{eqnarray}
for $h\in {\bf h},\; n\in \Z$. 
Then using induction we get
\begin{eqnarray}
\Psi(h(m)u)={1\over \ell}h(m)\Psi(u)
\;\;\;\mbox{ for }h\in {\bf h},\;m\in \Z,\; u\in V_{\ell K}.
\end{eqnarray}
That is,
\begin{eqnarray}\label{e4.67}
\Psi(Y(h,z)u)=Y(\Psi(h),z)\Psi(u)\;\;\;\mbox{ for }
h\in {\bf h},\; u\in V_{\ell K}.
\end{eqnarray}

{}From [FLM] we have
$V_{\ell K}=\C_{\e}[\ell K]\otimes M(1)$ and
\begin{eqnarray}
Y(e^{\a},z)=E^{-}(-\a,z)E^{+}(-\a,z)e^{\a}z^{\a(0)}\;\;\;\mbox{ for }
\a\in \ell K,
\end{eqnarray}
where
\begin{eqnarray}\label{e4.77}
e^{\a}\cdot e^{\b}=\e(\a,\b)e^{\a+\b}
\;\;\;\mbox{ for }\a,\b\in \ell K.
\end{eqnarray}
For $\a\in \ell K$, using (\ref{e4.31}) we have
\begin{eqnarray}\label{e4.70}
Y(\sigma_{\a}({\bf 1}),z)&=&E^{-}(-{1\over\ell}\a,z)E^{+}(-{1\over\ell}\a,z)
Y_{\Omega}(\sigma_{\a}({\bf 1}),z)z^{{1\over\ell}\a(0)}\nonumber\\
&=&E^{-}(-{1\over\ell}\a,z)E^{+}(-{1\over\ell}\a,z)
\sigma_{\a}z^{{1\over\ell}\a(0)}.
\end{eqnarray}
For $\a\in \ell K$, view $e^{\a}$ as an operator on $V_{\ell K}$. 
Since $\sigma_{\a}$
commutes with $h(n)$ for $h\in {\bf h},\; n\ne 0$ 
(by (\ref{ebrackethnpsia}) and (\ref{esigmapsi})), 
using (\ref{eprojsigma}) and (\ref{e4.77})
we get
\begin{eqnarray}
\Psi(e^{\a}u)=\sigma_{\a}\Psi(u)\;\;\;\mbox{ for }u\in V_{\ell K}.
\end{eqnarray}
Then using (\ref{e4.67})-(\ref{e4.70}) we get
\begin{eqnarray}
\Psi(Y(e^{\a},z)u)=Y(\Psi(e^{\a}),z)\Psi(u)
\;\;\;\;\mbox{ for }\a\in \ell K,\; u\in V_{\ell K}.
\end{eqnarray}
This proves the first assertion.
The second assertion can be proved similarly.
$\;\;\;\Box$

It is a well known fact (cf. [FZ]) that for $u,v\in V$,
$[Y(u,z_{1}),Y(v,z_{2})]=0$ if and only if $Y(u,z)v\in V[[z]]$,
or what is equivalent to, $u_{i}v=0$ for $i\in \N$.

\bp{pcoset1}
The subspace $\Omega_{V}^{0}$ is a vertex operator subalgebra of $V$
and
\begin{eqnarray}\label{e4.73}
V_{\ell K}=\{u\in V\;|\; [Y(u,z_{1}),Y(v,z_{2})]=0\;\;\;
\mbox{ for }v\in \Omega_{V}^{0}\}.
\end{eqnarray}
\ep

\pf It follows immediately from Theorem \ref{tcosetalgebra}
that $\Omega_{V}^{0}$ is a vertex operator subalgebra of $V$. 
(It also follows from a result of [FZ].) 

Set (cf. [DM])
\begin{eqnarray}
(\Omega_{V}^{0})^{c}=\{u\in V\;|\; [Y(u,z_{1}),Y(v,z_{2})]=0\;\;\;
\mbox{ for }v\in \Omega_{V}^{0}\}.
\end{eqnarray}
Since 
$h(i)\Omega_{V}^{0}=0$ for $h\in {\bf h},\; i\ge 0$, we have
${\bf h}\subset (\Omega_{V}^{0})^{c}$.
Recall that $Y_{\Omega}(v,z)=Y(v,z)$ for $v\in \Omega_{V}^{0}$. 
For $\a\in \ell K$, from Proposition \ref{psummary}, $\sigma_{\a}$ is an 
$\Omega_{V}^{0}$-automorphism of $V$.
Then for $v\in \Omega_{V}^{0},\; \a\in \ell K$,
\begin{eqnarray}
Y(v,z)\sigma_{\a}({\bf 1})
=\sigma_{\a}Y(v,z){\bf 1}\in V[[z]].
\end{eqnarray}
Hence $\sigma_{\a}({\bf 1})\in (\Omega_{V}^{0})^{c}$.
Since ${\bf h}$ and $\sigma_{\a}({\bf 1})$ for $\a\in \ell K$ 
generate $V_{\ell K}$ as a vertex algebra, we have 
$V_{\ell K}\subset (\Omega_{V}^{0})^{c}$.

Because $(\Omega_{V}^{0})^{c}$ is a completely reducible 
$\hat{\bf h}$-module and ${\bf h}\subset (\Omega_{V}^{0})^{c}$,
to prove $(\Omega_{V}^{0})^{c}\subset V_{\ell K}$ it suffices to prove
$$(\Omega_{V}^{0})^{c}\cap \Omega_{V}^{\b}\subset V_{\ell K}
\;\;\;\mbox{ for }\b\in L.$$
Let $u\in (\Omega_{V}^{0})^{c}\cap \Omega_{V}^{\b},\;\b\in L$.
Then $(\omega_{\Omega})_{i}u=0$ for $i\ge 0$.
In particular, $L_{\Omega}(-1)u=(\omega_{\Omega})_{0}u=0$.
It follows from Remark \ref{rvacuumlike} that
 $Y_{\Omega}(u,z)\in \End \Omega_{V}$ and
\begin{eqnarray}
f=Y_{\Omega}(u,z)(-1)^{{1\over\ell}\b (0)}
\in \Hom_{\Omega_{V}}^{{1\over\ell}\b}(\Omega_{V},\Omega_{V}).
\end{eqnarray}
Since by assumption $V$ is an irreducible $V$-module,
the linear map $\tilde{f}$ constructed in 
Proposition \ref{pisomorphismrelation}
must be a $V$-isomorphism from $V^{({1\over\ell}\b)}$ to $V$.
Thus ${1\over\ell}\b\in K$. By Schur lemma, 
$\tilde{f}\in \C\psi_{{1\over\ell}\b}$.
Consequently, 
$u=\tilde{f}({\bf 1})\in \C \sigma_{\b}({\bf 1})\in V_{\ell K}$.
The proof is complete. $\;\;\;\;\Box$

\br{rdm}
{\em Proposition \ref{pcoset1} states that $V_{\ell K}=\Omega_{V}^{c}$.
{}From this we have
\begin{eqnarray}
V_{\ell K}^{c}=(\Omega_{V}^{c})^{c}.
\end{eqnarray}
Since
\begin{eqnarray}
(V_{\ell K})^{c}\subset M(\ell)^{c}=\Omega_{V}^{0}\subset (\Omega_{V}^{c})^{c},
\end{eqnarray}
we have
\begin{eqnarray}
V_{\ell K}^{c}=\Omega_{V}.
\end{eqnarray}
Thus, in terms of a notion defined in [DM], $V_{\ell K}$ and 
$\Omega_{V}^{0}$ form a dual pair. }
\er

Let $A$ be a sublattice of $L_{0}\;(=\ell K\cap 2\ell L^{o}\subset \ell K)$ 
of the same rank. Then
\begin{eqnarray}
{1\over\ell}\<\a,\lambda\>\in \Z\;\;\;\mbox{ for }\a\in A,\; \lambda\in K^{o}.
\end{eqnarray}
Let $\{\b_{1},\dots, \b_{k}\}$ be a complete set of
representatives of cosets of $A$ in $K^{o}$.
Then
\begin{eqnarray}
V_{K^{o}}=\coprod_{i=1}^{k}V_{A+\b_{i}},
\end{eqnarray}
where $V_{A+\b_{i}}$ and $V_{A+\b_{j}}$ for $i\ne j$ are 
nonisomorphic irreducible $V_{A}$-modules.
For any weak $V$-module $W\in {\cal{C}}$, since $P(W)\subset K^{o}$
(Proposition \ref{pvacuumlike}), we have
\begin{eqnarray}
& &W=\coprod_{i=1}^{n}\Hom_{V_{A}}(V_{A+\b_{i}},W)\otimes V_{A+\b_{i}},\\
& &\Hom_{V_{A}}(V_{K^{o}},W)=\coprod_{i=1}^{n}\Hom_{V_{A}}(V_{A+\b_{i}},W).
\end{eqnarray}
Then $\Hom_{V_{A}}(V_{K^{o}},W)$ 
can be considered as the space of ``highest weight vectors'' of $V_{A}$ in $W$.
(Since there is no appropriate triangulated Lie algebra 
associated with $V_{A}$, there is a technical difficulty to define
the notion of a highest weight vector of $V_{A}$.)
In the following we shall identify $\Hom_{V_{A}}(V_{K^{o}},W)$ with
$\Omega_{W}^{A}$.

Let $W\in {\cal{C}}$. Since $V_{K^{o}}$ is a finitely generated 
$V_{A}$-module, by Proposition 4.8 of [Li7],
$\Hom_{V_{A}}(V_{K^{o}},W)$ 
is a natural module for $\Omega_{V}^{0}$ as a vertex algebra, where
\begin{eqnarray}
(Y(v,z)g)(u)=Y(v,z)g(u)\;\;\;\;\mbox{ for }v\in \Omega_{V}^{0},\;
g\in \Hom_{V_{A}}(V_{K^{o}},W),\; u\in V_{K^{o}}.
\end{eqnarray}

\bp{pnonabeliancoset}
Let $\{\b_{1},\dots, \b_{k}\}$ be a complete set of
representatives of cosets of $A$ in $K^{o}$.
For $W\in {\cal{C}}$, we define a linear map
\begin{eqnarray}
T: & &\Hom_{V_{A}}(V_{K^{o}},W)\rightarrow \Omega_{W}^{A}\nonumber\\
& &g\mapsto \sum_{i=1}^{k}g(e^{\b_{i}})+I_{W}^{A}.
\end{eqnarray}
Then $T$ is an $\Omega_{V}^{0}$-module isomorphism.
%Furthermore, if $A=2\ell K$,
%$T$ does not depend on the choices of representatives $\b_{i}$. 
\ep

\pf Since $V_{A}$, identified as a subalgebra of $V$, commutes with 
$\Omega_{V}^{0}$ (Proposition \ref{pcoset1}), 
it is easy to see that for any $u\in V_{K^{o}}$, 
the evaluation map $g\mapsto g(u)$ from $\Hom_{V_{A}}(V_{K^{o}},W)$
to $W$ is an $\Omega_{V}^{0}$-homomorphism.
Furthermore, if $u\in \Omega_{V_{K^{o}}}$, then the evaluation map
ranges in $\Omega_{W}$. Then it follows that $T$ is an 
$\Omega_{V}^{0}$-homomorphism. Now we need to show that $T$ is a
linear isomorphism.

Since $P(W)\subset K^{o}=\cup_{i=1}^{n}(\b_{i}+A)$ 
(Proposition \ref{pvacuumlike}), for convenience we may assume
\begin{eqnarray}
P(W)=\cup_{i=1}^{r}(\b_{i}+A).
\end{eqnarray}
Then the restriction of
the covering map $\xi$ (from $\Omega_{W}$ to $\Omega_{W}^{A}$)
\begin{eqnarray}
\xi:\hspace{1cm}
 \coprod_{i=1}^{r}\Omega_{W}^{\b_{i}}\rightarrow \Omega_{W}^{A}
\end{eqnarray}
is an $\Omega_{V}^{0}$-isomorphism.
For $g\in \Hom_{V_{A}}(V_{K^{o}},W)$, we have
\begin{eqnarray}
g(e^{\b_{i}})\in \Omega_{W}^{\b_{i}}\;\;\;\mbox{ for }i=1,\dots,n.
\end{eqnarray}
For $i>r$, since $\b_{i}\notin P(W)$, we must have
$g(e^{\b_{i}})=0$.
Then $T(g)=0$ if and only if $g(e^{\b_{i}})=0$ for $i=1,\dots,n$, 
which is equivalent to $g=0$ because
$e^{\b_{i}}$ $(i=1,\dots,n)$ generate $V_{K^{o}}$ as a $V_{A}$-module.
Then $T$ is injective.

On the other hand, let $0\ne w\in \Omega_{W}^{\b_{i}}$.
By Proposition \ref{platticesubalgebra},
 $w$ generates a $V_{A}$-submodule of $W$, which is isomorphic to
$V_{A+\b_{i}}$. Then there is $g\in \Hom_{V_{A}}(V_{K^{o}},W)$ such that
$T(g)=w+\Omega_{W}^{A}$. Then $T$ is onto. Therefore, $T$ is an 
$\Omega_{V}^{0}$-isomorphism.$\;\;\;\;\Box$
 
%Let $A=2\ell K$. 
%For $g\in \Hom_{V_{A}}(V_{K^{o}},W),\; \a\in \ell K,\; \b\in K^{o}$, 
%we have
%\begin{eqnarray}
%g(e^{2\a+\b})=g(e^{2\a}\cdot e^{\b})=g(Y_{\Omega}(e^{2\a},z)e^{\b})
%=Y_{\Omega}(\sigma_{2\a}({\bf 1}),z)g(e^{\b})=\sigma_{2\a}^{W}g(e^{\b}),
%\end{eqnarray}
%noting that ${1\over\ell}\<2\a,\b\>\in 2\Z$. Thus
%\begin{eqnarray}
%g(e^{2\a+\b})+I_{W}^{A}=g(e^{\b})+I_{W}^{A}.
%\end{eqnarray}
%It follows immediately that $T$ does not depend on the choices 
%of representatives $\b_{i}$.
%$\;\;\;\;\Box$

\br{rnonabelian}
{\em From Proposition \ref{pnonabeliancoset},
$\Omega_{W}^{A}$ can be identified as the ``vacuum space''
$\Hom_{V_{A}}(V_{K^{o}},W)$ of the subalgebra $V_{A}$ in $W$. 
In view of this,
the construction of $\Omega_{V}^{A}$ 
is essentially a non-abelian coset construction.}
\er

\section{Examples}
In this section we shall apply the results of Sections 3-5 
for $V=L(\ell,0)$, the vertex operator algebra associated to
an affine Lie algebra $\hat{\fg}$ with level $\ell$.

Let $\fg$ be a finite-dimensional simple Lie algebra, ${\bf h}$ be a
Cartan subalgebra and $\Phi$ be the set of roots. Let $\<\cdot,\cdot\>$ be 
the normalized killing form such that $\<\a,\a\>=2$ for a long root $\a$
where ${\bf h}^{*}$ is identified with ${\bf h}$ through
$\<\cdot,\cdot\>$. 
Let $\Q$ and $\Q^{\vee}$ be the root lattice and the co-root lattice,
respectively. Then $(\Q^{\vee})^{o}$ denoted by $P$ is the weight lattice
and $\Q^{o}=P^{\vee}$ is the co-weight lattice.

Let $\hat{\fg}$ be the affine Lie algebra:
\begin{eqnarray}
\hat{\fg}=\fg\otimes \C[t,t^{-1}]\oplus\C c,
\end{eqnarray}
where
\begin{eqnarray}	
&&[\hat{\fg},c]=0,\\
& &[a\otimes t^{m},b\otimes t^{n}]
=[a,b]\otimes t^{m+n}+m\<a,b\>\delta_{m+n,0}c
\;\;\;\;\mbox{ for }a,b\in \fg,\; m,n\in \Z.
\label{eaffinerelation}
\end{eqnarray}
As usual we also
use $a(n)$ for $a\otimes t^{n}$. For $n\in \Z$, we denote
\begin{eqnarray}
\fg (n)=\{ a(n)\;|\; a\in \fg\}.
\end{eqnarray}
For $a\in \fg$, define the generating series
\begin{eqnarray}
a(z)=\sum_{n\in \Z}(a\otimes t^{n})z^{-n-1}\in \hat{\fg}[[z,z^{-1}]].
\end{eqnarray}

\bd{dcategory1}
{\em For $\ell \in \C$, we define a category ${\cal{R}}_{\ell}$
of level-$\ell$ restricted $\hat{\fg}$-modules $W$ 
(cf. [K]) in the sense that $c$ acts as $\ell$ and 
for every $w\in W$, $\fg (n)w=0$ 
for $n$ sufficiently large.}
\ed

{\bf Assumption:}
Throughout this section we assume $\ell\in \C-\{0,-h^{\vee}\}$, 
where $h^{\vee}$ is the dual coxeter number of $\fg$.

Consider the generalized Verma $\hat{\fg}$-module
\begin{eqnarray}
M(\ell,0)=U(\hat{\fg})\otimes_{U(\fg\otimes \C[t]+\C c)}\C,
\end{eqnarray}
$\fg\otimes \C[t]$ 
acting trivially on $\C_{\ell}$ and $c$ acting as $\ell$. 
Denote by ${\bf 1}$
the highest weight vector $1\otimes 1$ of $M(\ell,0)$.
Let $L(\ell,0)$ be the (unique) irreducible quotient module 
of $M(\ell,0)$. We abuse ${\bf 1}$ for the image of ${\bf 1}$ 
in $L(\ell,0)$.
Identify $\fg$ as a subspace of $M(\ell,0)$ 
and $L(\ell,0)$ through the map $a\mapsto a(-1){\bf 1}$.
%For any $a\in \fg$, there is $b\in \fg$ such that $\<a,b\>\ne 0$.
%Then
%\begin{eqnarray}
%b(1)a(-1){\bf 1}=[b,a](0){\bf 1}+\ell \<b,a\>{\bf 1}=\ell \<b,a\>{\bf 1}\ne 0.
%\end{eqnarray}

Set
\begin{eqnarray}
\omega={1\over 2(\ell+h^{\vee})}\sum_{i=1}^{d}a_{i}(-1)a_{i}(-1){\bf 1}
\in M(\ell,0),
\end{eqnarray}
where $\{a_{1},\dots, a_{d}\}$ is any orthonormal basis of $\fg$.
Then we have (see [DL2], [FF], [FZ], [Li2] and [MP2]):

\bp{pfzli}
There exists a unique vertex operator algebra structure 
on $M(\ell,0)$ and $L(\ell,0)$ such that $Y(a,z)=a(z)$ 
for $a\in \fg$ with ${\bf 1}$ and $\omega$
being the vacuum vector and the Virasoro vector, respectively.
The central charge is
\begin{eqnarray}
c_{\ell}=\frac{\ell \dim \fg}{\ell+h^{\vee}}.
\end{eqnarray}
Furthermore, the category of modules for $M(\ell,0)$ viewed
as a vertex algebra is naturally isomorphic to the category 
${\cal{R}}_{\ell}$.
\ep

The vertex operator $Y(\omega,z)$ turns out to generate
the ``Segal-Sugawara'' realization of the Virasoro algebra 
(see [Su], [BH], [Se]).

Note that $\hat{\bf h}={\bf h}\otimes \C[t,t^{-1}]+\C c$
is a subalgebra of $\hat{\fg}$ and the subalgebra
$\hat{\bf h}_{\Z}=\hat{\bf h}^{+}+\hat{\bf h}^{-}+\C c$
is a Heisenberg algebra. We have
\begin{eqnarray}
M(\ell)\subset M(\ell,0), \;\; L(\ell,0).
\end{eqnarray}

\bd{dcategory2}
{\em We define a category ${\cal{C}}_{\ell}$
of level-$\ell$ restricted $\hat{\fg}$-modules $W$ such that
${\bf h}$ acts semisimply on $W$ and such that
$\dim U(\hat{\bf h}^{+})w<\infty$ for $w\in W$.}
\ed

Note that in view of Proposition \ref{pfzli},
the category ${\cal{C}}_{\ell}$ is exactly the category
${\cal{C}}$ defined in Section 3 for $V=M(\ell,0)$.

Now, we take $V=M(\ell,0)$ (or $L(\ell,0)$). 
We have
\begin{eqnarray}
\fg_{\a}\subset \Omega_{M(\ell,0)}^{\a}
\;\;\;\mbox{ for }\a\in \Phi
\end{eqnarray}
as
\begin{eqnarray}
h(n)a=h(n)a(-1){\bf 1}=\a(h)a(n-1){\bf 1}+a(-1)h(n){\bf 1}
=\delta_{n,0}\a(h)a
\end{eqnarray}
for $h\in {\bf h},\; n\ge 0,\; a\in \fg_{\a}.$
Since 
$M(\ell,0)$ is generated from ${\bf 1}$ by $\hat{\fg}$,
$M(\ell,0)$ is $\Q$-graded.
Then $L$ defined in Section 3 is exactly
the root lattice $\Q$ (identified as a subset of ${\bf h}$).

\bt{tcosetalgebraaffine}
There exists a unique generalized vertex algebra structure 
$Y_{\Omega}$ on $\Omega_{M(\ell,0)}$ and $\Omega_{L(\ell,0)}$ 
with $G=\Q, \; c(\cdot,\cdot)=1$
and $(\cdot,\cdot)$ defined by
\begin{eqnarray}\label{eapply1}
(\a,\b)={1\over\ell}\<\a,\b\>+2\Z
\;\;\;\mbox{ for }\a,\b\in \Q
\end{eqnarray}
such that $Y_{\Omega}({\bf 1},z)=1$ and
\begin{eqnarray}\label{eapply2}
Y_{\Omega}(a,z)
=E^{-}({1\over\ell}\a,z)a(z)E^{+}({1\over\ell}\a,z)
z^{-{1\over\ell}\a (0)}\;\;\;\mbox{ for }
a\in \fg_{\a},\; \a\in \Phi.
\end{eqnarray}
Furthermore, $\Omega_{M(\ell,0)}$ and $\Omega_{L(\ell,0)}$ are 
generated by $\fg_{\a}$ ($\a\in \Phi$), $\fg_{\a}$ is of weight 
$1-{1\over 2\ell}{\<\a,\a\>}$ and
the following relations hold
for $u\in \fg_{\a},\; v\in \fg_{\b},\;\a,\b\in \Phi$:
$$(z_{1}-z_{2})^{{1\over\ell}\<\a,\b\>}
Y_{\Omega}(u,z_{1})Y_{\Omega}(v,z_{2})
-(z_{2}-z_{1})^{{1\over\ell}\<\a,\b\>}
Y_{\Omega}(v,z_{2})Y_{\Omega}(u,z_{1})=$$
\begin{eqnarray}\label{ezalgebrarelations}
=\left\{\begin{array}{c}
z_{1}^{-1}\delta(z_{2}/z_{1})Y_{\Omega}([u,v],z_{2})
(z_{2}/z_{1})^{{1\over\ell}\a(0)}
\hspace{2cm} \mbox{ if }\a+\b\ne 0,\\
\ell \<u,v\> {\partial\over\partial z_{2}}
\left(z_{1}^{-1}\delta(z_{2}/z_{1})
(z_{2}/z_{1})^{{1\over \ell}\a(0)}\right)\hspace{2cm}\mbox{ if }\a+\b=0.
\end{array}\right.
\end{eqnarray}
\et

\pf It follows from Theorem \ref{tcosetalgebra} that
there exists a generalized vertex algebra structure 
$Y_{\Omega}$ on $\Omega_{M(\ell,0)}$ and $\Omega_{L(\ell,0)}$
such that (\ref{eapply1}) and (\ref{eapply2}) hold.

Since $\fg$ generates $M(\ell,0)$ as a vertex algebra and
$\fg/{\bf h}=\sum_{\a\in \Phi}\fg_{\a}$ is an ${\bf h}$-submodule
of $\Omega_{V}$, by Proposition \ref{pgenerators},
$\fg/{\bf h}$ generates $\Omega_{M(\ell,0)}$ as a 
generalized vertex algebra.
Then the uniqueness follows immediately.

For $a\in \fg_{\a},\; \a\in \Phi$, using Lemma \ref{lvirasoro} we get
\begin{eqnarray}
L_{\Omega}(0)a=(L(0)-L_{\bf h}(0))a
=\left(1-{1\over 2\ell}{\<\a,\a\>}\right)a.
\end{eqnarray}
That is, the $L_{\Omega}(0)$-weight of $a$ is $1-{1\over 2\ell}{\<\a,\a\>}$.

{}From the generalized Jacobi identity we get
\begin{eqnarray}\label{egcommformula}
& &(z_{1}-z_{2})^{{1\over\ell}\<\a,\b\>}
Y_{\Omega}(u,z_{1})Y_{\Omega}(v,z_{2})
-(z_{2}-z_{1})^{{1\over\ell}\<\a,\b\>}
Y_{\Omega}(v,z_{2})Y_{\Omega}(u,z_{1})\nonumber\\
&=&\Res_{z_{0}}z_{0}^{{1\over\ell}\<\a,\b\>}
z_{1}^{-1}\delta\left(\frac{z_{2}+z_{0}}{z_{1}}\right)
Y_{\Omega}(Y_{\Omega}(u,z_{0})v,z_{2})
\left(\frac{z_{2}+z_{0}}{z_{1}}\right)^{{1\over\ell}\a(0)}
\nonumber\\
&=&\Res_{z_{0}}z_{0}^{{1\over\ell}\<\a,\b\>}
Y_{\Omega}(Y_{\Omega}(u,z_{0})v,z_{2})
e^{z_{0}{\partial\over\partial z_{2}}}
\left(z_{1}^{-1}\delta(z_{2}/z_{1})(z_{2}/z_{1})^{{1\over\ell}\a(0)}\right).
\end{eqnarray}
Note that from the defining relations of $\hat{\fg}$ we have
\begin{eqnarray}
Y(u,z)v
=\ell \<u,v\>z^{-2}+[u,v]z^{-1}+\mbox{regular terms.}
\end{eqnarray}
Then
\begin{eqnarray}\label{e6.18}
z^{{1\over\ell}\<\a,\b\>}Y_{\Omega}(u,z)v&=&
E^{-}({1\over\ell}\a,z)Y(u,z)v\nonumber\\
&=&\ell \<u,v\>z^{-2}+([u,v]-\<u,v\> \a)z^{-1}+\mbox{regular terms.}
\end{eqnarray}
Note also that if $\a+\b\ne 0$, we have $\<u,v\>=0$
and if $\a+\b=0$, we have $[u,v]-\<u,v\>\a=0$. 
Then (\ref{ezalgebrarelations}) immediately follows from 
(\ref{egcommformula}) and (\ref{e6.18}). $\;\;\;\;\Box$

\br{ronmodule}
{\em Note that for $W\in {\cal{C}}_{\ell}$, by Proposition \ref{pfzli}
$W$ is an $M(\ell,0)$-module. Hence, by Theorem \ref{tcosetmodule},
$W$ is an $\Omega_{M(\ell,0)}$-module.
Then the same proof shows that (\ref{ezalgebrarelations}) holds on
$W$.}
\er

\br{rlwlp1} {\em Let $W\in {\cal{C}}_{\ell}$.
Recall that for $a\in \fg_{\a},\; \a\in \Phi$,
\begin{eqnarray}
Z(a,z)=Y_{\Omega}(a,z)z^{{1\over \ell}\a(0)}.
\end{eqnarray}
Multiplying (\ref{ezalgebrarelations}) from right by 
$z_{1}^{{1\over\ell}\a(0)}z_{2}^{{1\over\ell}\b(0)}$ and using
(\ref{ehmyomega}) we obtain
$$(1-z_{2}/z_{1})^{\<\a,\b\>/\ell}Z(u,z_{1})Z(v,z_{2})
-(1-z_{1}/z_{2})^{\<\a,\b\>/\ell}Z(v,z_{2})Z(u,z_{1})=$$
\begin{equation}\label{ezalgebra}
=\left\{
\begin{array}{c}
z_{1}^{-1}\delta(z_{2}/z_{1})Z([u,v],z_{2})
\hspace{4cm}\mbox{ if }\a+\b\ne 0\\
\<u,v\>\left(z_{1}^{-1}\delta(z_{2}/z_{1})\a z_{2}^{-1}
+ \ell {\partial \over\partial z_{2}}z_{1}^{-1}\delta(z_{2}/z_{1})\right)
\;\;\;\;\mbox{ if }\a+\b=0.
\end{array}\right.
\end{equation}
This is a well known result due to [LP1].}
\er

Let $\b\in \Q^{o}$ and let $W$ be a restricted $\hat{\fg}$-module 
$W$ of level $\ell$. From Proposition \ref{pfzli}, $W$
is a natural $M(\ell,0)$-module. Then we have an 
$M(\ell,0)$-module $W^{(\b)}$. From Proposition \ref{pfzli} again,
$W^{(\b)}$ is a $\hat{\fg}$-module of level $\ell$.

Define a linear endomorphism $\theta_{\b}$ of $\hat{\fg}$ by
\begin{eqnarray}
& &\theta_{\b}(c)=c,\\
& &\theta_{\b}(a(z))=z^{\<\a,\b\>}a(z)\;\;\;\mbox{ for }a\in \fg_{\a},
\; \a\in \Phi,\\
& &\theta_{\b}(h(z))=h(z)+c \<\b,h\> z^{-1}
\;\;\;\mbox{ for }h\in {\bf h}.
\end{eqnarray}
It is easy to prove (cf. [FS], [Li4]) that $\theta_{\b}$ is an automorphism
of $\hat{\fg}$. 
Then for any representation $\rho$ of $\hat{\fg}$ on 
$W$, $\rho \theta_{\b}$ is also a representation of $\hat{\fg}$ on 
$W$. It follows from [DLM1] and [Li4] (Proposition 3.5) that 
the representation $\rho \theta_{\b}$ of $\hat{\fg}$ on 
$W$ is canonically equivalent to the representation
of $\hat{\fg}$ on  $W^{(\b)}$. As an immediate corollary
of Proposition \ref{pirrational} we have:

\bp{pequivalence2}
Let $\ell$ be non-rational. Then the map $W\mapsto \Omega_{W}$ gives rise to
a one-to-one correspondence between the set of $\Q^{o}$-orbits
of the set of equivalence classes of irreducible 
$\hat{\fg}$-modules in ${\cal{C}}_{\ell}$ and 
the set of equivalence classes of irreducible
$\Omega_{L(\ell,0)}$-modules. $\;\;\;\;\Box$
\ep

It was known that $M(\ell,0)=L(\ell,0)$ is simple. Then we immediately
have:

\bc{csimple}
If $\ell$ is not rational, the generalized vertex algebra 
$\Omega_{M(\ell,0)}$ $(=\Omega_{L(\ell,0)})$ is simple 
in the sense that the adjoint module is irreducible. $\;\;\;\;\Box$
\ec

Now, let $\ell$ be a positive integer. Set
\begin{eqnarray}
P_{\ell}=\{\lambda\in P_{+}\;|\; \<\lambda, \theta\>\le \ell\},
\end{eqnarray}
where $P_{+}$ is the set of all dominant integral weights of $\fg$ 
(cf. [K]).
Then for $\lambda\in P_{\ell}$,
$L(\ell,\lambda)$ is an integral module, or a standard module.
A known fact about the vertex operator algebra $L(\ell,0)$ 
(cf. [DL2], [DLM2], [FZ], [Li2]) is
that the category of all weak $L(\ell,0)$-modules is semisimple
and that irreducible weak $L(\ell,0)$-modules are exactly
standard or highest weight irreducible integrable $\hat{\fg}$-modules 
of level $\ell$.

In view of Theorem \ref{tequivalence} we immediately have:

\bp{pmain}
Let $\ell$ be a positive integer. Then
every $\Omega_{L(\ell,0)}$-${\bf h}$-module is completely reducible and
the set of $\Omega_{L(\ell,\lambda)}$ for $\lambda\in P_{\ell}$
is a complete set of representatives of equivalence classes of irreducible
$\Omega_{L(\ell,0)}$-${\bf h}$-modules. $\;\;\;\;\Box$
\ep

It was proved in [Li6] that 
\begin{eqnarray}
K=\{\a\in \Q^{o}\;|\; L(\ell,0)^{(\a)}\simeq L(\ell,0)\}=\Q^{\vee}.
\end{eqnarray}
Then we have
\begin{eqnarray}
& &K^{0}=(\Q^{\vee})^{o}=P,\\
& &G=K^{o}/2\ell K=P/2\ell \Q^{\vee}
\end{eqnarray}
(cf. [DL2]).
Then any $\Omega_{L(\ell,0)}^{2\ell \Q^{\vee}}$-module is
$G$-graded. 
In view of Theorem \ref{tmain} we immediately have:

\bt{tmainaffine}
Let $\ell$ be a positive integer. Then $\Omega_{L(\ell,0)}^{2\ell\Q^{\vee}}$
is a simple generalized vertex operator algebra with 
$G=P/2\ell \Q^{\vee}$. Furthermore, the category of 
$G$-graded modules for $\Omega_{L(\ell,0)}^{2\ell\Q^{\vee}}$
as a generalized vertex algebra is semisimple and
the map $W\mapsto \Omega_{W}^{2\ell\Q^{\vee}}$ gives
rise to a one-to-one correspondence between the set
of $\Q^{o}$-orbits of the set of equivalence classes of
irreducible $L(\ell,0)$-modules and the set of equivalence classes of
$G$-graded irreducible $\Omega_{L(\ell,0)}^{2\ell\Q^{\vee}}$-modules.
$\;\;\;\;\Box$
\et

\end{document}